\documentclass[10pt]{article}
\usepackage{amsmath}
\usepackage{latexsym}
\usepackage{amssymb}
\usepackage{amsfonts}
\usepackage{amsthm}
\title{CLAUSIUS/COSSERAT/MAXWELL/WEYL EQUATIONS: \\  THE VIRIAL THEOREM REVISITED}
\author{J.-F. Pommaret \\ CERMICS, Ecole des Ponts ParisTech,\\ 6/8 Av. Blaise Pascal, 77455 Marne-la-Vall\'ee Cedex 02, France \\
E-mail: jean-francois.pommaret@wanadoo.fr, pommaret@cermics.enpc.fr \\
URL: http://cermics.enpc.fr/$\sim$pommaret/home.html }
\date{  }
\begin{document}
\maketitle

\noindent
{\bf ABSTRACT} \\

In 1870, R. Clausius found the {\it virial theorem} which amounts to introduce the trace of the stress tensor when studying the foundations of thermodynamics, as a way to relate the absolute temperature of an ideal gas to the mean kinetic energy of its molecules.   \\
\hspace*{5mm}In 1901, H. Poincar\'{e} introduced a {\it duality principle} in analytical mechanics in order to study lagrangians invariant under the action of a Lie group of transformations. In 1909, the brothers E. and F. Cosserat discovered {\it another approach for studying the same problem though using quite different equations}. In 1916, H. Weyl considered again {\it the same problem} for the conformal group of transformations, obtaining {\it at the same time} the Maxwell equations and an additional specific equation also involving the trace of the impulsion-energy tensor. Finally, having in mind the space-time formulation of electromagnetism and the Maurer-Cartan equations for Lie groups, gauge theory has been created by C.N. Yang and R.L. Mills in 1954 as a way to introduce in physics the differential geometric methods available at that time, {\it independently of any group action}, contrary to all the previous approaches.\\
\hspace*{5mm}The main purpose of this paper is to revisit the mathematical foundations of thermodynamics and gauge theory by using new differential geometric methods coming from the formal theory of systems of partial differential equations and Lie pseudogroups, mostly developped by D.C Spencer and coworkers around 1970. In particular, we {\it justify} and {\it extend} the virial theorem, showing that the Clausius/Cosserat/Maxwell/Weyl equations are nothing else but the formal adjoint of the {\it Spencer operator} appearing in the canonical {\it Spencer sequence} for the conformal group of space-time and are thus {\it totally dependent on the group action}. The duality principle also appeals to the formal adjoint of a linear differential operator used in differential geometry and to the extension modules used in homological algebra.\\

\hspace*{10mm}


\vspace*{40mm}

\noindent
{\bf 1) INTRODUCTION}  \\

There are many ways to define the concept of " {\it temperature} " in thermodynamics or thermostatics. A very useful one depends on the properties of the so-called {\it ideal gases} under a pressure not exceeding the atmospheric pressure, summarized by the following three experimental laws:  \\

\noindent
1) {\it The Boyle-Mariotte law} : Discovered by Boyle in England (1662), it has been rediscovered by Mariotte in France (1676). For a given mass of a gas at a constant temperature, say the molar mass $M$, the product of the pressure $P$ by the volume $V$ occupied by this gas is (approximatively) constant.  \\
\noindent
2) {\it The Gay-Lussac-Charles law} : Established around 1800 after the works of Gay-Lussac, Charles and Dalton, it says that, under the conditions of the preceding law, the product $PV$ does not depend on the gas but only on the temperature.  \\
\noindent
3) {\it The Avogadro-Amp\`{e}re law} : Stated around 1810 by Avogadro, it says that the product $PV$ for a given gas at a given temperature is proportional to the number of moles of the gas or to the number of molecules as a mole is made by $N$ molecules where $N=6,0225.10^{23}$ is the {\it Avogadro number}.  \\

As a byproduct, an ideal gas is such that $PV=nRT$ where $n$ is the number of moles and $k=R/N$ is the {\it Boltzmann constant} while $T$ is the {\it ideal gas scale of temperature}, also called {\it absolute temperature}. \\

The {\it first principle of thermostatics} says that the exchange of work $\delta W= - PdV$ plus the exchange of heat $\delta Q=CdT+L dV$ of the system with its surroundings is a {\it total differential}, that is there exists a function $U=U(T,V)$ called {\it internal energy}, such that $dU=\delta W + \delta Q$. Accordingly, the properties of ideal gases are complemented by another experimental law.  \\

\noindent
4) {\it The Joule law} : Stated by Joule in 1845 who introduced on this occasion the concept of internal energy, it says  that the internal energy $U$ of an ideal gas only depends on the temperature, that is $U=U(T)$. \\

This law has been checked by means of various expansion experiments realized by Gay-Lussac (1806), Joule (1845) and Hirn (1865). The idea is to consider an adiabatic cylinder separated in the middle by a wall with a tap which is suddenly opened or by a glass window which is suddenly broken. One part is filled with a gas at temperature $T$ while the other part is empty. At the end of the experiment, which is therefore done without any exchange of heat or work with the surroundings, one checks that the final temperature of the expanded gas is again $T$. As the new volume is twice the initial volume, the law follows with quite a good precision (apart for helium discovered later on).\\

The {\it second principle of thermostatics} says that the $1$-form $\delta Q$ admits an {\it integrating factor} which is a function of the absolute temperature {\it only}, that is one can find a function $\theta=\theta (T)$ and a function $S=S(T,V)$ called {\it entropy} such that $\delta Q=\theta dS$ has the integrating factor $1/\theta$. More generally, if $\delta Q=\theta dS={\theta}'dS'$ for two arbitrary $\theta (V,T)$ and ${\theta}'(V,T)$, we get $S'=h(S)$ and thus $1/{\theta}'=(\partial h(S)/\partial S)(1/ \theta)$. \\

In the case of an ideal gas, $dU=CdT+(L-P)dV \Rightarrow \partial C/\partial V - \partial (L-P)/\partial T=0$ while $dS=(C/\theta)dT+(L/\theta)dV \Rightarrow \partial (C/\theta)/\partial V - \partial (L/\theta)/\partial T=0$. First of all, it follows from the Joule law that $L=P$ on one side and thus $C=C(T)$ on the other side. As a byproduct, $C/\theta$ only depends on $T$ and $P/\theta = RT/\theta V$ must only depend on $V$, that is $T/\theta =c=cst$ or $T=c \theta$, a result showing that the ideal gas scale of temperature $T$ can be used in place of $\theta$ by choosing $c=1$ or, equivalently, that {\it the absolute temperature is only defined up to a scaling factor}. It also follows that we may choose $U=U(S,V)$ with $dU=TdS-PdV$ and that the so-called {\it free energy} $F=U-TS$ introduced by Helmholtz is such that $dF=-PdV-SdT$, a result leading therefore to a function $F=F(T,V)$ allowing to {\it define} $S=-\partial F/\partial T$ and thus $U=F-T\partial F/\partial T$ as a way to bypass  the principles by means of a mechanical approach to thermodynamics along the {\it helmholtz analogy} that we now recall. Indeed, in the lagrangian approach to analytical mechanics that we shall see thereafter, if one has functions $q(t)$ of time, for example positions $x(t),y(t),z(t)$ of points in cartesian space, and a {\it lagrangian} $L(t,q,\dot{q})$ where $\dot{q}=dq/dt$, the variational calculus applied to $\int L(t,q,\dot{q})dt$ may produce the Euler-Lagrange (EL) equations $\frac{d}{dt}(\frac{\partial L}{\partial \dot{q}})-\frac{\partial L}{\partial q}=0$. Introducing the {\it hamiltonian} $H=\dot{q}\frac{\partial L}{\partial \dot{q}}-L$, we get $\frac{dH}{dt}=\frac{\partial L}{\partial t}$ and thus the conservation of energy $H=cst$ whenever $\frac{\partial L}{\partial t}=0$. Accordingly, if one could find a function $q(t)$ such that $T=\dot{q}$, then one could recover the previous formulas on the condition to choose $L=-F$ (See [22],[25] for more details). \\

The following three examples are among the best ones we have been able to find in order to understand why exhibiting an integrating factor may not at all be as simple as what is claimed in most textbooks.\\

\noindent
{\bf EXAMPLE 1.1}: ({\it Ideal Gas}) With volume $V$, absolute temperature $T$, pressure $P$, entropy $S$ and internal energy $U$ for one mole of a perfect gas, we obtain $dU=\delta W + \delta Q $ with $\delta W=-PdV$ and $ \delta Q=CdT+LdV$ where $C=C_V$ is the heat capacity at constant volume and $PV=RT$ for one mole. Replacing and writing that $dU$ and $dS=(1/T)\delta Q$ are closed $1$-forms, we obtain successively $\partial C/ \partial V-\partial (L-P)/ \partial T=0$ and $\partial (C/T)/ \partial V-\partial (L/T)/ \partial T=0$, that is to say $L=P$ {\it and} $C=C(T)$. We get therefore $\delta Q=CdT+PdV$ and thus $dU=CdT$. However, when $C$ is a constant as in the case of an ideal gas, looking for a general integrating factor of the form $A(V,T)/T$, the $1$-form $(CA/T)dT+(RA/V)dV$ must be closed and thus $(C/T)\partial A/\partial V-(R/V)\partial A/\partial T=0$, a result leading to $A=A(VT^{\alpha})$ where $\alpha = C/R$ with $R=C_P-C_V=(\gamma -1)C$ according to the Mayer's relation. Of course, we find the well known integrating factor $1/T$ leading to $S=Rlog(VT^{\alpha})$ and $F=CT(1 - log(T)) - RT log(V)$, but we could also use the other integrating factor $VT^{\alpha -1}$ leading to $S'=RVT^{\alpha}$ and get $S'=R\hspace{1mm}exp(S/R)$. If we look for an integrating factor depending only on $T$, we can only have $c/T$ whith an arbitrary non-zero constant $c$ used in order to fix the absolute temperature up to a change of scale.   \\                                                                                                                                                                                                                                                                                                                                                                                                                                                                                                                                                                                                                                                                                                                                                                                                                                                                                                                                                                                                                                                                                                                                                                                                                                                                                                                                                                                                                                                                                                                                                   

\noindent
{\bf EXAMPLE 1.2}:({\it Black Body}) Using the same notations, we have now $U=\alpha VT^4$ and $P=\frac{1}{3}\alpha T^4 \Rightarrow \delta Q=dU-\delta W=dU+PdV= 4 \alpha VT^3dT+\frac{4}{3}\alpha T^4dV$. Looking again for any integrating factor of the same form $A(V,T)/T$ as before, we should obtain $3V\partial A/\partial V-T\partial A/\partial T$ and thus $A=A(VT^3)$. Of course, $1/T$ is the standard integrating factor leading to 
$S=\frac{4}{3}\alpha VT^3$ and $F= - \frac{1}{3}\alpha V T^4$. However, $VT^3/T=VT^2$ is also an integrating factor with $VT^2\delta Q=4\alpha V^2T^5dT+\frac{4}{3}\alpha VT^6dV=d(\frac{2}{3}\alpha V^2T^6)=dS' $ and we get $S'=\frac{3}{8\alpha}S^2$.  \\                                                                                                                                        

\noindent
{\bf EXAMPLE 1.3}: (Compare to [13], p 117) Two {\it different} ideal gases, one mole each, with respective heat capacities $C'$ at constant volume $V'$ and $C"$ at constant volume $V"$ such that $C'\neq C"$ are contained in a cylinder and separated by an adiabatic movable piston. We shall prove that there cannot be any integrating factor for the exchange of heat $\delta Q=\delta Q' + \delta Q"$ of this system. Using the first law of thermodynamics as in the previous examples, we have $\delta Q'=C'dT'+P'dV', \delta Q"=C"dT"+P"dV"$. However, for a reversible transformation, the piston must be in mechanical equilibrium and thus $P=P'=P"$. Now, we have $PV'=RT', PV"=RT"$ and we obtain therefore $PdV'=RdT'-(RT'/P)dP, PdV"=RdT"-(RT"/P)dP$ for the system now described by the only three state variables $T',T", P$. Accordingly, we get the $1$-form $\alpha=\delta Q= (C'+R)dT'+(C"+R)dT"-(R/P)(T'+T")dP$. Taking the exterior derivative, we get $d\alpha=(R/P)dP\wedge (dT' + dT")$ and thus $\alpha \wedge d\alpha= (R/P)(C"-C' )    dT'\wedge dT" \wedge dP \neq 0$. Accordingly, integrating factors do not exist in general for systems which are not in thermal equilibrium. \\

It remains to relate this {\it macroscopic aspect} of thermodynamics that we have presented with its {\it microscopic aspect}, in particular with the kinetic theory of gases. For this, assuming the {\it molecular chaos}, namely that the gases are made by a juxtaposition of individual molecules of mass $m$ with $M=Nm$, we assume that, at a given time, the directions of the speeds have a random distribution in space, that the size of the molecules is small compared to their respective average distance and that the average speed is bigger when the temperature is higher. We also assume that there is no interactions apart very negligible attractive forces compared to the repulsive shock forces existing on vey small distances. As a byproduct, the pressure is produced by the only forces acting on the wall of a containing volume $V$ limited by a surface $S$ with outside normal vector $\vec{n}$ which are made by the molecules hitting the surface. We explain the way followed by Clausius.  \\

If $O$ is a fixed point inside $V$, for example the origin of a cartesian frame, and $M$ an arbitrary point (for a few lines), we have:  \\
\[ \frac{d(\stackrel{\longrightarrow}{OM})^2}{dt} = 2 \stackrel{\longrightarrow}{OM}.\frac{d\stackrel{\longrightarrow}{OM}}{dt}= 2\stackrel{\longrightarrow}  {OM}.\vec{v} \hspace{3mm} \Longrightarrow  \hspace{3mm}
\frac{d^2(\stackrel{\longrightarrow}{OM})^2}{dt^2} = 2( \frac{d\stackrel{\longrightarrow}{OM}}{dt})^2 + 2 \stackrel{\longrightarrow}{OM}.\frac{d^2\stackrel{\longrightarrow}{OM}}{dt^2}                 \]        
Multiplying by $\frac{m}{4}$, we recognize, in the right member, the kinetic energy of a molecule and the force $\vec{F} = m \frac{d^2\stackrel{\longrightarrow}{OM} }{dt^2}$ acting on this molecule at time $t$. Summing on all the molecules contained in $V$ while taking into account the fact that the sum $\sum m (\stackrel{\longrightarrow}{OM})^2$ is constant when the statistical equilibrium is achieved, we obtain the formula:  \\
\[    \sum \frac{1}{2}m(\vec{v})^2= - \sum \frac{1}{2}\stackrel{\longrightarrow}{OM}.{\vec{F} } \]
where the term on the right side is called {\it virial} of the gas. In the case of an ideal gas, the forces are annihilated two by two apart from the ones existing on $S$. However, the force produced by the pressure on a small part $dS$ of $S$ is known to be $d\vec{F}= - P\vec{n} dS$. Taking into account that $P$ is constant inside $V$ and on $S$, the total kinetic energy contained in $V$ is thus equal to the half of: \\
\[ {\int}_SP\stackrel{\longrightarrow}{OM}.\vec{n}dS= {\int}_Vdiv(P\stackrel{\longrightarrow}{OM})dV=P{\int}_Vdiv(\stackrel{\longrightarrow}{OM})dV=3PV\]
after using the Sokes formula because ${\partial}_1x^1+{\partial}_2x^2+{\partial}_3x^3=3$. Introducing the mean quadratic speed $u$ such that $\Sigma (\vec{v})^2=N u^2$ for a mole of gas with $N$ molecules and mass $M=Nm$, we obtain therefore $PV=\frac{1}{3}\Sigma m(\vec{v})^2= \frac{1}{3} Nmu^2=\frac{1}{3} Mu^2$ and recover the experimental law found by Boyle and Mariotte. As a byproduct, we find $\frac{1}{2}mu^2= \frac{3}{2}kT$ for the mean translational kinetic energy of a molecule.  \\

In order to start establishing a link between the virial theorem that we have exhibited and group theory, let us recall that the stress equation of continuum mechanics is ${\partial}_r{\sigma}^{ri}=f^i$ when the ambient space is ${\mathbb{R}}^3$ with cartesian coordinates and that the stress in a liquid or a gas is the $3\times 3$ diagonal matrix with diagonal terms equal to $-P$. Using the only infinitesimal generator $\theta = x^i{\partial}_i$ of the dilatation group while raising or lowering the indices by means of the euclideam metric of ${\mathbb{R}}^3$, we obtain (Compare to (29)+(30) in [9]):  \\
\[   x^i{\partial}_r{\sigma}^r_i= {\partial}_r(x^i{\sigma}^r_i)- {\sigma}^r_r=x^if_i   \]
as a way to exhibit the trace of the stress tensor $\sigma$ but, of course, it remains to justify this purely technical computation by means of group theoretical arguments.\\

We conclude this paragraph with a few comments on the so-called {\it axiomatic thermodynamics} initiated by P. Duhem (1861-1916) around 1892-1894 in ([8])  and then by C. CarathŽodory in 1909 ([5]) (See the pedagogical review made by M. Born in 1921 [4]).\\
A first comment concerns the use of differential forms (See a forthcoming paragraph for definitions), introduced by E. Cartan in 1899 but only used in physics and particularly in thermodynamics after decades. If $\alpha = \delta Q$ and $\beta = \delta W$ are respectively the exchange of heat and work of the system with its surroundings, one must never forget that any finite heat $Q$ and work $W$ obtained by integration is counted positively if it is {\it provided to} the system (One of the best references we know is by far [13]). In this framework, the {\it first principle} amounts to $\alpha + \beta = dU$ where $U$ is the {\it internal energy} or, equivalently, $d(\alpha +\beta)=0$. As for the {\it second principle} amounting to the existence of an " {\it integrating factor} " for $\alpha$, that is the possibility to write $\delta Q=TdS$, it is well known that it is equivalent to the existence of a $1$-form $\varphi=\frac{1}{T}dT$ such that $d\alpha = \varphi \wedge \alpha$ when $n\geq 2$ or simply to the condition $\alpha \wedge d\alpha=0$ when $n\geq 3$ ([20], Th 6.4.6, p 245). Equivalently, we may use in both cases the {\it Frobenius theorem} saying that, for any couple of vector fields $\xi, \eta \in T$ such that $i(\xi)\alpha= 0, i(\eta)\alpha =0$ where $i( )$ is the {\it  interior product} of a vector with a form, then $i([\xi,\eta])\alpha =0$ because $(i(\xi)di(\eta)-i(\eta)d(i(\xi))\alpha=i(\xi)i(\eta)d\alpha + i([\xi,\eta])\alpha , \forall \xi,\eta \in T, \forall \alpha \in T^*$ from the definition of the exterior derivative on $1$-forms. However, what is surprisingly {\it not known at all} is the link existing between such conditions and group theory. We start with the following key definition (See Section 2B and [23] for more details):  \\

\noindent
{\bf DEFINITION 1.4}: A {\it Lie pseudogroup of transformations} is a group of transformations solutions of a (linear or even non-linear) system of ordinary or partial differential equations called system of {\it finite Lie equations}.  \\

\noindent
{\bf EXAMPLE 1.5}: When $n=1$ and we consider transformations $y=f(x)$ of the real line, the {\it affine group} of transformations is defined by the linear system $y_{xx}=0$ with jet notations saying that any transformation is such that ${\partial}_{xx}f(x)=0$ while the {\it projective group} of transformations is defined by the nonlinear system $\frac{y_{xxx}}{y_x}-\frac{3}{2}(\frac{y_{xx}}{y_x})^2=0$ with a similar comment. In both cases we have indeed a {\it Lie group of transformations} depending on a finite number of constant parameters, namely $y=ax+b$ in the first case and $y=\frac{ax+b}{cx+d}$ in the second case. Accordingly, the respective geometric object the invariance of which is characterizing the corresponding Lie pseudogroup are surely not made by tensors because the defining finite Lie equations are not first order. However, in the case of transformations of the plane $(x^1,x^2) \rightarrow (y^1,y^2)$ satisfying $y^2y^1_1-y^1y^2_1=x^2, y^2y^1_2-y^1y^2_2=-x^1\Rightarrow y^1_1y^2_2-y^1_2y^2_1=1$, no explicit integration can be obtained in order to provide general solutions but another way is to say that the corresponding Lie pseudogroup preserves the $1$-form $\alpha=x^2dx^1-x^1dx^2$ and thus the $2$-form $d\alpha=dx^1\wedge dx^2$ as we have indeed $y^2dy^1-y^1dy^2=x^2dx^1-x^1dx^2$ and thus also $dy^1\wedge dy^2=dx^1\wedge dx^2$. The Lie pseudogroup is thus preserving the {\it geometric object} $\omega =(\alpha, d\alpha)$ made by a $1$-form and a $2$-form. More generally, we may consider the Lie pseudogroup preserving the geometric object $\omega =(\alpha, \beta)$ where $\alpha$ is a $1$-form and $\beta$ is a $2$-form. As $d\alpha$ is also preserved, if we want that the system behaves at least like the preceding one, that is {\it cannot have any zero order equation}, we {\it must} have $d\alpha =c\beta$ for some arbitrary constant $c$. The two pseudogroups defined by $\omega \rightarrow c$ and $\bar{\omega}\rightarrow \bar{c}$ can be exchanged by a change of variables bringing $\omega$ to $\bar{\omega}$ {\it if and only if} $\bar{c}=c$. This situation is the simplest example of the celebrated {\it formal equivalence problem} ([20],[21]). \\

\noindent
{\bf EXAMPLE 1.6}: As a more general situation of a Lie pseudogroup of transformations of space with $n=3$ also involving differential forms, let us consider the $1$-forms $\alpha= x^3dx^1$ and $\beta = dx^2+ x^1dx^3$ with $\alpha \wedge d\alpha=0, \alpha + \beta = d(x^2 + x^1x^3)$. The Lie pseudogroup preserving $\alpha$ and $\beta$ also preserves $\gamma=d\alpha=-d\beta$ with $\alpha \wedge \gamma = 0$. It is easily seen to be made by the following transformations:  \\
\[   y^1=f(x^1), \hspace{5mm}   y^2=x^2+(x^1-\frac{f(x^1)}{f'(x^1)})x^3+a, \hspace{5mm} y^3=x^3/f'(x^1)               \]
where $f(x^1)$ is an arbitrary invertible function of $x^1$ only and we have set $f'(x^1)={\partial}_1f(x^1)$ while $a$ is an arbitrary constant, because we obtain at once $y^2+y^1y^3=x^2+x^1x^3 + a$ for an arbitrary constant $a$. An elementary but quite tedious computation similar to the previous one or to the ones that can be found in ([29], [30]) shows that solving the formal equivalence problem for $\omega=(\alpha, \beta)$ depends on the following {\it structure equations}:  \\
\[  d\alpha = c_1\alpha \wedge \beta + c'_1 \gamma, \hspace{5mm} d \beta = c_2 \alpha \wedge \beta + c'_2 \gamma, \hspace{5mm} d \gamma = c_3\beta \wedge \gamma Ê\]
because $\alpha \wedge \gamma=0,\hspace{3mm}  \alpha \wedge \beta \neq 0, \hspace{3mm}\gamma= dx^3\wedge dx^1\neq 0, \hspace{3mm}\beta \wedge \gamma =dx^1 \wedge dx^2 \wedge dx^3\neq 0 $. Closing this exterior system by taking again the exterior derivative, we get:  \\
\[     c'_1(c_1-c_3)=0, \hspace{8mm}  c'_1c_2-c'_2c_3=0     \]
In the present situation, we have $c_1=0, c'_1=1, c_2=0, c'_2=-1, c_3=0$. Eliminating $\gamma$, we get the only conditions:ÊÊ\\
\[    d(\alpha + \beta)=0, \hspace{1cm}  \alpha \wedge d\alpha=0  \]
that is exactly the conditions to be found in thermodynamics through a forthcoming example.\\
We invite the reader to choose:  \\
\[  \bar{\alpha}=dx^1, \hspace{1cm} \bar{\beta}=dx^2-x^3dx^1, \hspace{1cm} \bar{\gamma}= dx^1\wedge dx^3  \]
in order to obtain ${\bar{c}}_1=0, {\bar{c}}'_1=0, {\bar{c}}_2=0, {\bar{c}}'_2=1, {\bar{c}}_3=0$ with therefore no possibility to solve the equivalence problem $\omega \rightarrow \bar{\omega}$.  \\

Setting finally $\alpha=\delta Q=TdS, \beta=\delta W$ with $\alpha + \beta = dU$, the {\it Helmholtz postulate}, first stated in ([8]), assumes that it is always possible to choose the $n$ state variables, called {\it normal variables}, in such a way that $dT$ does not appear in $\delta W$. {\it This is a crucial assumption indeed} because, introducing the {\it free energy} $F=U-TS$, we get $dF=\delta W-TdS \Rightarrow S= - \frac{\partial F}{\partial T}$. We recall and improve the following result already provided in 1983 ([21], p 712-715) but never acknowledged up to now.\\

\noindent
{\bf THEOREM 1.7}: Helmholtz postulate is a theorem whenever $\alpha \wedge \beta \neq 0$.  \\

\noindent
{\bf Proof}: Let us prove first that, setting $\alpha=T(x)dS$ with $S=x^1$, it is always possible to choose the state variables in such a way that $dx^1$ does not appear in $\delta W$. \\
Starting with $n=2$, we get $\alpha \wedge \beta=\alpha\wedge (\alpha+\beta)=T(x)dx^1\wedge dU\neq 0$, implying that $U$ must not involve only $x^1$ and we may introduce the new variables $y^1=x^1, y^2=U(x)$ in such a way that $\beta= dy^2-T(y)dy^1$. Let now $v(y)$ be a {\it non-constant} orbital integral of the ordinary differential equation $dy^2/dy^1=T(y^1,y^2)$ satisfying therefore $\frac{\partial v}{\partial y^1}+T(y)\frac{\partial v}{\partial y^2}=0$. It follows that $dv=\frac{\partial v}{\partial y^1}dy^1+\frac{\partial v}{\partial y^2}dy^2=\frac{\partial v}{\partial y^2}\beta$ with $\frac{\partial v}{\partial y^2}\neq 0$ because otherwise we should also have $\frac{\partial v}{\partial y^1}=0$ and $v$ should be constant. Using the new variables $z^1=y^1, z^2=v(y^1,y^2)$, we have a jacobian $\frac{\partial (z^1,z^2)}{\partial (y^1,y^2)}=\frac{\partial v}{\partial y^2}\neq 0$ and obtain at once $\alpha=T(z)dz^1, \beta=b(z)dz^2$. Finally, setting $\bar{\alpha}= - SdT$, we have $d\bar{\alpha}=dT\wedge dS=d\alpha$ but now $\beta +\bar{\alpha}=\beta + \alpha - d(TS)=d(U-TS)=dF$, that is we may exchange $\alpha,\beta,U$ with $\bar{\alpha}, \beta, F$ and repeat the same procedure with $T$ in place of $S$ and $F$ in place of $U$, obtaining therefore the desired result. \\
Similarly, when $n\geq 3$, we can choose the new variables $y^1=x^1, y^2=U, y^3=x^3,...,y^n=x^n$ and obtain $\alpha=T(y)dy^1, \beta=dy^2-T(y)dy^1$. Considering now $y^3,...,y^n$ like parameters, we may use the same argument as above and substitute $\frac{\partial v}{\partial y^1}= - T(y)\frac{\partial v}{\partial y^2}$ in order to get:  \\
\[ dv=\frac{\partial v}{\partial y^1}dy^1+ \frac{\partial v}{\partial y^2}dy^2+ \frac{\partial v}{\partial y^3}dy^3+ ...+ \frac{\partial v}{\partial y^n}dy^n
       = \frac{\partial v}{\partial y^2}\beta+ \frac{\partial v}{\partial y^3}dy^3+ ...+ \frac{\partial v}{\partial y^n}dy^n   \]
Choosing $z^1=y^1, z^2=v(y), z^3=y^3,...,z^n=y^n$, we obtain $\alpha=T(z)dz^1, \beta= b_2(z)dz^2 + ... +b_n(z)dz^n$. The final exchange may be done as before.  \\
\hspace*{12cm}   Q.E.D.   \\

\noindent
{\bf EXAMPLE 1.8}: In the case of an ideal gas with $n=2$, we may choose $y^1=S, y^2=U=CT$ and we have $\alpha = Tdy^1, \beta =dy^2-Tdy^1$. Meanwhile, we have also $TdS=CdT+PdV=dU+PdV \Rightarrow dV=\frac{T}{P}dS- \frac{1}{P}dU$. It follows that $\frac{\partial V}{\partial y^1}=\frac{T}{P}, \frac{\partial V}{\partial y^2}= - \frac{1}{P}$ and thus $\frac{\partial V}{\partial y^1} + T \frac{\partial V}{\partial y^2}=0$. Accordingly, $V(S,U)$ can be chosen to be the desired orbital integral, a result highly not evident at first sight but explaining the notations. \\

\noindent
{\bf EXAMPLE 1.9}: With $n=3$ and local coordinates $x=(x^1,x^2,x^3)$ for the state variables, let us consider an abstract system with $\delta Q=\alpha=x^3dx^1, \delta W=\beta= dx^2+x^1dx^3$. We have indeed $d\alpha=dx^3\wedge dx^1=\frac{1}{x^3}dx^3\wedge \alpha \Rightarrow \alpha \wedge d\alpha=0$ and $\alpha + \beta =d(x^2+x^1x^3)$. We may therefore set $S=x^1, T=x^3, U=x^2+x^1x^3$ and the existence of the integrating factor is compatible with the change of scale allowing to define $T$. However, we should get $F=U-TS=x^2 \Rightarrow dF=dx^2$ and we should be tempted to conclude with a contradiction as we should get $S=-\frac{\partial F}{\partial T}=-{\partial}_3F=0\neq x^1$. However, {\it things are much more subtle} when dealing with normal variables as it has been largely emphasized by Duhem in ([8]) but totally absent from the survey reference ([3]). Indeed, we have now $dF=\delta W-SdT \Rightarrow \delta W=dF+SdT=dx^2+x^1dx^3$ in a coherent way with the definition of $\beta$. Accordingly, the correct way is thus to say that the formula $S=-\frac{\partial F}{\partial T}$ is no longer true because $\delta W$ now contains $dT$ or, equivalently, that {\it the state variables x are not normal}. However, exchanging $U$ and $F$, it follows from our proof of the Helmholtz postulate that {\it it is always possible to obtain normal state variables} $y=(y^1,y^2,y^3)$. For this, we just need to set $y^1=x^1, y^2=x^2+x^1x^3,y^3=x^3 \Leftrightarrow x^1=y^1, x^2=y^2-y^1y^3, x^3=y^3$ and obtain $\delta Q=\alpha= y^3dy^1, \delta W= \beta= dy^2-y^3dy^1$ where $\beta$ {\it does not contain} $dy^3=dx^3=dT$ {\it any longer}. Meanwhile, $dF=\delta W-SdT=d(y^2-y^1y^3)=dx^2$ as before but now $F=y^2-y^1y^3$ is such that $S=x^1=y^1= -\frac{\partial F}{\partial y^3}=-\frac{\partial F}{\partial T}$ as we wished.  \\

\noindent
{\bf 2) MATHEMATICAL TOOLS}\\

\noindent
{\bf A) LIE GROUPS}\\

The word "{\it group}" has been introduced for the first time in 1830 by E. Galois and this concept slowly passed from algebra (groups of permutations) to geometry (groups of transformations). It is only in 1880 that S. Lie studied the groups of transformations depending on a finite number of parameters and now called {\it Lie groups of transformations}. We now describe in a modern language the procedure followed by Poincar\'{e} in [19], both with the corresponding {\it dual variational framework}. We invite the reader to look at ([25], [26], [30], [31]) in order to discover its link with {\it homological algebra} and the {\it extension functor}.\\

 Let $X$ be a manifold with local coordinates $x=(x^1, ... , x^n)$ and $G$ be a {\it Lie group}, that is another manifold with local coordinates $a=(a^1, ... , a^p)$ called {\it parameters}, with a {\it composition} $G\times G \rightarrow G: (a,b)\rightarrow ab$, an {\it inverse} $G \rightarrow G: a \rightarrow a^{-1}$ and an {\it identity} $e\in G$ satisfying:\\
\[(ab)c=a(bc)=abc,\hspace{1cm} aa^{-1}=a^{-1}a=e,\hspace{1cm} ae=ea=a,\hspace{1cm} \forall a,b,c \in G \]
Then $G$ is said to {\it act} on $X$ if there is a map $X\times G \rightarrow X: (x,a) \rightarrow y=ax=f(x,a)$ such that $(ab)x=a(bx)=abx, \forall a,b\in G, \forall x\in X$ and, for simplifying the notations, we shall use global notations even if only local actions are existing. The action is said to be {\it effective} if $ax=x, \forall x\in X\Rightarrow a=e$. A subset $S\subset X$ is said to be {\it invariant} under the action of $G$ if $aS\subset S,\forall a\in G$ and the {\it orbit} of $x\in X$ is the invariant subset $Gx=\{ax\mid a\in G\}\subset X$. If $G$ acts on two manifolds $X$ and $Y$, a map $f:X\rightarrow Y$ is said to be {\it equivariant} if $f(ax)=af(x), \forall x\in X, \forall a\in G$. For reasons that will become clear later on, it is often convenient to introduce the {\it graph} $X\times G\rightarrow X\times X: (x,a)\rightarrow (x,y=ax)$ of the action. In the product $X\times X$, the first factor is called the {\it source} while the second factor is called the {\it target}. \\

We denote as usual by $T=T(X)$ the {\it tangent bundle} of $X$, by $T^*=T^*(X)$ the {\it cotangent bundle}, by ${\wedge}^rT^*$ the {\it bundle of r-forms} and by $S_qT^*$ the {\it bundle of q-symmetric tensors}. Moreover, if  $\xi,\eta\in T$ are two vector fields on $X$, we may define their {\it bracket} $[\xi,\eta]\in T$ by the local formula $([\xi,\eta])^i={\xi}^r{\partial}_r{\eta}^i-{\eta}^s{\partial}_s{\xi}^i$ leading to the {\it Jacobi identity} $[\xi,[\eta,\zeta]]+[\eta,[\zeta,\xi]]+[\zeta,[\xi,\eta]]=0, \forall \xi,\eta,\zeta \in T$ allowing to define a {\it Lie algebra}. We have also the useful formula $[T(f)(\xi),T(f)(\eta)]=T(f)([\xi,\eta])$ where $T(f):T(X)\rightarrow T(Y)$ is the tangent mapping of a map $f:X\rightarrow Y$. If $\xi \in T$ and $f\in C^{\infty}(X)$, we set $\xi . f={\xi}^i{\partial}_if$ and, if $\omega \in {\wedge}^rT^*$, we denote by $i(\xi)\omega \in {\wedge}^{r-1}T^*$ the {\it interior product} of $\omega$ by $\xi$. Finally, when $I=\{ i_1< ... < i_r\}$ is a multi-index, we may set $dx^I=dx^{i_1}\wedge ... \wedge dx^{i_r}$ and introduce the {\it exterior derivative} $d:{\wedge}^rT^*\rightarrow {\wedge}^{r+1}T^*:\omega={\omega}_Idx^I \rightarrow d\omega={\partial}_i{\omega}_Idx^i\wedge dx^I$ with $d^2=d\circ d\equiv 0$ because ${\partial}_{ij}{\omega}_I dx^i \wedge dx^j \wedge dx^I \equiv 0$, in the {\it Poincar\'{e} sequence}:\\
\[  {\wedge}^0T^* \stackrel{d}{\longrightarrow} {\wedge}^1T^* \stackrel{d}{\longrightarrow} {\wedge}^2T^* \stackrel{d}{\longrightarrow} ... \stackrel{d}{\longrightarrow} {\wedge}^nT^* \longrightarrow 0  \]

In order to fix the notations, we quote without any proof a few results that will be of constant use in the sequel (See [23] for more details).\\

According to the {\it first fundamental theorem of Lie}, the orbits $x=f(x_0,a)$ satisfy the system of PD equations $\partial x^i/\partial a^{\sigma}= {\theta}^i_{\rho}(x){\omega}^{\rho}_{\sigma}(a)$  with $det(\omega)\neq 0$. The vector fields ${\theta}_{\rho}={\theta}^i_{\rho}(x){\partial}_i$ are called {\it infinitesimal generators} of the action and are linearly independent over the constants when the action is effective. In a rough symbolic way, we have $x=ax_0\Rightarrow dx=dax_0=daa^{-1}x$ and $daa^{-1}=\omega=({\omega}^{\tau}={\omega}^{\tau}_{\sigma}(a)da^{\sigma})$ is thus a family of right invariant 1-forms on $G$ with value in ${\cal{G}}=T_e(G)$ the tangent space to $G$ at the identity $e\in G$, called {\it Maurer-Cartan} (MC) {\it forms}.  \\

Then, according to the {\it second fundamental theorem of Lie}, if ${\theta}_1,...,{\theta}_p$ are the infinitesimal generators of the effective action of a lie group $G$ on $X$, then $[{\theta}_{\rho},{\theta}_{\sigma}]=c^{\tau}_{\rho\sigma}{\theta}_{\tau}$ where the $c=(c^{\tau}_{\rho\sigma}= - c^{\tau}_{\sigma\rho})$ are the {\it structure constants} of a Lie algebra of vector fields which can be identified with ${\cal{G}}$ by using the action as we already did. Equivalently, introducing the non-degenerate inverse matrix $\alpha={\omega}^{-1}$ of right invariant vector fields on $G$, we obtain from crossed-derivatives the {\it compatibility conditions} (CC) for the previous system of partial differential (PD) equations called {\it Maurer-Cartan} (MC) {\it equations}, namely: \\
\[  \frac{\partial {\omega}^{\tau}_s}{\partial a^r} -\frac{ \partial {\omega}^{\tau}_r}{ \partial a^s} + c^{\tau}_{\rho\sigma} {\omega}^{\rho}_r {\omega}^{\sigma}_s  = 0 \hspace{5mm} \Leftrightarrow \hspace{5mm} d{\omega}^{\tau}+\frac{1}{2} c^{\tau}_{\rho\sigma}{\omega}^{\rho}\wedge {\omega}^{\sigma} = 0  \]
({\it care to the sign used}) or equivalently $[{\alpha}_{\rho},{\alpha}_{\sigma}]=c^{\tau}_{\rho\sigma} {\alpha}_{\tau} $. \\
 
Finally, using again crossed-derivatives, we obtain the corresponding {\it integrability conditions} (IC) on the structure constants:\\
\[  c^{\tau}_{\rho\sigma}+c^{\tau}_{\sigma \rho}=0, \hspace{1cm} c^{\lambda}_{\mu\rho}c^{\mu}_{\sigma\tau}+c^{\lambda}_{\mu\sigma}c^{\mu}_{\tau\rho}+c^{\lambda}_{\mu\tau}c^{\mu}_{\rho\sigma}=0  \]
also called {\it Jacobi conditions}. The Cauchy-Kowaleski theorem finally asserts that one can construct an analytic group $G$ such that ${\cal{G}}=T_e(G)$ by recovering the MC forms from the MC equations, a result amounting to the {\it third fundamental theorem of Lie}.\\

\noindent
{\bf EXAMPLE 2.A.1}: Considering the affine group of transformations of the real line $y=a^1x+a^2$, the orbits are defined by $x=a^1x_0+a^2$, a definition leading to $dx=da^1x_0+da^2$ and thus $dx=((1/a^1)da^1)x+(da^2-(a^2/a^1)da^1)$. We obtain therefore ${\theta}_1=x{\partial}_x, {\theta}_2={\partial}_x \Rightarrow [{\theta}_1,{\theta}_2]=-{\theta}_2$ and ${\omega}^1=(1/{a^1})da^1, {\omega}^2=da^2-(a^2/{a^1})da^1\Rightarrow d{\omega}^1=0, d{\omega}^2-{\omega}^1\wedge{\omega}^2=0 \Leftrightarrow [{\alpha}_1,{\alpha}_2]=-{\alpha}_2$ with ${\alpha}_1=a^1{\partial}_1+a^2{\partial}_2, {\alpha}_2={\partial}_2$.\\
 
\noindent
{\bf EXAMPLE 2.A.2}: If $x=a(t)x_0+b(t)$ with $a(t)$ a time depending orthogonal matrix ({\it rotation}) and $b(t)$ a time depending vector ({\it translation}) describes the movement of a rigid body in ${\mathbb{R}}^3$, then the projection of the {\it absolute speed} $v=\dot{a}(t)x_0+\dot{b}(t)$ in an orthogonal frame fixed in the body is the so-called {\it relative speed} $a^{-1}v=a^{-1}\dot{a}x_0+a^{-1}\dot{b}$ and the kinetic energy/Lagrangian is a quadratic function of the $1$-forms $A=(a^{-1}\dot{a}$, $a^{-1}\dot{b})$. Meanwhile, taking into account the preceding example, the {\it Eulerian speed}  $v=v(x,t)=\dot{a}a^{-1}x+\dot{b}-\dot{a}a^{-1}b$ only depends on the 1-forms $B=(\dot{a}a^{-1}, \dot{b}-\dot{a}a^{-1}b)$. We notice that $a^{-1}\dot{a}$ and $\dot{a}a^{-1}$ are both $3\times 3$ skewsymmetric time depending matrices that may be quite different.\\

The above particular case, well known by anybody studying the analytical mechanics of rigid bodies, can be generalized as follows. If $X$ is a manifold and $G$ is a lie group ({\it not acting necessarily on} $X$ {\it now}), let us consider maps $a:X\rightarrow G: (x)\rightarrow (a(x))$ or equivalently sections of the trivial (principal) bundle $X\times G$ over $X$, namely maps $X \rightarrow X\times G:(x) \rightarrow (x,a(x))$. If $x+dx$ is a point of $X$ "close " to $x$, then $T(a)$ will provide a point $a+da=a+\frac{\partial a}{\partial x}dx$ "close " to $a$ on $G$. We may bring $a$ back to $e$ on $G$ by acting on $a$ with $a^{-1}$, {\it  either on the left or on the right}, getting therefore a $1$-form $a^{-1}da=A$ or $daa^{-1}=B$ with value in $\cal{G}$. As $aa^{-1}=e$ we also get $a^{-1}da=-(da^{-1})a=-dbb^{-1}$ if we set $b=a^{-1}$ as a way to link $A$ with $B$. When there is an action $y=ax$, we have $x=a^{-1}y=by$ and thus $dy=dax=daa^{-1}y$, a result leading to the equivalent formulas:\\
\[   a^{-1}da=A=({A}^{\tau}_i(x)dx^i=-{\omega}^{\tau}_{\sigma}(b(x)){\partial}_ib^{\sigma}(x)dx^i)  \]
\[   daa^{-1}=B=({B}^{\tau}_i(x)dx^i={\omega}^{\tau}_{\sigma}(a(x)){\partial}_ia^{\sigma}(x)dx^i)  \]
Introducing the induced bracket $[A,A](\xi,\eta)=[A(\xi),A(\eta)]\in {\cal{G}}, \forall \xi,\eta\in T$, we may define the {\it curvature} $2$-form $dA-[A,A]=F\in {\wedge}^2T^*\otimes {\cal{G}}$ by the local formula ({\it care again to the sign}):\\
\[     {\partial}_iA^{\tau}_j(x)-{\partial}_jA^{\tau}_i(x)-c^{\tau}_{\rho\sigma}A^{\rho}_i(x)A^{\sigma}_j(x)=F^{\tau}_{ij}(x)  \]
This definition can also be adapted to $B$ by using $dB+[B,B]$ and we obtain:\\

\noindent
{\bf THEOREM 2.A.3}: There is a {\it nonlinear gauge sequence}:\\
\[  \begin{array}{ccccc}
X\times G & \longrightarrow & T^*\otimes {\cal{G}} &\stackrel{MC}{ \longrightarrow} & {\wedge}^2T^*\otimes {\cal{G}}  \\
a                & \longrightarrow  &    a^{-1}da=A         &    \longrightarrow & dA-[A,A]=F
\end{array}   \]

Choosing $a$ "close" to $e$, that is $a(x)=e+t\lambda(x)+...$ and linearizing as usual, we obtain the linear operator $d:{\wedge}^0T^*\otimes {\cal{G}}\rightarrow {\wedge}^1T^*\otimes {\cal{G}}:({\lambda}^{\tau}(x))\rightarrow ({\partial}_i{\lambda}^{\tau}(x))$ leading to:\\

\noindent
{\bf COROLLARY 2.A.4}: There is a {\it linear gauge sequence}:\\ 
\[  {\wedge}^0T^*\otimes {\cal{G}}\stackrel{d}{\longrightarrow} {\wedge}^1T^*\otimes {\cal{G}} \stackrel{d}{\longrightarrow} {\wedge}^2T^*\otimes{\cal{G}} \stackrel{d}{\longrightarrow} ... \stackrel{d}{\longrightarrow} {\wedge}^nT^*\otimes {\cal{G}}\longrightarrow  0   \]
which is the tensor product by $\cal{G}$ of the Poincar\'{e} sequence:\\

 It remains to introduce the previous results into a variational framework. The procedure has been found in 1901 by H. Poincar\'{e} who introduced a {\it duality principle} in analytical mechanics in order to study lagrangians invariant under the action of a Lie group of transformations ([19]). This method has been used later on by G. Birkhoff in 1954 ([2]) and V. Arnold in 1966 ([1]), each one omitting to quote the previous results.  \\
 
For this, we may consider a lagrangian on $T^*\otimes \cal{G}$, that is an {\it action} $W=\int w(A)dx$ where $dx=dx^1\wedge ...\wedge dx^n$ and vary it. With $A=a^{-1}da=-dbb^{-1}$ we may introduce $\lambda=a^{-1}\delta a=-\delta bb^{-1}\in {\cal{G}}={\wedge}^0T^*\otimes {\cal{G}}$ with local coordinates ${\lambda}^{\tau}(x)=-{\omega}^{\tau}_{\sigma}(b(x))\delta b^{\sigma}(x)$ and we obtain $ \delta A=d\lambda - [A,\lambda]  $ that is $\delta A^{\tau}_i={\partial}_i\lambda^{\tau}-c^{\tau}_{\rho\sigma}A^{\rho}_i{\lambda}^{\sigma}$ in local coordinates. Then, setting $\partial w/\partial A={\cal{A}}=({\cal{A}}^i_{\tau})\in {\wedge}^{n-1}T^*\otimes \cal{G}$, we get:\\
\[  \delta W=\int {\cal{A}}\delta Adx=\int {\cal{A}}(d\lambda-[A,\lambda])dx  \]
and, after integration by part, the Euler-Lagrange (EL) {\it relative equations} ([22],[23]):\\
\[     {\partial}_i{\cal{A}}^i_{\tau}+c^{\sigma}_{\rho\tau}A^{\rho}_i{\cal{A}}^i_{\sigma}=0    \]
Such a linear operator for $\cal{A}$ has {\it non constant coefficients linearly depending on} $A$ {\it and the structure constants}. Setting $\delta aa^{-1}=\mu\in {\cal{G}}$, we get $\lambda=a^{-1}(\delta aa^{-1})a=Ad(a)\mu$ while, setting $a \rightarrow a'=ab$, we get the {\it gauge transformation} $A \rightarrow A'=(ab)^{-1}d(ab)=b^{-1}a^{-1}(dab+adb)=Ad(b)A+b^{-1}db, \forall b\in G$. Setting $b=e+t\lambda+...$ with $t\ll 1$, then $\delta A$ becomes an infinitesimal gauge transformation. However, setting now $a \rightarrow a'=ca$, we get $ A'=a^{-1}c^{-1}(dca+cda)=a^{-1}(c^{-1}dc)a+A$ and thus $ \delta A=Ad(a)d\mu$ when $c=aba^{-1}=e+t\mu +...$ with $t\ll 1$ (See [23], p 180, 424 for more details and computations using local coordinates). We may also notice that $aa^{-1}=e \Rightarrow \delta a a^{-1}+a \delta (a^{-1})=0 \Rightarrow \delta (a^{-1})=-a^{-1}\delta a a^{-1}$ and thus:  \\
\[  \begin{array}{rcl}
\delta A  & = & -a^{-1}\delta aa^{-1}da+a^{-1}d((\delta aa^{-1})a)   \\
                & = & -a^{-1}\delta aa^{-1}da +a^{-1}d(\delta aa^{-1})a+a^{-1}\delta aa^{-1}da  \\
                & = &  Ad(a)d\mu
    \end{array}    \]
Therefore, introducing by duality $\cal{B}$ such that ${\cal{B}}\mu= {\cal{A}}\lambda$, we get the divergence-like {\it absolute equations} $ {\partial}_i{\cal{B}}^i_{\sigma}=0$. When $n=1$, we recognize at once the Birkhoff-Arnold dynamics of a rigid body, with time $t$ as independent variable, or the Kirchhoff-Love theory of a thin elastic beam, with curvilinear abcissa $s$ along the beam as independent variable. \\

\noindent
{\bf REMARK 2.A.5}: As the passage from $A$ to $B$, that is from left invariance to right invariance is not easy to achieve in actual practice, we indicate a way to simplify the use of the adjoint mapping (Compare to [23], Proposition 10, p 180). Indeed, working formally, from $\delta A= d\lambda-cA\lambda$, we may define on $G$ a square matrix acting on $\cal{G}$ and define $\mu$ by $\lambda= M\mu$. Substituting, we obtain $\delta A=d (M\mu)-cAM\mu=Md \mu+ (dM-cAM)\mu$ and thus $\delta A=Md\mu \Leftrightarrow dM-cAM=0$, that is to say $M(b)$ must be a solution of the linear system of PD equations $\frac{\partial M^{\tau}_{\mu}}{\partial b^{r}} + c^{\tau}_{\rho \sigma}{\omega}^{\rho}_r(b)M^{\sigma}_{\mu}=0$. It just remains to prove that this system is involutive by computing the crossed derivatives. An easy but tedious computation provides: \\
\[  (c^{\tau}_{\rho\sigma}(\frac{\partial {\omega}^{\rho}_s}{\partial b^r}-\frac{\partial {\omega}^{\rho}_r}{\partial b^s}+c^{\rho}_{\alpha \beta}{\omega}^{\alpha}_r{\omega}^{\beta}_s) - (c^{\tau}_{\rho\sigma}c^{\rho}_{\alpha \beta}+c^{\tau}_{\rho\alpha}c^{\rho}_{\beta\sigma}+c^{\tau} _{\rho\beta}c^{\rho}_{\sigma \alpha}){\omega}^{\alpha}_r{\omega}^{\beta}_s)M^{\sigma}_{\mu}=0   \]
{\it Using both the} MC {\it equations and the Jacobi conditions} achieves the proof of this technical but quite useful result.  \\

We may therefore ask: \\

\noindent
{\bf PROBLEM}: HOW IS IT POSSIBLE AND WHY IS IT EVEN NECESSARY TO INTRODUCE DIFFERENT EQUATIONS WITHIN THE SAME 
GROUP BACKGROUND.  \\

\noindent
{\bf B) LIE PSEUDOGROUPS}\\

We start recalling a few notations and definitions about fibered manifolds and their jet bundles (See [20] and [22] for more details). In particular, if ${\cal{E}}\rightarrow X:(x,y)\rightarrow (x)$ is a fibered manifold with changes of local coordinates having the form $\bar{x}=\varphi (x), \bar{y}=\psi (x,y)$, we shall denote by $J_q({\cal{E}})\rightarrow X:(x,y_q)\rightarrow (x) $ the $q$-{\it jet bundle} of ${\cal{E}}$ with local coordinates $(x^i,y^k_{\mu})$ for $i=1,...,n$, $k=1,...,m$, $0\leq \mid \mu \mid\leq q$ and $y^k_0=y^k$. We may consider sections $f_q:(x)\rightarrow (x,f^k(x), f^k_i(x), f^k_{ij}(x), ...)=(x,f_q(x))$ transforming like the sections $j_q(f):(x) \rightarrow (x,f^k(x), {\partial}_if^k(x), {\partial}_{ij}f^k(x), ...)=(x,j_q(f)(x))$ where both $f_q$ and $j_q(f)$ are over the section $f:(x)\rightarrow (x,y^k=f^k(x))=(x,f(x))$ of $\cal{E}$ transforming like $\bar{f}(\varphi (x)) = \psi (x, f(x))$. If $T(\cal{E})$ has local coordinates $(x,y;u,v)$, we shall denote by $V(\cal{E})$ the {\it vertical bundle} of $\cal{E}$, namely the sub-vector bundle of $T(\cal{E})$ with local coordintates $(x,y;0,v)$. The (nonlinear) {\it Spencer operator} just allows to distinguish a section $f_q$ from a section $j_q(f)$ by introducing a kind of "{\it difference} " through the operator $D:J_{q+1}({\cal{E}})\rightarrow T^*\otimes V(J_q({\cal{E}})): f_{q+1}\rightarrow j_1(f_q)-f_{q+1}$ with local components $({\partial}_if^k(x)-f^k_i(x), {\partial}_if^k_j(x)-f^k_{ij}(x),...) $ and more generally $(Df_{q+1})^k_{\mu,i}(x)={\partial}_if^k_{\mu}(x)-f^k_{\mu+1_i}(x)$. If $m=n$ and ${\cal{E}}=X\times X$ with source projection, we denote by ${\Pi}_q={\Pi}_q(X,X)\subset J_q(X\times X)$ the open sub-bundle locally defined by $det(y^k_i)\neq 0$ and we shall set $\Delta=det({\partial}_if^k(x))$. Also, if $\cal{E}$ and $\cal{F}$ are two fibered manifolds over $X$ with local coordinates $(x,y)$ and $(x,z)$ respectively, we shall denote by ${\cal{E}}{\times}_X{\cal{F}}$ their {\it fibered product} over $X$ with local coordinates $(x,y,z)$. Finally, if $E$ is a {\it vector bundle} over $X$ with transition rules having the form $\bar{x}=\varphi(x), \bar{y}=A(x)y$, we shall denote by $E^*$ the vector bundle obtained from $E$ by inverting the transition matrices, exactly like $T^*$ is obtained from $T$. \\

In 1890, Lie discovered that {\it Lie groups of transformations} were examples of {\it Lie pseudogroups of transformations} along the following definition which expands the preliminary Definition 1.4: \\
 
 \noindent
{\bf DEFINITION 2.B.1}: A {\it Lie pseudogroup of transformations} $\Gamma\subset aut(X)$ is a group of transformations solutions of a system of OD or PD equations such that, if $y=f(x)$ and $z=g(y)$ are two solutions, called {\it finite transformations}, that can be composed, then $z=g\circ f(x)=h(x)$ and $x=f^{-1}(y)=g(y)$ are also solutions while $y=x$ is the {\it identity} solution denoted by $id=id_X$ and we shall set $id_q=j_q(id)$.  In all the sequel we shall suppose that $\Gamma$ is {\it transitive} that is $\forall x,y\in X, \exists f\in \Gamma, y=f(x)$ \\

We notice that an action $y=f(x,a)$ provides a Lie pseudogroup by eliminating the $p$ parameters $a$ among the equations $y_q=j_q(f)(x,a)$ obtained by successive differentiations with respect to $x$ only when $q$ is large enough. The system ${\cal{R}}_q\subset {\Pi}_q$ of OD or PD equations thus obtained may be quite nonlinear and of high order. {\it The concept of parameters is not existing in this new framework} and thus no one of the methods already presented may be used any longer. Setting $f(x)=f(x,a(x))$ and $f_q(x)=j_q(f)(x,a(x))$, we obtain $a(x)=a=cst \Leftrightarrow f_q=j_q(f)$ because 
$Df_{q+1}=j_1(f_q)-f_{q+1}=(\partial f_q(x,a(x))/\partial a^{\tau}){\partial}_ia^{\tau}(x)dx^i$ as a $1$-form and the matrix involved has rank $p$ in the following commutative diagram:  \\
\[  \begin{array}{rcccl}
0\rightarrow  & \hspace{5mm} X  \times  G &  =  &  \hspace{4mm}{\cal{R}}_q  &  \rightarrow 0  \\
      &   a=cst \uparrow \downarrow \uparrow a(x)  &    &  j_q(f) \uparrow \downarrow \uparrow  f_q &  \\
         &  \hspace{4mm}X     &  =   &    \hspace{4mm}X   &   
\end{array}  \]

More generally, looking now for transformations "close" to the identity, that is setting $y=x+t\xi(x)+...$ when $t\ll 1$ is a small constant parameter and passing to the limit $t\rightarrow 0$, we may linearize any (nonlinear) {\it system of finite Lie equations} in order to obtain a (linear) {\it system of infinitesimal Lie equations} $R_q\subset J_q(T)$ for vector fields. Such a system has the property that, if $\xi,\eta$ are two solutions, then $[\xi,\eta]$ is also a solution. Accordingly, the set $\Theta\subset T$ of its solutions satisfies $[\Theta,\Theta]\subset \Theta$ and can therefore be considered as the Lie algebra of $\Gamma$. \\

Looking at the way a vector field and its derivatives are transformed under any $f\in aut(X)$ while replacing $j_q(f)$ by $f_q$, we obtain:\\
\[  {\eta}^k(f(x))=f^k_r(x){\xi}^r(x) \Rightarrow {\eta}^k_u(f(x))f^u_i(x)=f^k_r(x){\xi}^r_i(x)+f^k_{ri}(x){\xi}^r(x)  \]
and so on, a result leading to:\\

\noindent
{\bf LEMMA 2.B.2}: $J_q(T)$ is {\it associated} with ${\Pi}_{q+1}={\Pi}_{q+1}(X,X)$ that is we can obtain a new section ${\eta}_q=f_{q+1}({\xi}_q)$ from any section ${\xi}_q \in J_q(T)$ and any section $f_{q+1}\in {\Pi}_{q+1}$ by the formula:\\
\[ d_{\mu}{\eta}^k\equiv {\eta}^k_rf^r_{\mu}+ ...=f^k_r{\xi}^r_{\mu}+  ...  +f^k_{\mu+1_r}{\xi}^r , \forall 0\leq {\mid}\mu {\mid}\leq q\]
where the left member belongs to $V({\Pi}_q)$. Similarly $R_q\subset J_q(T)$ is associated with ${\cal{R}}_{q+1}\subset {\Pi}_{q+1}$.\\

We now need a few basic definitions on {\it Lie groupoids} and {\it Lie algebroids} that will become substitutes for Lie groups and Lie algebras. The first idea is to use the chain rule for derivatives $j_q(g\circ f)=j_q(g)\circ j_q(f)$ whenever $f,g\in aut(X)$ can be composed and to replace both $j_q(f)$ and $j_q(g)$ respectively by $f_q$ and $g_q$ in order to obtain the new section $g_q\circ f_q$. This kind of "{\it composition}" law can be written in a pointwise symbolic way by introducing another copy $Z$ of $X$ with local coordinates $(z)$ as follows:\\
\[ {\gamma}_q:{\Pi}_q(Y,Z){\times}_Y{\Pi}_q(X,Y)\rightarrow {\Pi}_q(X,Z):((y,z,\frac{\partial z}{\partial y},...),(x,y,\frac{\partial y}{\partial x},...)\rightarrow (x,z,\frac{\partial z}{\partial y}\frac{\partial y}{\partial x},...)      \]
We may also define $j_q(f)^{-1}=j_q(f^{-1})$ and obtain similarly an "{\it inversion}" law.\\

\noindent
{\bf DEFINITION 2.B.3}: A fibered submanifold ${\cal{R}}_q\subset {\Pi}_q$ is called a {\it system of finite Lie equations} or a {\it Lie groupoid} of order $q$ if we have an induced {\it source projection} ${\alpha}_q:{\cal{R}}_q\rightarrow X$, {\it target projection} ${\beta}_q:{\cal{R}}_q\rightarrow X$, {\it composition} ${\gamma}_q:{\cal{R}}_q{\times}_X{\cal{R}}_q\rightarrow {\cal{R}}_q$, {\it inversion} ${\iota}_q:{\cal{R}}_q\rightarrow {\cal{R}}_q$ 
and {\it identity} $id_q:X\rightarrow {\cal{R}}_q$. In the sequel we shall only consider {\it transitive} Lie groupoids such that the map $({\alpha}_q,{\beta}_q):{\cal{R}}_q\rightarrow X\times X $ is an epimorphism. One can prove that the new system ${\rho}_r({\cal{R}}_q)={\cal{R}}_{q+r}=J_r({\cal{R}}_q)\cap {\Pi}_{q+r}\subset J_r({\Pi}_q)$ obtained by differentiating $r$ times all the defining equations of ${\cal{R}}_q$ is a Lie groupoid of order $q+r$. \\

Using the {\it algebraic bracket} $\{ j_{q+1}(\xi),j_{q+1}(\eta)\}=j_q([\xi,\eta]), \forall \xi,\eta\in T$, we may  obtain by bilinearity a {\it differential bracket} on $J_q(T)$ extending the bracket on $T$:\\
\[   [{\xi}_q,{\eta}_q]=\{{\xi}_{q+1},{\eta}_{q+1}\}+i(\xi)D{\eta}_{q+1}-i(\eta)D{\xi}_{q+1}, \forall {\xi}_q,{\eta}_q\in J_q(T) \]
which does not depend on the respective lifts ${\xi}_{q+1}$ and ${\eta}_{q+1}$ of ${\xi}_q$ and ${\eta}_q$ in $J_{q+1}(T)$. One can prove that this bracket on sections satisfies the Jacobi identity and we set: \\

\noindent
{\bf DEFINITION 2.B.4}: We say that a vector subbundle $R_q\subset J_q(T)$ is a {\it system of infinitesimal Lie equations} or a {\it Lie algebroid} if $[R_q,R_q]\subset R_q$, that is to say $[{\xi}_q,{\eta}_q]\in R_q, \forall {\xi}_q,{\eta}_q\in R_q$. Such a definition can be tested by means of computer algebra. We shall also say that $R_q$ is {\it transitive} if we have the short exact sequence $0\rightarrow R^0_q \rightarrow R_q \stackrel{{\pi}^q_0}{ \rightarrow} T  \rightarrow 0$. In that case, a {\it splitting} of this sequence, namely a map ${\chi}_q:T \rightarrow R_q$ such that $ {\pi}^q_0\circ {\chi}_q=id_T$ or equivalently a section ${\chi}_q\in T^*\otimes R_q$ over $id_T\in T^*\otimes T$, is called a $R_q$-{\it connection} and its {\it curvature} ${\kappa}_q\in {\wedge}^2T^*\otimes R^0_q$ is defined by ${\kappa}_q(\xi,\eta)=[{\chi}_q(\xi),{\chi}_q(\eta)]-{\chi}_q([\xi,\eta]), \forall \xi,\eta \in T$.\\

\noindent
{\bf PROPOSITION 2.B.5}: If $[R_q,R_q]\subset R_q$, then $[R_{q+r},R_{q+r}]\subset R_{q+r}, \forall r\geq 0$.  \\

\noindent
{\bf Proof}: When $r=1$, we have ${\rho}_1(R_q)=R_{q+1}=\{ {\xi}_{q+1}\in J_{q+1}(T)\mid {\xi}_q\in R_q, D{\xi}_{q+1}\in T^*\otimes R_q\}$ and we just need to use the two following formulas showing how the Spencer operator acts on the various brackets (See [GB1] [deform] for more details):  \\
\[  i(\zeta)D\{{\xi}_{q+1},{\eta}_{q+1}\}=\{i(\zeta)D{\xi}_{q+1},{\eta}_q\}+\{{\xi}_q,i(\zeta)D{\eta}_{q+1}\} ,\hspace {4mm} \forall \zeta \in T \]  
\[  i(\zeta)D[{\xi}_{q+1},{\eta}_{q+1}]=[i(\zeta)D{\xi}_{q+1},{\eta}_q]+[{\xi}_q,i(\zeta)D{\eta}_{q+1}]+i(L({\eta}_1)\zeta)D{\xi}_{q+1}-i(L({\xi}_1)\zeta)D{\eta}_{q+1}   \]
because the right member of the second formula is a section of $R_q$ whenever ${\xi}_{q+1},{\eta}_{q+1}\in R_{q+1}$. The first formula may be used when $R_q$ is formally integrable. \\
\hspace*{12cm}   Q.E.D.  \\

\noindent
{\bf EXAMPLE 2.B.6}: With $n=1, q=3, X=\mathbb{R}$ and evident notations, the components of $[{\xi}_3,{\eta}_3]$ at order zero, one and two are defined by the totally unusual successive formulas:\\
\[    [\xi,\eta]=\xi{\partial}_x\eta-\eta{\partial}_x\xi     \]
\[    ([{\xi}_1,{\eta}_1])_x=\xi{\partial}_x{\eta}_x-\eta{\partial}_x{\xi}_x    \]
\[    ([{\xi}_2,{\eta}_2])_{xx}={\xi}_x{\eta}_{xx}-{\eta}_x{\xi}_{xx}+\xi{\partial}_x{\eta}_{xx}-\eta{\partial}_x{\xi}_{xx}   \]
\[    ([{\xi}_3,{\eta}_3])_{xxx}=2{\xi}_x{\eta}_{xxx}-2{\eta}_x{\xi}_{xxx}+\xi {\partial}_x{\eta}_{xxx}-\eta{\partial}_x{\xi}_{xxx}   \]
For affine transformations, ${\xi}_{xx}=0,{\eta}_{xx}=0\Rightarrow ([{\xi}_2,{\eta}_2])_{xx}=0$ and thus $[R_2,R_2]\subset R_2$.\\
For projective transformations, ${\xi}_{xxx}=0,{\eta}_{xxx}=0 \Rightarrow ([{\xi}_3,{\eta}_3])_{xxx}=0$ and thus $[R_3,R_3]\subset R_3$.  \\

The next definition will generalize the definition of the {\it classical Lie derivative}:    \\
\[  {\cal{L}}(\xi)\omega =(i(\xi)d +di(\xi))\omega=\frac{d}{dt}j_q(exp\hspace{2mm} t\xi)^{-1}(\omega){\mid}_{t=0}.  \]

\noindent
{\bf DEFINITION 2.B.7}: We say that a vector bundle $F$ is {\it associated} with $R_q$ if there exists a first order differential operator $L({\xi}_q):F \rightarrow  F$ called {\it formal Lie derivative} and such that:  \\
1)  $L({\xi}_q+{\eta}_q)=L({\xi}_q)+L({\eta}_q)     \hspace{2cm} \forall {\xi}_q,{\eta}_q\in R_q$.  \\
2)  $L(f{\xi}_q)=fL({\xi}_q)  \hspace{2cm}   \forall {\xi}_q\in R_q, \forall f\in C^{\infty}(X)$.   \\
3)  $[L({\xi}_q),L({\eta}_q)]=L({\xi}_q)\circ L({\eta}_q)-L({\eta}_q)\circ L({\xi}_q)=L([{\xi}_q,{\eta}_q])  \hspace{2cm}  \forall {\xi}_q,{\eta}_q\in R_q$.   \\
4)  $ L({\xi}_q)(f\eta)=fL({\xi}_q)\eta + ({\xi}.f)\eta \hspace{15mm} \forall {\xi}_q\in R_q,\forall f\in C^{\infty}(X), \forall \eta \in F$. \\

As a byproduct, if $E$ and $F$ are associated with $R_q$, we may set on $E\otimes F$:Ê\\
\[    L({\xi}_q)(\eta \otimes \zeta)=L({\xi}_q)\eta\otimes \zeta+\eta\otimes L({\xi}_q)\zeta \hspace{2cm} \forall {\xi}_q\in R_q,\forall \eta \in E, \forall \zeta \in F   \] 
If $\Theta \subset T$ denotes the solutions of $R_q$, then we may set ${\cal{L}}(\xi)=L(j_q(\xi)), \forall \xi\in \Theta$ but no explicit computation can be done when $\Theta$ is infinite dimensional.  \\

\noindent
{\bf PROPOSITION 2.B.8}: $J_q(T)$ is associated with $J_{q+1}(T)$ if we define:  \\
\[     L({\xi}_{q+1}){\eta}_q=\{{\xi}_{q+1},{\eta}_{q+1}\}+i(\xi)D{\eta}_{q+1}=[{\xi}_q,{\eta}_q]+i(\eta)D{\xi}_{q+1}        \]
and thus $R_q$ is associated with $R_{q+1}$.  \

\noindent
{\bf Proof}: It is easy to check the properties 1, 2, 4 and it only remains to prove property 3 as follows.\\
\[  \begin{array}{rcl}
[L({\xi}_{q+1}),L({\eta}_{q+1})]{\zeta}_q & = & L({\xi}_{q+1})(\{{\eta}_{q+1},{\zeta}_{q+1}\}+i(\eta)D{\zeta}_{q+1})-L({\eta}_{q+1})(\{{\xi}_{q+1},{\zeta}_{q+1}\}+i(\xi)D{\zeta}_{q+1})   \\
       &  =  & \{{\xi}_{q+1},\{{\eta}_{q+2},{\zeta}_{q+2}\}\}-\{{\eta}_{q+1},\{{\xi}_{q+2},{\zeta}_{q+2}\}\}  \\
   &   & +\{{\xi}_{q+1},i(\eta)D{\zeta}_{q+2}\}-\{{\eta}_{q+1},i(\xi)D{\zeta}_{q+2}\}   \\
   &   & +i(\xi)D\{{\eta}_{q+2},{\zeta}_{q+2}\}-i(\eta)D\{{\xi}_{q+2},{\zeta}_{q+2}\}   \\
   &   & +i(\xi)D(i(\eta)D{\zeta}_{q+2})-i(\eta)D(i(\xi)D{\zeta}_{q+2})    \\
   &= & \{\{{\xi}_{q+2},{\eta}_{q+2}\},{\zeta}_{q+1}\}+\{i(\xi)D{\eta}_{q+2},{\zeta}_{q+1}\}-\{i(\eta)D{\xi}_{q+2},{\zeta}_{q+1}\}      \\
    &   &+i([\xi,\eta])D{\zeta}_{q+1}   \\
   &= & \{[{\xi}_{q+1},{\eta}_{q+1}],{\zeta}_{q+1}\}+i([\xi,\eta])D{\zeta}_{q+1}
   \end{array}   \]
by using successively the Jacobi identity for the algebraic bracket and the last proposition.Ê\\
\hspace*{12cm}   Q.E.D.   \\

\noindent
{\bf EXAMPLE 2.B.9}: $T$ and $T^*$ both with any tensor bundle are associated with $J_1(T)$. For $T$ we may define $L({\xi}_1)\eta=[\xi,\eta]+i(\eta)D{\xi}_1=\{{\xi}_1,j_1(\eta)\}$. We have ${\xi}^r{\partial}_r{\eta}^k-{\eta}^s{\partial}_s{\xi}^k+{\eta}^s({\partial}_s{\xi}^k-{\xi}^k_s)=-{\eta}^s{\xi}^k_s+{\xi}^r{\partial}_r{\eta}^k$ and the four properties of the formal Lie derivative can be checked directly. Of course, we find back ${\cal{L}}(\xi)\eta=[\xi,\eta], \forall \xi,\eta \in T$. We let the reader treat similarly the case of $T^*$. \\

 \noindent
{\bf THEOREM 2.B.10} : There is a {\it first nonlinear Spencer sequence}:\\
\[ 0\longrightarrow aut(X) \stackrel{j_{q+1}}{\longrightarrow} {\Pi}_{q+1}(X,X)\stackrel{\bar{D}}{\longrightarrow}T^*\otimes J_q(T)\stackrel{{\bar{D}}'}{\longrightarrow} {\wedge}^2T^*\otimes J_{q-1}(T)  \]
which is {\it locally exact} if $\Delta \neq 0$, with {\it restriction}:  \\
\[ 0\longrightarrow \Gamma \stackrel{j_{q+1}}{\longrightarrow} {\cal{R}}_{q+1} \stackrel{\bar{D}}{\longrightarrow} T^*\otimes R_q\stackrel{{\bar{D}}'}{\longrightarrow}{\wedge}^2T^*\otimes J_{q-1}(T)  \]

\noindent
{\bf Proof}: There is a canonical inclusion ${\Pi}_{q+1}\subset J_1({\Pi}_q)$ defined by $y^k_{\mu,i}=y^k_{\mu+1_i}$ and the composition $f^{-1}_{q+1}\circ j_1(f_q)$ is a well defined section of $J_1({\Pi}_q)$ over the section $f^{-1}_q\circ f_q=id_q$ of ${\Pi}_q$ like $id_{q+1}$. The difference ${\chi}_q=\bar{D}f_{q+1}=f^{-1}_{q+1}\circ j_1(f_q)-id_{q+1}$ is thus a section of $T^*\otimes V({\Pi}_q)$ over $id_q$ and thus of $T^*\otimes J_q(T)$. For $q=1$, setting ${\chi}_0=A-id\in T^*\otimes T$ and $g_1=f^{-1}_1$, we get:\\
\[ {\chi}^k_{,i}=g^k_l{\partial}_if^l-{\delta}^k_i=A^k_i-{\delta}^k_i,\hspace{5mm} {\chi}^k_{j,i}=g^k_l({\partial}_if^l_j-A^r_if^l_{rj})  \]
We shall prove later on the useful formula $ f^k_r{\chi}^r_{\mu,i}+...+f^k_{\mu+1_r}{\chi}^r_{,i}={\partial}_if^k_{\mu}-f^k_{\mu+1_i}$ allowing to determine ${\chi}_q$ inductively.\\
We have ${\bar{D}}'{\chi}_q(\xi,\eta)\equiv D{\chi}_q(\xi,\eta)-\{{\chi}_q(\xi),{\chi}_q(\eta)\}=0 $ and provide the only formulas that will be used later on and can be checked directly by the reader:\\
\[  {\partial}_i{\chi}^k_{,j}-{\partial}_j{\chi}^k_{,i}-{\chi}^k_{i,j}+{\chi}^k_{j,i}-({\chi}^r_{,i}{\chi}^k_{r,j}-{\chi}^r_{,j}{\chi}^k_{r,i})=0  \eqno{(1)}  \]
\[  {\partial}_i{\chi}^k_{l,j}-{\partial}_j{\chi}^k_{l,i}-{\chi}^k_{li,j}+{\chi}^k_{lj,i}-({\chi}^r_{,i}{\chi}^k_{lr,j}+{\chi}^r_{l,i}{\chi}^k_{r,j}-{\chi}^r_{l,j}{\chi}^k_{r,i}-{\chi}^r_{,j}{\chi}^k_{lr,i})=0  \eqno{(2)}  \]
In these sequences, the kernels are taken with respect to the zero section of the vector bundles involved. We finally notice that the condition $det(A)\neq 0$ amounts to $\Delta=det({\partial}_if^k)\neq 0$ because $det(f^k_i)\neq 0$ by assumption. One can prove by induction that the first nonlinear Spencer sequence is locally exact if $det(A)\neq 0$, that is any section of $T^*\otimes J_q(T)$ killed by ${\bar{D}}'$ is locally the image by $\bar{D}$ of a section of ${\Pi}_{q+1}$, contrary to its restriction (See [23], p 215 for more details and compare to [14], p 162, 195). Also, introducing the vector bundle $C_1=T^* \otimes R_q/ \delta g_{q+1}$, we have $ det(A)\neq 0 \Rightarrow \exists {\chi}_q={\tau}_q\circ A$ with ${\tau}^k_{\mu,i}={\chi}^k_{\mu,r}(A^{-1})^r_i=-g^k_lf^l_{\mu + 1_i}+ ...  $ and $\bar{D}$ induces a nonlinear operator ${\bar{D}}_1:R_q \rightarrow C_1$, a result that will be generalized later on in the linear framework. The brothers Cosserat were speaking about the {\it lagrangian field} ${\chi}_q$ and the {\it eulerian field} ${\tau}_q$ defined in ([7], \S 71, (70)+(71) $\leftrightarrow$ (72)+(73), p 190). {\it This is a subtle confusion} because the true {\it eulerian field} ${\sigma}_q= - {\bar{D}}f^{-1}_{q+1}$, obtained by exchanging source with target, cannot be expressed from ${\chi}_q$ by means of linear algebra (See [22], p 303 for more details).\\
\hspace*{12cm}  Q.E.D.  \\

\noindent
{\bf REMARK 2.B.11}: Rewriting the previous formulas with $A$ instead of ${\chi}_0$ we get:  \\
\[ {\partial}_iA^k_j-{\partial}_jA^k_i-A^r_i{\chi}^k_{r,j}+A^r_j{\chi}^k_{r,i}=0  \eqno{(1*)}  \]
\[ {\partial}_i{\chi}^k_{l,j}-{\partial}_j{\chi}^k_{l,i}-{\chi}^r_{l,i}{\chi}^k_{r,j}+{\chi}^r_{l,j}{\chi}^k_{r,i}-A^r_i{\chi}^k_{lr,j}+A^r_j{\chi}^k_{lr,i}=0 \eqno{(2*)} \]
When $q=1$ and $ g_2=0$, we find back {\it exactly} all the formulas presented by E. and F. Cosserat in [22], p 123 and [34]) (Compare to [14]). We finally notice that ${\chi}'_q= - {\chi}_q$ is a $R_q$-connection if and only if $A=0$, a result in contradiction with the use of connections in physics (Compare to [14], p 162, 195). However, when $A=0$, we have ${\chi}'_0(\xi)=\xi$ and thus (exercise):  \\
\[  \begin{array}{rcl}
 {\bar{D}}'{\chi}_{q+1} & = & (D{\chi}_{q+1})(\xi,\eta)-([{\chi}_q(\xi),{\chi}_q(\eta)] + i(\xi)D({\chi}_{q+1}(\eta)) - i(\eta)D({\chi}_{q+1}(\xi))) \\
    &   =   & - [{\chi}_q(\xi),{\chi}_q(\eta)]  - {\chi}_q([\xi,\eta])  \\
    &  =  & - {\kappa}'_q(\xi,\eta)
    \end{array}   \]
does not depend on the lift of ${\chi}_q$.  \\

\noindent
{\bf THEOREM 2.B.12}: In the case of a lie group of transformations, the nonlinear Spencer sequence is isomorphic to the nonlinear gauge sequence when $q$ is large enough and we have the following commutative diagram ([22], [23]):    \\
\[ \begin{array}{rccccc}
  &   X\times G & \rightarrow & T^*\otimes {\cal{G}} &\stackrel{MC}{ \rightarrow} &  {\wedge}^2T^*\otimes {\cal{G}}  \\
    & \downarrow  &   &  \downarrow  &  & \downarrow  \\
    0\rightarrow \Gamma \rightarrow & {\cal{R}}_q  & \stackrel{\bar{D}}{\rightarrow} & T^*\otimes R_q & \stackrel{{\bar{D}}'}{\rightarrow}  & {\wedge}^2T^* \otimes R_q  
    \end{array}     \]
{\it The action is essential in the Spencer sequence but disappears in the gauge sequence}.\\

Introducing now the Lie algebra ${\cal{G}}=T_e(G)$ and the Lie algebroid $R_q\subset J_q(T)$, namely the linearization of ${\cal{R}}_q$ at the $q$-jet of the identity $y=x$, we get the commutative and exact diagram:\\

\[  \begin{array}{rclccl}
0\rightarrow  &  \hspace{13mm} X   \hspace{6mm} \times & \cal{G} & = & \hspace{4mm} R_q & 
\rightarrow 0 \\
  &  \lambda=cst \uparrow \downarrow \uparrow \lambda(x) & & & j_q(\xi) \uparrow \downarrow \uparrow {\xi}_q &  \\
    &  \hspace{4mm} X &  & = & \hspace{4mm} X  & 
    \end{array}  \]
\noindent
where the upper isomorphism is described by ${\lambda}^{\tau}(x)\rightarrow {\xi}^k_{\mu}(x)={\lambda}^{\tau}(x){\partial}_{\mu}{\theta}^k_{\tau}(x)$ for $q$ large enough. Applying the Spencer operator, we finally obtain ${\partial}_i{\xi}^k_{\mu}(x)-{\xi}^k_{\mu +1_i}(x)={\partial}_i{\lambda}^{\tau}(x){\partial}_{\mu}{\theta}^k_{\tau}(x)$ and get:\\

\noindent
{\bf COROLLARY 2.B.13}: The linear Spencer sequence is isomorphic to the tensor product of the Poincar\'{e} sequence by $\cal{G}$ in the following commutative diagram:   \\
\[ \begin{array}{rlclcl}
  &   {\wedge}^0T^*\otimes {\cal{G}} & \longrightarrow & {\wedge}^1T^*\otimes {\cal{G}} & \longrightarrow &  {\wedge}^2T^*\otimes {\cal{G}}  \\
    & \hspace{9mm} \downarrow  &   & \hspace{9mm} \downarrow  &  &\hspace{9mm} \downarrow  \\
    0\rightarrow \Theta \stackrel{j_q}{\longrightarrow} & {\wedge}^0T^*\otimes R_q  & \stackrel{D}{\longrightarrow} & {\wedge}^1T^*\otimes R_q & \stackrel{D}{\longrightarrow}  & {\wedge}^2T^* \otimes R_q  
    \end{array}     \]
where the vertical isomorphisms are induced by the previous diagram.  \\

When $E$ is a vector bundle over $X$ and $D:J_{q+1}(E) \rightarrow T^*\otimes J_q(E)$ is the corresponding (linear) Spencer operator, we denote by $\delta :S_{q+1}T^*\otimes E \rightarrow T^*\otimes S_qT^*\otimes E$ the {\it Spencer map} induced (up to sign) by applying $D$ to the short exact sequence  $0 \rightarrow S_{q+1}T^*\otimes E \rightarrow J_{q+1}(E) \rightarrow J_q(E)\rightarrow 0$. We can extend the Spencer operator to an operator $D:{\wedge}^rT^*\otimes J_{q+1}(E) \rightarrow {\wedge}^{r+1}T^*\otimes J_q(E): \alpha\otimes {\xi}_{q+1}\rightarrow d\alpha \otimes {\xi}_q+(-1)^r\alpha\wedge D{\xi}_{q+1}$ and the corresponding Spencer map $\delta : {\wedge}^rT^*\otimes S_{q+1}T^*\otimes E \rightarrow {\wedge}^{r+1}T^*\otimes S_qT^*\otimes E$ is defined by  $(\delta \omega)^k_{\mu}=dx^i\wedge {\omega}^k_{\mu+1_i}$ . For any linear system $R_q\subset J_q(E)$, we may define the $r$-{\it prolongation} ${\rho}_r(R_q)=R_{q+r}=J_r(R_q)\cap J_{q+r}(E)\subset J_r(J_q(E))$ and the {\it symbol} ${\rho}_r(g_q)=g_{q+r}=R_{q+r}\cap S_{q+r}T^*\otimes E$ both with the restrictions $D:{\wedge}^rT^*\otimes R_{q+1} \rightarrow {\wedge}^{r+1}T^*\otimes R_q$ and $\delta : {\wedge}^rT^*\otimes g_{q+1} \rightarrow {\wedge}^{r+1}T^*\otimes g_q $. It is finally easy to verify that $D^2=0 \Rightarrow {\delta}^2=0$ 
([20], [25]). \\

\noindent
{\bf DEFINITION 2.B.14}: A system $R_q\subset J_q(E)$ is said to be {\it formally integrable} if all the equations of order $q+r$ that can be obtained from the system are obtained by differentiating $r$ times {\it only} the equations of order $q$ defining $R_q$ or, equivalently, if the maps $R_{q+r+1} \rightarrow R_{q+r}$ are epimorphisms $\forall r\geq 0$. Its symbol $g_q\subset S_qT^*\otimes E$ is said to be {\it finite type} if $g_{q+r}=0$ for $r$ large enough, $l$-{\it acyclic} if all the sequences $... \stackrel{\delta}{\rightarrow} {\wedge}^sT^* \otimes g_{q+r} \stackrel{\delta}{\rightarrow} ...  $ are exact $\forall r\geq 0, \forall s=1,...,l$ and {\it involutive} if it is $n$-acyclic. A finite type symbol $g_q$ is involutive if and only if $g_q=0$. Finally, a system is said to be {\it involutive} if it is formally integrable and if its symbol is involutive. Such crucial properties can now be checked by means of computer algebra techniques based on the {\it Janet/Goldschmidt/Spencer criterion} saying roughly that $R_q$ is formally integrable whenever $g_q$ is involutive or even $2$-acyclic and ${\pi}^{q+1}_q :R_{q+1} \rightarrow R_q$ is an epimorphism. Otherwise, one may start afresh with $R^{(1)}_q={\pi}^{q+1}_q(R_{q+1})\subset R_q$ and so on, till the criterion could be used ([12], [20], [25]).   \\

\noindent
{\bf EXAMPLE 2.B.15}: Linearizing the finite Lie equations of Example ..., we find a system $R_1\subset J_1(T)$ defined by the two first order equations $x^2{\xi}^1_1-x^1{\xi}^2_1+{\xi}^2=0, x^2{\xi}^1_2-x^1{\xi}^2_2-{\xi}^1=0$. In such an example, $g_1$ is involutive (exercise) but the system is not formally integrable because, using crossed derivatives, one can obtain the new first order equation ${\xi}^1_1+{\xi}^2_2=0$. The combined first order system, namely the projection $R^{(1)}_1$ of $R_2$ into $R_1$, is involutive with the same solutions.  \\

\noindent
{\bf REMARK 2.B.16}: The (first) {\it linear Spencer sequence} $0 \rightarrow E \stackrel{j_{q+1}}{\longrightarrow}J_{q+1}(E)\stackrel{D}{\longrightarrow}T^*\otimes J_q(E)\stackrel{D}{\longrightarrow} {\wedge}^2T^*\otimes J_{q-1}(E)$ and its restriction $0 \rightarrow \Theta \stackrel{j_{q+1}}{\longrightarrow} R_{q+1}\stackrel{D}{\longrightarrow} T^*\otimes R_q\stackrel{D}{\longrightarrow}{\wedge}^2T^*\otimes J_{q-1}(E)$ are not very useful in actual practice because the operator $D$ is not involutive and even not formally integrable. Indeed, from the first order equations ${\partial}_i{\xi}^k-{\xi}^k_i=0$, we obtain, by using crossed derivatives, the new first order equations ${\partial}_i{\xi}^k_j-{\partial}_j{\xi}^k_i=0$. \\

For any involutive system $R_q\subset J_q(E)$  the {\it Janet bundles} $F_r={\wedge}^rT^*\otimes J_q(E)/({\wedge}^rT^*\otimes R_q+\delta ({\wedge}^{r-1}T^*\otimes S_{q+1}T^*\otimes E))$ and {\it Spencer bundles} $C_r={\wedge}^rT^*\otimes R_q/\delta({\wedge}^{r-1}T^*\otimes g_{q+1})\subset C_r(E)={\wedge}^rT^*\otimes J_q(E)/\delta({\wedge}^{r-1}T^*\otimes S_{q+1}T^*\otimes E)$ are related by the short exact sequences $0 \rightarrow C_r \longrightarrow C_r(E) \stackrel{{\Phi}_r}{\longrightarrow} F_r \rightarrow 0 $ where the epimorphisms ${\Phi}_r$ for $ r=0, 1,...,n$ are induced by the epimorphism $\Phi={\Phi}_0:C_0(E)=J_q(E)\rightarrow J_q(E)/R_q=F_0$. In the commutative diagram below where {\it all the operators are induced by 
$D$}, the (second) {\it  linear Spencer sequence} is the kernel of the projection of the {\it linear hybrid sequence} onto the {\it linear Janet sequence}:  \\
 \[ \begin{array}{rcccccccccl}
 &&&&& 0 &&0&&0&  \\
 &&&&& \downarrow && \downarrow &   & \downarrow &  \\
  & 0& \rightarrow & \Theta &\stackrel{j_q}{\longrightarrow}&C_0 &\stackrel{D_1}{\longrightarrow}& C_1 &\stackrel{D_2}{\longrightarrow} ... \stackrel{D_n}{\longrightarrow}& C_n &\rightarrow 0 \\
  &&&&& \downarrow & & \downarrow &  &\downarrow &     \\
   & 0 & \rightarrow & E & \stackrel{j_q}{\longrightarrow} & C_0(E) & \stackrel{D_1}{\longrightarrow} & C_1(E) &\stackrel{D_2}{\longrightarrow}  ... \stackrel{D_n}{\longrightarrow} & C_n(E) &   \rightarrow 0 \\
   & & & \parallel && \hspace{5mm}\downarrow {\Phi}_0 & &\hspace{5mm} \downarrow {\Phi}_1 &  & \hspace{5mm}\downarrow {\Phi}_n & \\
   0 \rightarrow & \Theta &\rightarrow & E & \stackrel{\cal{D}}{\longrightarrow} & F_0  & \stackrel{{\cal{D}}_1}{\longrightarrow} & F_1 & \stackrel{{\cal{D}}_2}{\longrightarrow} ... \stackrel{{\cal{D}}_n}{\longrightarrow} & F_n & \rightarrow  0 \\
   &&&&& \downarrow & & \downarrow &    &\downarrow &   \\
   &&&&& 0 && 0 & &0 &  
   \end{array}     \]
In particular, if $E=T$ and $R_q\subset J_q(T)$ is a transitive involutive system of infinitesimal Lie equations, the Janet bundles are associated with 
$R_q$. If moreover $g_q=0$, then, whenever the dimension of the underlying Lie group is increasing, the dimensions of the Janet bundles are decreasing while the dimensions of the Spencer bundles are increasing by the same amount. We obtain therefore the following picture: \\
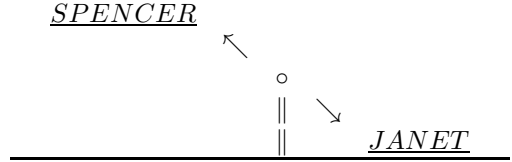
\[  \begin{array}{rcccccl}
   & \underline{SPENCER} & & & & & \\
   &  &\nwarrow & & & &  \\
   &  & & \circ & & &   \\
   &  & & \parallel & \searrow & &   \\
   &  & & \parallel & & \underline{JANET} &    \\
   \hline
\end{array}  \]
showing why, in some virtual sense, {\it Janet and Spencer are playing at see-saw} ([12]$\leftrightarrow$[33]). This picture will give the key for all the applications we shall present in the next section.\\

\noindent
{\bf EXAMPLE 2.B.17}: When $n=3$ and $E=X\times \mathbb{R}$, the second order system $R_2\subset J_2(E)$ defined by the three PD equations $y_{33}=0, y_{23}-y_{11}=0, y_{22}=0$ is trivially formally integrable because it is homogeneous but is not involutive because its symbol $g_2$ with $dim(g_2)=6-3=3$ is finite type with $dim(g_3)=1$ and $g_{4+r}=0, \forall r\geq 0$. Accordingly, we have $dim(R_2)=1+3+3=7$ while $dim(R_{3+r})=8=2^n, \forall r\geq 0$ ([16], p 79). We let the reader prove as an exercise of linear algebra that $g_3$ is $2$-acyclic by showing the exactness of the $\delta$-sequence $0 \rightarrow {\wedge}^2T^*\otimes g_3 \stackrel{\delta}{\rightarrow} {\wedge}^3T^*\otimes g_2 \rightarrow 0$ and we may consider the first prolongation $R_3\subset J_3(E)$ defined by the following $12$ PD equations:  \\
\[  \left\{  \begin{array}{l}
{\phi}^1 \equiv y_{333}=0  \\
{\phi}^2 \equiv y_{233}=0, {\phi}^3\equiv y_{223}=0, {\phi}^4\equiv y_{222}=0 \\
{\phi}^5 \equiv y_{133}=0, {\phi}^6\equiv y_{123}-y_{111}=0, {\phi}^7\equiv y_{122}=0, {\phi}^8\equiv y_{113}=0, {\phi}^9\equiv y_{112}=0 \\
{\phi}^{10}\equiv y_{33}=0, {\phi}^{11}\equiv y_{23}-y_{11}=0, {\phi}^{12}\equiv y_{22}=0
\end{array}
\right. \fbox{$\begin{array}{lll}
1 & 2 & 3 \\
1 & 2 & \bullet \\
1 & \bullet & \bullet \\
\bullet & \bullet & \bullet
\end{array}$}  \]
{\it In this particular situation}, that is when $g_3$ is already $2$-acyclic though NOT involutive, it is known that the generating {\it compatibility conditions} (CC) are first order (See [23], p 120) and described by the following $21$ PD equations:  \\
\[  \left\{  \begin{array}{l}
{\psi}^1\equiv d_3{\phi}^2-d_2{\phi}^1=0,...,{\psi}^5\equiv d_3{\phi}^6-d_1{\phi}^2+d_1{\phi}^8,..., {\psi}^8\equiv d_3{\phi}^9-d_2{\phi}^8=0 \\
{\psi}^9\equiv d_2{\phi}^5-d_1{\phi}_2=0, {\Psi}^{10}\equiv d_2{\phi}^6-d_1{\phi}^3+d_1{\phi}^9=0,..., {\psi}^{12}\equiv d_2{\phi}^9-d_1{\phi}^7=0\\
{\psi}^{13}\equiv d_3{\phi}^{10}-{\phi}^1=0, {\psi}^{14}\equiv d_3{\phi}^{11}-{\phi}^2, {\Psi}^{15}\equiv d_3{\phi}^{12}-{\phi}^3=0  \\
{\psi}^{16}\equiv d_2{\Phi}^{10}-{\Phi}^2=0, {\psi}^{17}\equiv d_2{\phi}^{11}-{\Phi}^3+{\Phi}^5=0, {\psi}^{18}\equiv d_2{\phi}^{12}-{\Phi}^4=0 \\
{\psi}^{19}\equiv d_1{\phi}^{10}-{\phi}^5=0, {\psi}^{20}\equiv d_1{\phi}^{11}-{\phi}^6=0, {\psi}^{21}\equiv d_1{\phi}^{12}-{\phi}^7=0
\end{array}
\right. \fbox{$\begin{array}{lll}
1 & 2 & 3 \\
1 & 2 & \bullet \\
1 & 2 & 3 \\
1 & 2& \bullet \\
1 & \bullet & \bullet 
\end{array}$}  \]
Each dot is producing one CC {\it apart from one} as we may verify the relation:  \\
\[  d_3{\psi}^{12}-d_2{\psi}^8+d_1{\psi}^6\equiv d_{22}{\phi}^8-d_{11}{\phi}^3   \]
and check therefore the remaining $13-1=12$ first order CC:  \\
\[  \left\{  \begin{array}{l}
{\theta}^1\equiv d_3{\psi}^9-d_2{\psi}^4+d_1{\psi}^9=0,...,{\theta}^3\equiv d_3{\psi}^{11}-d_2{\psi}^6+d_1{\psi}^{11}=0   \\
{\theta}^4\equiv d_3{\psi}^{16}-d_2{\psi}^{13}+{\psi}^1=0,...,{\theta}^6\equiv d_3{\psi}^{18}-d_2{\psi}^{15}+{\psi}^3=0   \\
{\theta}^7\equiv d_3{\psi}^{19}-d_1{\psi}^{13}+{\psi}^4=0, ..., {\theta}^9\equiv d_3{\psi}^{21}-d_1{\psi}^{15}+{\psi}^6=0 \\
{\theta}^{10}\equiv d_2{\psi}^{19}-d_1{\psi}^{16}+{\psi}^9=0,...,{\theta}^{12}\equiv d_2{\psi}^{21}-d_1{\psi}^{18}+{\psi}^{11}=0
\end{array}
\right. \fbox{$\begin{array}{lll}
1 & 2 & 3 \\
1 & 2 & 3 \\
1 & 2 & 3 \\
1 & 2 & \bullet 
\end{array}$}  \] 
It is quite a pure chance that this system is involutive with the following $3$ first order CC:  \\
\[  \left\{  \begin{array}{l}
d_3{\theta}^{10}-d_2{\theta}^7+d_1{\theta}^4-{\theta}^1=0, ... ,d_3{\theta}^{12}-d_2{\theta}^9+d_1{\theta}^6-{\theta}^3=0
\end{array}
\right. \fbox{$\begin{array}{lll}
1 & 2 & 3
\end{array}$}  \] 
THE FOLLOWING ABSOLUTELY NONTRIVIAL POINT WILL BE CRUCIAL FOR UNDERSTANDING THE STRUCTURE OF THE CONFORMAL LIE EQUATIONS LATER ON. \\
Indeed, with $q=3$ and $g_4=0$, we can define the Spencer bundles to be $C_r={\wedge}^rT^*\otimes R_3$, construct {\it in any case} the Janet sequence for the trivially involutive operator $j_3$ and obtain the following contradictory diagram where $dim(F_3)=2$ instead of the awaited $3$ (!):  \\
   \[ \begin{array}{rccccccccccccl}
 &&&&& 0 &&0&&0&  &0  & \\
 &&&&& \downarrow && \downarrow && \downarrow &    & \downarrow &  \\
  & 0& \rightarrow& \Theta &\stackrel{j_3}{\rightarrow}&8 &\stackrel{D_1}{\rightarrow}& 24 &\stackrel{D_2}{\rightarrow} & 24 &\stackrel{D_3}{\rightarrow} & 8    &\rightarrow 0   \\
  &&&&& \downarrow & & \downarrow & & \downarrow & &\downarrow &   \\
   & 0 & \rightarrow & 1 & \stackrel{j_3}{\rightarrow} & 20 & \stackrel{D_1}{\rightarrow} & 45 &\stackrel{D_2}{\rightarrow} & 36 &\stackrel{D_3}{\rightarrow} & 10  &   \rightarrow 0 \\
   & & & \parallel && \hspace{5mm}\downarrow {\Phi}_0 & &\hspace{5mm} \downarrow {\Phi}_1 & & \hspace{5mm}\downarrow {\Phi}_2 &  & \hspace{5mm}\downarrow {\Phi}_3 &  \\
   0 \rightarrow & \Theta &\rightarrow & 1 & \stackrel{\cal{D}}{\rightarrow} & 12  & \stackrel{{\cal{D}}_1}{\rightarrow} & 21 & \stackrel{{\cal{D}}_2}{\rightarrow} & 12 & \stackrel{{\cal{D}}_3}{\rightarrow} &2 & \rightarrow  0   \\
   &&&&& \downarrow & & \downarrow & & \downarrow &   &\downarrow &     \\
   &&&&& 0 && 0 && 0 &&0 & 
   \end{array}     \]
The explanation needs difficult homological algebra even on this elementary example which could be nevertheless treated by means of computer algebra while using quite large matrices. Indeed, starting from the short exact sequence $0\rightarrow R_3 \rightarrow J_3(E) \stackrel{\Phi}{\longrightarrow} F_0 \rightarrow 0$ with fiber dimensions $0 \rightarrow 8 \rightarrow 20 \rightarrow 12 \rightarrow 0$ and using $3$ prolongations in order to "reach" $F_3$, we get the following jet sequence of vector bundles, in fact the same that should be produced by any symbolic package:  \\
\[   0 \rightarrow R_6 \rightarrow J_6(E) \rightarrow J_3(F_0) \rightarrow J_2(F_1) \rightarrow J_1(F_2) \rightarrow F_3 \rightarrow 0  \]
with respective fiber dimensions:  \\
\[   0 \rightarrow 8 \rightarrow 84 \rightarrow 240 \rightarrow 210 \rightarrow 48 \rightarrow dim(F_3) \rightarrow 0  \]
Accordingly, {\it if the sequence were exact}, using the Euler-Poincar\'{e} formula ([15], Lemma 2.2, p 206), we should get $dim(F_3)=48-210+240-84+8=2$, a result showing that {\it the sequence cannot be exact}. Knowing why it is not exact and what is the resulting cohomology needs the following diagram obtained by induction, where {\it all the rows are exact but perhaps the upper one}:  \\
\[  \begin{array}{rcccccccccl}
   & 0 & & 0  & & 0  &  & 0  & & 0  \\
   & \downarrow  &  & \downarrow  & & \downarrow  & & \downarrow   & & \downarrow &  \\
0 \rightarrow  & g_6 & \rightarrow  & S_6 T^* & \rightarrow & S_3T^*\otimes F_0 & \rightarrow & S_2T^*\otimes F_1  & \rightarrow &T^*\otimes F_2  \\
  &\hspace{2mm} \downarrow \delta &  &\hspace{2mm} \downarrow \delta & &\hspace{2mm} \downarrow \delta &  & \hspace{2mm}\downarrow\delta &  & \parallel \\
 0 \rightarrow &  T^*\otimes g_5  & \rightarrow& T^*\otimes S_5T^*& \rightarrow  & T^* \otimes S_2T^*\otimes F_0& \rightarrow &T^* \otimes T^*\otimes F_1&   \rightarrow & T^*\otimes F_2  \\
   &\hspace{2mm} \downarrow  \delta &   & \hspace{2mm} \downarrow \delta &  &\hspace{2mm} \downarrow \delta  &  & \hspace{2mm}\downarrow \delta &  & \downarrow   \\
 0 \rightarrow & {\wedge}^2T^*\otimes  g_4  &\rightarrow  & {\wedge}^2 T^*\otimes S_4 T^* & \rightarrow  &{\wedge}^2T^* \otimes T^* \otimes F_0  &\rightarrow& {\wedge}^2T^*\otimes F_1 & \rightarrow & 0 \\
   & \hspace{2mm} \downarrow \delta &  &\hspace{2mm}  \downarrow \delta &  &\hspace{2mm}  \downarrow\delta  & & \downarrow  \\
 0 \rightarrow & {\wedge}^3 T^*\otimes g_3 & \rightarrow & {\wedge}^3T ^* \otimes S_3 T^*  & \rightarrow &{\wedge}^3T^*\otimes F_0&  &0  \\
   & \downarrow &  & \downarrow &  & \downarrow     \\
   &  0  &    &  0 &    & 0 &   &  & 
  \end{array}   \]
As $g_4=g_5=g_6=0$ and $dim ({\wedge}^3T^*\otimes g_3)=dim(g_3)=1$, a chase using the standard {\it snake lemma} of homological algebra ([32], p 174) proves that the upper seqence is not exact at $S_2T^*\otimes F_1$ with cohomology of dimension $1$. Hence, the previous sequence is not exact at $J_2(F_1)$, that is with $dim(im(J_3(F_0)\rightarrow J_2(F_1)))=240-84+8=164$ while $dim(ker(J_2(F_1) \rightarrow J_1(F_2)))= 164+1=165$ and we have indeed $48-210+165=3$. \\
The explanation of this tricky situation is not easy to grasp by somebody not familiar with homological algebra. Indeed, let us apply the $\delta$-map inductively to the short exact sequence $0\rightarrow g_{q+r} \rightarrow S_{q+r}T^*\otimes E \rightarrow h_r \rightarrow 0$ and consider the right part of the diagram thus obtained where the middle row is exact (See [23], p 151,152 for more details): \\
\[  \begin{array}{rlcll}
       & \hspace{13mm} 0  &  & \hspace{9mm} 0  &  \\
       & \hspace{13mm}\downarrow &  & \hspace{9mm} \downarrow &   \\
... \rightarrow &{\wedge}^{n-1}T^*\otimes g_{q+1} & \stackrel{\delta}{\rightarrow} & {\wedge}^nT^*\otimes g_q & \rightarrow 0  \\
 & \hspace{13mm}\downarrow &  & \hspace{9mm} \downarrow &   \\
... \rightarrow &{\wedge}^{n-1}T^*\otimes S_{q+1}T^*\otimes E & \stackrel{\delta}{\rightarrow} & {\wedge}^nT^*\otimes S_q T^*\otimes E& \rightarrow 0  \\
   & \hspace{13mm}\downarrow &  & \hspace{9mm} \downarrow &   \\
... \rightarrow &{\wedge}^{n-1}T^*\otimes h_1 & \stackrel{\delta}{\rightarrow} & {\wedge}^nT^*\otimes F_0 &   \\
   & \hspace{13mm}\downarrow &  & &   \\
 & \hspace{13mm} 0  &  &  &  
\end{array}   \]
Cutting the diagram, we may consider the following quotient diagram:  \\
\[ \begin{array}{rlclccr}
  & \hspace{16mm}0   &  &\hspace{9mm} 0  &  &  &    \\
  &  \hspace{16mm} \downarrow & &  \hspace{9mm} \downarrow &  &     &   \\
0 \rightarrow & \delta ({\wedge}^{n-1}T^*\otimes g_{q+1})  & \rightarrow & {\wedge}^nT^*\otimes R_q & \rightarrow & C_n  & \rightarrow 0  \\
  &\hspace{16mm} \downarrow &  & \hspace{9mm} \downarrow &  & \downarrow      &   \\
0 \rightarrow & \delta ({\wedge}^{n-1}T^*\otimes S_{q+1}T^*\otimes E) & \rightarrow & {\wedge}^nT^*\otimes J_q(E) & \rightarrow & C_n(E)& \rightarrow 0  \\
   & \hspace{16mm}\downarrow & &\hspace{9mm} \downarrow &  & \downarrow     &   \\
 0\rightarrow & \delta ({\wedge}^{n-1}T^*\otimes h_1) &\rightarrow & {\wedge}^nT^*\otimes F_0 & \rightarrow & F_n  & \rightarrow 0  \\
    &\hspace{16mm} \downarrow & &  \hspace{9mm}\downarrow &  & \downarrow      &   \\
  &  \hspace{16mm} 0  &  &  \hspace{9mm}  0  &  & 0  &
  \end{array}  \]
 When $g_q$ is $n-1$-acyclic but NOT $n$-acyclic, then $h_1$ is NOT $n-1$-acyclic and a chase is showing that the left vertical column is not exact at the central vector bundle. Using again the snake lemma, there is no way to get an upper injective map in the right vertical column. In the present situation with $n=3$ and $q=3$, as $g_3$ is $2$-acyclic but NOT $3$-acyclic and $g_4=0$, we have indeed $dim(F_3)=3$ because $dim(\delta({\wedge}^2T^*\otimes h_1))=10-1=9$ in a coherent way with explicit computations.  Accordingly, the only correct diagram allowing to deal with exact sequences on the jet level is the following one where all the operators involved are involutive, the $n=4$ vertical sequences are short exact sequences and $1-27+60-46+12=8-24+24-8=0$ ([15], Lemma 2.2, p 206): \\
  \[ \begin{array}{rccccccccccccl}
 &&&&& 0 &&0&&0&  &0  & \\
 &&&&& \downarrow && \downarrow && \downarrow &    & \downarrow &  \\
  & 0& \rightarrow& \Theta &\stackrel{j_4}{\rightarrow}&8 &\stackrel{D_1}{\rightarrow}& 24 &\stackrel{D_2}{\rightarrow} & 24 &\stackrel{D_3}{\rightarrow} & 8    &\rightarrow 0   \\
  &&&&& \downarrow & & \downarrow & & \downarrow & &\downarrow &   \\
   & 0 & \rightarrow & 1 & \stackrel{j_4}{\rightarrow} & 35 & \stackrel{D_1}{\rightarrow} & 84 &\stackrel{D_2}{\rightarrow} & 70 &\stackrel{D_3}{\rightarrow} & 20  &   \rightarrow 0 \\
   & & & \parallel && \hspace{5mm}\downarrow {\Phi}_0 & &\hspace{5mm} \downarrow {\Phi}_1 & & \hspace{5mm}\downarrow {\Phi}_2 &  & \hspace{5mm}\downarrow {\Phi}_3 &  \\
   0 \rightarrow & \Theta &\rightarrow & 1 & \stackrel{\cal{D}}{\rightarrow} & 27  & \stackrel{{\cal{D}}_1}{\rightarrow} & 60 & \stackrel{{\cal{D}}_2}{\rightarrow} & 46 & \stackrel{{\cal{D}}_3}{\rightarrow} & 12 & \rightarrow  0   \\
   &&&&& \downarrow & & \downarrow & & \downarrow &   &\downarrow &     \\
   &&&&& 0 && 0 && 0 &&0 & 
   \end{array}     \]
As a byproduct, ONE MUST CONSTRUCT THE JANET AND SPENCER SEQUENCES FOR AN INVOLUTIVE SYSTEM IN ORDER TO CONNECT THEM CONVENIENTLY.\\

\noindent
{\bf 3) APPLICATIONS} \\

Looking back to the end of Section 2A, it remains to graft a variational procedure adapted to the results of Section 2B. Similarly, as a major result first discovered in specific cases by the brothers Cosserat in 1909 and by Weyl in 1916, we shall prove and apply the following key result:  \\

\noindent
THE PROCEDURE ONLY DEPENDS ON THE DUAL OF THE SPENCER OPERATOR. \\

In order to prove this result, if $f_{q+1},g_{q+1},h_{q+1} \in {\Pi}_{q+1}$ can be composed in such a way that $g'_{q+1}=g_{q+1}\circ f_{q+1} = f_{q+1}\circ h_{q+1}$, we get:\\
\[ \begin{array}{rcccl}
{\bar{D}}g'_{q+1}&=&f^{-1}_{q+1}\circ g^{-1}_{q+1}\circ j_1(g_q)\circ j_1(f_q)-id_{q+1}  &  =  & f^{-1}_{q+1}\circ {\bar{D}}g_{q+1}\circ j_1(f_q)+{\bar{D}}f_{q+1} \\
       & = & h^{-1}_{q+1}\circ f^{-1}_{q+1}\circ j_1(f_q)\circ j_1(h_q) - id_{q+1}   &  =  &  h^{-1}_{q+1} \circ {\bar{D}}f_{q+1} \circ j_1(h_q) + \bar{D} h_{q+1}
 \end{array}    \]
Using the local exactness of the first nonlinear Spencer sequence or ([23], p 219), we may state:  \\
 
\noindent
{\bf LEMMA 3.1}: For any section $f_{q+1}\in {\cal{R}}_{q+1}$, the {\it finite gauge transformation}:\\
\[   {\chi}_q \in T^*\otimes R_q  \longrightarrow  {\chi}'_q= f^{-1}_{q+1}\circ {\chi}_q\circ j_1(f_q)+{\bar{D}}f_{q+1} \in T^* \otimes R_q \]
exchanges the solutions of the {\it field equations} ${\bar{D}}'{\chi}_q=0$.  \\

\noindent
{\bf LEMMA 3.2}: Passing to the limit {\it over the source} with $h_{q+1}=id_{q+1}+t {\xi}_{q+1}+ ... $ for $t\rightarrow 0$, we get an  {\it infinitesimal gauge transformation} leading to the {\it infinitesimal variation}:  \\
\[        \delta {\chi}_q= D{\xi}_{q+1}+ L(j_1({\xi}_{q+1})){\chi}_q     \]
which {\it does not depend on the parametrization} of ${\chi}_q$. \\

\noindent
{\bf LEMMA 3.3}: Passing to the limit {\it over the target} with ${\chi}_q=\bar{D}f_{q+1}$ and $g_{q+1}=id_{q+1}+ t {\eta}_{q+1}+ ... $, we get the other {\it infinitesimal variation}:\\
\[         \delta {\chi}_q= f^{-1}_{q+1}\circ D{\eta}_{q+1}\circ j_1(f_q)     \]
which {\it depends on the parametrization} of ${\chi}_q$. \\

We obtain in particular:   \\
\[   \begin{array}{rcl}
\delta {\chi}^k_{,i}  &  =  &  ({\partial}_i{\xi}^k-{\xi}^k_i) + ({\xi}^r{\partial}_r {\chi}^k_{,i} + {\chi}^k_{,r}{\partial}_i{\xi}^r - {\chi}^r_{,i} {\xi}^k_r )  \\
 \delta {\chi}^r_{r,i}  & =  &  ({\partial}_i{\xi}^r_r-{\xi}^r_{ri})+({\xi}^r{\partial}_r{\chi}^s_{s,i}+{\chi}^s_{s,r}{\partial}_i{\xi}^r+{\chi}^s_{,i}{\xi}^r_{rs})  
 \end{array}    \]
a result showing the importance of the Spencer operator. In the case of the Killing system $R_1$ with $g_2=0$, {\it these variations are exactly the ones provided by the brothers Cosserat} in ([7], (49)+(50), p 124, with a printing mistake corrected on p 128), replacing a $3\times 3$ skewsymetric matrix by the corresponding vector in ${\mathbb{R}}^3$. \\

These two explicit general formulas of the lemma cannot be found somewhere else (The reader may compare them to the ones obtained in [14] by means of the so-called " diagonal " method that cannot be applied to the study of explicit examples). The following unusual difficult proposition generalizes well known variational techniques used in continuum mechanics and will be crucially used for applications:  \\

\noindent
{\bf PROPOSITION 3.4}: The same variation is obtained whenever ${\eta}_q=f_{q+1}({\xi}_q+{\chi}_q(\xi))$ with ${\chi}_q=\bar{D}f_{q+1}$, a transformation which only depends on $j_1(f_q)$ and is invertible if and only if $det(A)\neq 0$.\\

\noindent
{\bf Proof}: First of all, setting ${\bar{\xi}}_q={\xi}_q + {\chi}_q(\xi)$, we get $\bar{\xi}=A(\xi)$ for $q=0$, a transformation which is invertible if and only if $det(A)\neq 0$. In the nonlinear framework, we have to keep in mind that there is no need to vary the object $\omega$ which is given but only the need to vary the section $f_{q+1}$ as we already saw, using ${\eta}_q\in R_q$ {\it over the target} or ${\xi}_q\in R_q$ {\it over the source}. With ${\eta}_q=f_{q+1}({\xi}_q)$, we obtain for example: \\
\[ \begin{array}{rcccl}
  \delta f^k & = & {\eta}^k & = & f^k_r{\xi}^r \\
   \delta f^k_i & = & {\eta}^k_uf^u_i & = & f^k_r{\xi}^r_i+f^k_{ri}{\xi}^r \\
  \delta f^k_{ij} & = & {\eta}^k_{uv}f^u_if^v_j+{\eta}^k_uf^u_{ij} & = & f^k_r{\xi}^r_{ij} + f^k_{ri}{\xi}^r_j+f^k_{rj}{\xi}^r_i+f^k_{rij}{\xi}^r 
  \end{array}  \]
and so on. Introducing the formal derivatives $d_i$ for $i=1,...,n$, we have:  \\
\[  \delta f^k_{\mu}={\zeta}^k_{\mu}(f_q,{\eta}_q)= d_{\mu}{\eta}^k={\eta}^k_uf^u_{\mu} + ... = f^k_r{\xi}^r_{\mu} + ... + f^k_{\mu +1_r} {\xi}^r       \]
We shall denote by $\sharp({\eta}_q)={\zeta}^k_{\mu}(y_q,{\eta}_q)\frac{\partial}{\partial y^k_{\mu}}\in V({\cal{R}}_q) $ with ${\zeta}^k={\eta}^k$ the corresponding vertical vector field, namely:    \\
\[   \sharp({\eta}_q)= 0\frac{\partial}{\partial x^i}+{\eta}^k(y)\frac{\partial}{\partial y^k}+({\eta}^k_u(y)y^u_i)\frac{\partial}{\partial y^k_i}+({\eta}^k_{uv}(y)y^u_iy^v_j+{\eta}^k_u(y)y^u_{ij} )\frac{\partial}{\partial y^k_{ij}}+ ...  \]
However, the standard prolongation of an infinitesimal change of source coordinates described by the horizontal vector field $\xi$, obtained by replacing all the derivatives of $\xi$ by a section ${\xi}_q \in R_q$ over $\xi \in T$, is the vector field: \\
\[   \flat({\xi}_q)={\xi}^i(x)\frac{\partial}{\partial x^i}+ 0\frac{\partial}{\partial y^k} - (y^k_r{\xi}^r_i(x))\frac{\partial}{\partial y^k_i}-(y^k_r{\xi}^r_{ij}(x)+y^k_{rj}{\xi}^r_i(x)+y^k_{ri}{\xi}^r_j(x)) \frac{\partial}{\partial y^k_{ij}}+ ...                                                             \]
It can be proved that $[\flat({\xi})_q,\flat({\xi}'_q]=\flat([{\xi}_q,{\xi}'_q]), \forall {\xi}_q,{\xi}'_q\in R_q$ {\it over the source}, with a similar property for $\sharp(.)$ {\it over the target} ([23]). However, $\flat({\xi}_q)$ {\it is not a vertical vector field and cannot therefore be compared to} $\sharp({\eta}_q)$.The solution of this problem explains a strange comment made by Weyl in ([36], p 289 + (78), p 290) and which became a founding stone of classical gauge theory. Indeed, ${\xi}^r_r$ is {\it not} a scalar because ${\xi}^k_i$ is {\it not} a $2$-tensor. However, when $A=0$, then $-{\chi}_q$ is a $R_q$-connection and ${\bar{\xi}}^r_r={\xi}^r_r+{\chi}^r_{r,i}{\xi}^i$ is a true scalar that may be set equal to zero in order to obtain ${\xi}^r_r=-{\chi}^r_{r,i}{\xi}^i$, a fact explaining why the EM-potential is considered as a connection in quantum mechanics instead of using the second order jets ${\xi}^r_{ri}$ of the conformal system, with a {\it shift by one step in the physical interpretation of the Spencer sequence} (See [22] for more historical details).\\
The main idea is to consider the vertical vector field $T(f_q)(\xi) - \flat({\xi}_q)\in V({\cal{R}}_q)$ whenever $y_q=f_q(x)$. Passing to the limit $t\rightarrow 0$ in the formula $g_q\circ f_q=f_q\circ h_q$, we first get $g\circ f = f\circ h \Rightarrow f(x)+t\eta (f(x)) + ... = f(x + t \xi(x) + ... )$. Using the chain rule for derivatives and substituting jets, we get successively:    \\
\[ \delta f^k(x)={\eta}^k(f(x))={\xi}^r{\partial}_r f^k, \hspace{2mm}  \delta f^k_i={\xi}^r{\partial}_rf^k_i + f^k_r {\xi}^r_i,\hspace{2mm}  \delta f^k_{ij}={\xi}^r{\partial}_rf^k_{ij}+f^k_{rj}{\xi}^r_i + f^k_{ri}{\xi}^r_j + f^k_r{\xi}^r_{ij}           \]
and so on, replacing ${\xi}^rf^k_{\mu + 1_r}$ by ${\xi}^r{\partial}_rf^k_{\mu}$ in ${\eta}_q=f_{q+1}({\xi}_q)$ in order to obtain:  \\
\[   \delta f^k_{\mu} = {\eta}^k_rf^r_{\mu} + ... ={\xi}^i({\partial}_if^k_{\mu}-f^k_{\mu+1_i})+ f^k_{\mu +1_r}{\xi}^r +  ... +f^k_r{\xi}^r_{\mu}   \]
where the right member only depends on $j_1(f_q)$ when $\mid\mu\mid=q$. \\
Finally, we may write the symbolic formula $f_{q+1}({\chi}_q)=j_1(f_q)-f_{q+1}=Df_{q+1}\in T^*\otimes V({\cal{R}}_q)$ in the explicit form:\\
\[        f^k_r{\chi}^r_{\mu,i} + ... +f^k_{\mu +1_r}{\chi}^r_{,i} = {\partial}_if^k_{\mu}-f^k_{\mu +1_i}   \]
Substituting in the previous formula provides ${\eta}_q=f_{q+1}({\xi}_q+ {\chi}_q(\xi))$ and we just need to replace $q$ by $q+1$ in order to achieve the proof. Replacing in the previous variations and using all the formulas involving the Spencer operator and the algebraic bracket that have been already exhibited, we let the reader prove as an exercise that we have equivalently: \\
\[  \delta {\chi}_q=D{\bar{\xi}}_{q+1} -\{{\chi}_{q+1}(.),{\bar{\xi}}_{q+1}\}    \]
We obtain in particular:  \\
\[ \begin{array}{rcl}
 \delta {\chi}^k_{,i} & = & ({\partial}_i{\bar{\xi}}^k-{\bar{\xi}}^k_i) - ({\chi}^r_{,i}{\bar{\xi}}^k_r - {\chi}^k_{r,i}{\bar{\xi}}^r)  \\
 \delta {\chi}^k_{j,i} & = & ({\partial}_i{\bar{\xi}}^k_j - {\bar{\xi}}^k_{ij}) - ({\chi}^r_{j,i}{\bar{\xi}}^k_r+{\chi}^r_{,i}{\bar{\xi}}^k_{jr} - {\chi}^k_{rj,i} - {\chi}^k_{r,i}) 
 \end{array}  \] 
Checking directly the proposition is not evident even when $q=0$ as we have:  \\
\[  (\frac{\partial {\eta}^k}{\partial y^u}-{\eta}^k_u){\partial}_if^u = f^k_r [({\partial}_i{\bar{\xi}}^r-{\bar{\xi}}^r_i) - ({\chi}^s_{,i}{\bar{\xi}}^r_s - {\chi}^r_{s,i}{\bar{\xi}}^s)]    \]
but cannot be done by hand when $q\geq 1$.Ê \\
\hspace*{12cm}  Q.E.D.   \\

We recall that the linear Spencer sequence for a Lie group of transformations $G\times X\rightarrow X$, which {\it essentially} depends on the action because infinitesimal generators are needed, is locally isomorphic to the linear gauge sequence which does not depend on the action any longer as it is the tensor product of the Poincar\'{e} sequence by the Lie algebra ${\cal{G}}$ of $G$. Accordingly, the main idea will be to introduce and compare the three following Lie groups of transformations but other subgroups of the conformal group may be considered, like the {\it optical subgroup} which is a maximal subgroup with $10$ parameters, contrary to the Poincar\'{e} subgroup which is not maximal:ÊÊ\\

\noindent
$\bullet$ The {\it Poincare group} of transformations leading to the {\it Killing system} $R_2$:  \\
\[ {\Omega}_{ij}\equiv (L({\xi}_1)\omega)_{ij}\equiv {\omega}_{rj}(x){\xi}^r_i+{\omega}_{ir}(x){\xi}^r_j+{\xi}^r{\partial}_r{\omega}_{ij}(x)=0 \]
\[{\Gamma}^k_{ij}\equiv (L({\xi}_2)\gamma)^k_{ij}\equiv {\xi}^k_{ij}+{\gamma}^k_{rj}(x){\xi}^r_i+{\gamma}^k_{ir}(x){\xi}^r_j-{\gamma}^r_{ij}(x){\xi}^k_r+{\xi}^r{\partial}_r{\gamma}^k_{ij}(x)=0 \]
\noindent 
$\bullet$ The {\it Weyl group} of transformations leading to the system ${\tilde{R}}_2$:   \\
\[ (L({\xi}_1)\omega)_{ij}\equiv {\omega}_{rj}(x){\xi}^r_i+{\omega}_{ir}(x){\xi}^r_j+{\xi}^r{\partial}_r{\omega}_{ij}(x)=A(x){\omega}_{ij}(x) \]
\[{\Gamma}^k_{ij}\equiv (L({\xi}_2)\gamma)^k_{ij}\equiv {\xi}^k_{ij}+{\gamma}^k_{rj}(x){\xi}^r_i+{\gamma}^k_{ir}(x){\xi}^r_j-{\gamma}^r_{ij}(x){\xi}^k_r+{\xi}^r{\partial}_r{\gamma}^k_{ij}(x)=0 \]
\noindent
$\bullet$ The {\it conformal group} of transformations leading to the {\it conformal Killing system} ${\hat{R}}_2$: \\
\[ (L({\xi}_1)\omega)_{ij}\equiv {\omega}_{rj}(x){\xi}^r_i+{\omega}_{ir}(x){\xi}^r_j+{\xi}^r{\partial}_r{\omega}_{ij}(x)=A(x){\omega}_{ij}(x) \]
\[ \begin{array}{rcl}
(L({\xi}_2)\gamma)^k_{ij} & \equiv & {\xi}^k_{ij}+{\gamma}^k_{rj}(x){\xi}^r_i+{\gamma}^k_{ir}(x){\xi}^r_j-{\gamma}^r_{ij}(x){\xi}^k_r+{\xi}^r{\partial}_r{\gamma}^k_{ij}(x) \\
  & =  & {\delta}^k_iA_j(x)+{\delta}^k_jA_i(x)-{\omega}_{ij}(x){\omega}^{kr}(x)A_r(x) 
  \end{array}   \]
where one has to eliminate the arbitrary function $A(x)$ and $1$-form $A_i(x)dx^i$ for finding sections, replacing the {\it ordinary Lie derivative} ${\cal{L}}(\xi)$ by the {\it formal Lie derivative} $L({\xi}_q)$, that is replacing $j_q(\xi)$ by ${\xi}_q$ when needed. According to the structure of the above Medolaghi equations, it is important to notice that $\Omega=L({\xi}_1)\omega \in S_2T^*$ and that $\Gamma=L({\xi}_2)\gamma \in S_2T^*\otimes T$. Moreover, as another way to consider the Christoffel symbols, $(\delta, -\gamma)=({\delta}^k_i,-{\gamma}^k_{ij})$ is a $R_1$-connection and thus also a ${\hat{R}}_1$-connection because $R_1\subset {\tilde{R}}_1 = {\hat{R}}_1$.  \\

$\bullet$ We make a few comments on the relationship existing between these systems. \\

First of all, when $\omega=({\omega}_{ij}(x)={\omega}_{ji}(x))$ is a non-degenerate metric, the corresponding Christoffel symbols are $\gamma = ({\gamma}^k_{ij}(x)=\frac{1}{2}{\omega}^{kr}(x)({\partial}_i{\omega}_{rj}(x) + {\partial}_j{\omega}_{ri}(x) - {\partial}_r{\omega}_{ij}(x))={\gamma}^k_{ji}(x))$. We have the relations $R_1 \subset {\tilde{R}}_1={\hat{R}}_1$ and obtain therefore $R_2={\rho}_1(R_1)$, $ {\tilde{R}}_2\subset {\rho}_1({\tilde{R}}_1)$, ${\hat{R}}_2={\rho}_1({\hat{R}}_1)$, a result leading to the {\it strict inclusions} $R_2\subset{\tilde{R}}_2\subset{\hat{R}}_2$ with respective fiber dimensions $10<11<15$ when $n=4$ and $\omega$ is the Minkowski metric with signature $(1,1,1,-1)$. \\
Secondly, if we want to deal with geometric objects in both cases, we have to introduce the {\it symmetric tensor density} ${\hat{\omega}}_{ij}={\omega}_{ij}/{\mid det(\omega)\mid}^{1/n}$ and the second order object ${\hat{\gamma}}^k_{ij}={\gamma}^k_{ij}-\frac{1}{n}({\delta}^k_i{\gamma}^r_{rj}+{\delta}^k_j{\gamma}^r_{ri} - {\omega}_{ij}{\omega}^{ks}{\gamma}^r_{rs}$ such that $\mid det(\hat{\omega})\mid =1,{\hat{\gamma}}^r_{ri}=0$, in such a way that ${\cal{R}}_2=\{f_2\in {\Pi}_2\mid f^{-1}_1({\hat{\omega}}=\hat{\omega}, f^{-1}_2(\hat{\gamma})=\hat{\gamma}\} $. It follows that ${\hat{g}}_1$ is defined by the equations ${\omega}_{rj}{\xi}^r_j+{\omega}_{rj}{\xi}^r_i-\frac{2}{n}{\omega}_{ij}{\xi}^r_r=0$ while ${\hat{g}}_2$ is defined by the equations ${\xi}^k_{ij}=\frac{1}{n}({\delta}^k_i{\xi}^r_{rj}+{\delta}^k_j{\xi}^r_{ri}-{\omega}_{ij}{\omega}^{ks}{\xi}^r_{rs})=0 $ which only depend on $\omega$ and no longer on $\hat{\omega}$. Only the first of the three following technical lemmas is known ([21], p 624-628):   \\

\noindent
{\bf LEMMA 3.5}: ${\hat{g}}_1$ is finite type with ${\hat{g}}_3=0$ when $n\geq 3$.  \\

\noindent
{\bf Proof}: The symbol ${\hat{g}}_3$ is defined by the equations ${\xi}^k_{ijt}-\frac{1}{n}({\delta}^k_i{\xi}^r_{jrt}+{\delta}^k_j{\xi}^r_{irt}-{\omega}_{ij}{\omega}^{ks}{\xi}^r_{rst}) =0$. Summing on $k$ and $t$, we get ${\xi}^r_{rij}-\frac{1}{n}(2{\xi}^r_{rij} - {\omega}_{ij}{\omega}^{st}{\xi}^r_{rst})=0$. Multiplying by ${\omega}^{ij}$ and summing on $i$ and $j$, we get ${\omega}^{ij}{\xi}^r_{rij}-\frac{2}{n}{\omega}^{ij}{\xi}^r_{rij}+{\omega}^{ij}{\xi}^r_{rij}=0$, that is to say ${\omega}^{ij}{\xi}^r_{rij}=0$ whenever $n\geq 2$. Substituting, we obtain $(n-2){\xi}^r_{rij}=0$ and thus ${\xi}^r_{rij}=0$ when $n\geq 3$, a result finally leading to ${\xi}^k_{ijt}=0$ and thus ${\hat{g}}_3=0, \forall n\geq 3$. In this case, it is important to notice that the third order jets only vanish when $\gamma=0$ locally or, equivalently, when $\omega$ is locally constant, for example when $n=4$ and $\omega$ is the Minkowski metric of space-time. \\
\hspace*{12cm}  Q.E.D.  \\

\noindent
{\bf LEMMA 3.6}: ${\hat{g}}_2$ is $2$-acyclic when $n\geq 4$.  \\

\noindent
{\bf Proof}: As ${\hat{g}}_{3+r}=0, \forall r\geq 0$, we have only to prove the injectivity of the map $\delta$ in the sequence: \\
\[ 0 \rightarrow {\wedge}^2T^*\otimes {\hat{g}}_2 \stackrel{\delta}{\rightarrow} {\wedge}^3T^*\otimes {\hat{g}}_1 \]
and thus to solve the linear system: \\
\[  {\xi}^k_{i\alpha,\beta\gamma} +{\xi}^k_{i\beta,\gamma\alpha}+{\xi}^k_{i\gamma,\alpha\beta} =0   \]
Substituting, we get the alternate sum over the cycle, where $\delta$ is again the Kronecker symbol:  \\
\[    {\cal{C}}(\alpha\beta\gamma)({\delta}^k_i{\xi}^r_{r\alpha,\beta\gamma}+{\delta}^k_{\alpha}{\xi}^r_{ri,\beta\gamma}-{\omega}_{i\alpha}{\omega}^{ks}{\xi}^r_{rs,\beta\gamma})=0   \]
Summing on $k$ and $i$, we get:  \\
\[   {\cal{C}}(\alpha\beta\gamma){\xi}^r_{r\alpha,\beta\gamma}=0 \hspace{5mm} \Rightarrow \hspace{5mm} {\cal{C}}(\alpha\beta\gamma)({\delta}^k_{\alpha}{\xi}^r_{ri,\beta\gamma}-{\omega}_{i\alpha}{\omega}^{ks}{\xi}^r_{rs,\beta\gamma})=0  \]
that is to say:  \\
\[  {\delta}^k_{\alpha}{\xi}^r_{ri,\beta\gamma}+{\delta}^k_{\beta}{\xi}^r_{ri,\gamma\alpha}+{\delta}^k_{\gamma}{\xi}^r_{ri,\alpha\beta}
-{\omega}^{ks}({\omega}_{i\alpha}{\xi}^r_{rs,\beta\gamma}+{\omega}_{i\beta}{\xi}^r_{rs,\gamma\alpha}+{\omega}_{i\gamma}{\xi}^r_{rs,\alpha\beta})=0  \]
Summing now on $k$ and $\alpha$, we get:  \\
\[  (n-3){\xi}^r_{ri,\beta\gamma}-{\omega}^{st}({\omega}_{i\beta}{\xi}^r_{rs,\gamma t}+{\omega}_{i\gamma}{\xi}^r_{rs,t\beta}=0  \]
Multiplying by ${\omega}^{ij}$ and summing on $i$, we get:  \\
\[  (n-3){\omega}^{ij}{\xi}^r_{ri,\beta\gamma}-{\omega}^{st}({\delta}^j_{\beta}{\xi}^r_{rs,\gamma t}+{\delta}^j_{\gamma}{\xi}^r_{rs,t\beta})=0  \]
Summing on $j$ and $\beta$, we finally obtain: \\
\[  2(n-2){\omega}^{ij}{\xi}^r_{ri,j\gamma}=0 \hspace{5mm} \Rightarrow \hspace{5mm} {\xi}^r_{ri,\beta\gamma}=0\hspace{5mm} \Rightarrow \hspace{5mm} {\xi}^k_{ij,\beta\gamma}=0, \hspace{5mm}  \forall n\geq 4\]
Accordingly, the linear system has the only zero solution and ${\hat{g}}_2$ is thus $2$-acyclic $\forall n\geq 4$, a quite deep reason for which {\it space-time has formal properties that are not satisfied by space alone}.  \\
\hspace*{12cm}  Q.E.D.  \\

\noindent
{\bf LEMMA 3.7}: ${\hat{g}}_2$ is $3$-acyclic when $n\geq 5$.     \\

\noindent
{\bf Proof}: As ${\hat{g}}_{3+r}=0, \forall r\geq 0$, we have only to prove the injectivity of the map $\delta$ in the sequence:  \\
\[    0 \rightarrow {\wedge}^3T^*\otimes {\hat{g}}_2 \stackrel{\delta}{\longrightarrow} {\wedge}^4T^*\otimes T^*\otimes T   \]
and thus to solve the linear system:  \\
\[   {\xi}^k_{i\alpha,\beta \gamma\delta}-{\xi}^k_{i\beta,\gamma\delta\alpha}+{\xi}^k_{i\gamma,\delta\alpha\beta}-{\xi}^k_{i\delta,\alpha\beta\gamma} = 0  \]
Substituting, we get the alternate sum over the cycle (care to the Kronecker symbol $\delta$):  \\
\[  {\cal{C}}(\alpha\beta\gamma\delta)({\delta}^k_i{\xi}^r_{r\alpha,\beta\gamma\delta}+{\delta}^k_{\alpha}{\xi}^r_{ri,\beta\gamma\delta}-{\omega}_{i\alpha}{\omega}^{ks}{\xi}^r_{rs,\beta\gamma\delta})=0    \]
Contracting in $k$ and $i$ the previous formula, we get:  \\
\[  {\cal{C}}(\alpha\beta\gamma\delta){\xi}^r_{r\alpha,\beta\gamma\delta}=0\hspace{5mm} \Rightarrow \hspace{5mm}  {\cal{C}}(\alpha\beta\gamma\delta)({\delta}^k_{\alpha}{\xi}^r_{ri,\beta\gamma\delta}-{\omega}_{i\alpha}{\omega}^{ks}{\xi}^r_{rs,\beta\gamma\delta})=0          \]
Contracting now in $k$ and $\alpha$, we get:  \\
\[  n{\xi}^r_{ri,\beta\gamma\delta}-{\xi}^r_{ri,\gamma\delta\beta}+{\xi}^r_{ri, \delta\gamma\beta}-{\xi}^r_{ri,\delta\beta\gamma} -{\xi}^r_{ri,\beta\gamma\delta}+{\omega}^{st}({\omega}_{i\beta}{\xi}^r_{rs,\gamma\delta t}-{\omega}_{i\gamma}{\xi}^r_{rs,\delta t \beta}+{\omega}_{i \delta}{\xi}^r_{rs,t \beta\gamma})=0\]
and thus:  \\
\[  (n-4){\xi}^r_{ri,\beta\gamma\delta}+{\omega}^{st}({\omega}_{i\beta}{\xi}^r_{rs,\gamma\delta t}+{\omega}_{i\gamma}{\xi}^r_{rs,\delta \beta t}+{\omega}_{i \delta}{\xi}^r_{rs, \beta\gamma t})=0\]
that we may transform into:  \\
\[    (n-4){\omega}^{ij}{\xi}^r_{ri,\beta\gamma\delta}+{\omega}^{st}({\delta}^j_{\beta}{\xi}^r_{rs,\gamma\delta t}+{\delta}^j_{\gamma}{\xi}^r_{rs,\delta \beta t}+{\delta}^j_{ \delta}{\xi}^r_{rs, \beta\gamma t})=0\]
Contracting in $l$ and $\beta$, we finally obtain:  \\
\[  2(n-3){\omega}^{ij}{\xi}^r_{ri,j\gamma\delta}=0 \hspace{5mm}\Rightarrow \hspace{5mm} (n-4){\xi}^r_{ri,\beta\gamma\delta}=0  \]
and ${\hat{g}}_2$ is thus $2$-acyclic for $n\geq 5$.  \\
\hspace*{12cm} Q.E.D.  \\

It follows from these lemmas that {\it we are exactly in the same situation as the one met in the previous example}, with a shift by one in the order of the operators involved. We may thus choose $C_r={\wedge}^rT^*\otimes {\hat{R}}_3\simeq {\wedge}^rT^*\otimes {\hat{R}}_2$ in the Spencer sequence:  \\
\[    0 \longrightarrow {\hat{\Theta}} \stackrel{j_3}{\longrightarrow} 15 \stackrel{D_1}{\longrightarrow} 60 \stackrel{D_2}{\longrightarrow} 90 \stackrel{D_3}{\longrightarrow} 60 \stackrel{D_4}{\longrightarrow} 15 \longrightarrow  0 \]
Each operator $D_r$ is thus induced by the Spencer operator $D:{\wedge}^{r-1}T^*\otimes {\hat{R}}_3\rightarrow {\wedge}^rT^*\otimes {\hat{R}}_2$ and is therefore a first order operator with constant coefficients, {\it both with its formal adjoint}. For later computations, the first Spencer operator in the sequence $J_3(E)\stackrel{D}{\longrightarrow}T^*\otimes J_2(E)\stackrel{D}{\longrightarrow} {\wedge}^2T^*\otimes J_1(E)$ can be described by the following images: \\
\[ {\partial}_i{\xi}^k-{\xi}^k_i=X^k_{,i} , \hspace{3mm}  {\partial}_i{\xi}^k_j-{\xi}^k_{ij}=X^k_{j,i} , \hspace{3mm}   {\partial}_i{\xi}^k_{lj}-{\xi}^k_{lij}=X^k_{lj,i}   \] 
while the second Spencer operator leads to the {\it identities}:\\
 \[   {\partial}_iX^k_{,j}-{\partial}_jX^k_{,i}+X^k_{j,i}-X^k_{i,j}=0, \hspace{5mm}  {\partial}_iX^k_{l,j}-{\partial}_jX^k_{l,i}+X^k_{lj,i}-X^k_{li,j}=0  \]

Finally, if ${\cal{D}}:E \longrightarrow F$ is a linear differential operator of order $q$, its {\it formal adjoint} $ad({\cal{D}}):{\wedge}^nT^*\otimes F^* \longrightarrow {\wedge}^nT^*\otimes E^*$ is again a linear differential operator of the same order $q$ that can be constructed by contraction with a test row $n$-form and integration by parts as usual by means of the Stokes formula. According to well known properties of the adjoint procedure, if ${\cal{D}}_1$ generates the CC of ${\cal{D}}$, then we have $ad({\cal{D}})\circ ad({\cal{D}}_1)=ad ({\cal{D}}_1\circ {\cal{D}})=ad (0)=0$ and thus $ad({\cal{D}})$ is surely among the CC of $ad({\cal{D}}_1)$ but may not generate them in general. By duality, this remark is at the origin of the difficult concept of {\it extension modules} in homological algebra and its application to the theory of differential modules. It can be proved that such a property does not depend on the differential sequence used, that is one can study ${\cal{D}}_{r-1}$ and ${\cal{D}}_r$ in the Janet sequence $\forall r\geq 1$ with ${\cal{D}}={\cal{D}}_0$ or, equivalently $D_r$ and $D_{r+1}$ in the Spencer sequence, {\it a first highly nontrivial result} ([25], [32]). In the case of the previous systems, as the Poincar\'{e} sequence is self-adjoint up to sign because $ad(grad)= - div$ when $X={\mathbb{R}}^3$, it follows that $ad({\cal{D}}_{r-1})$ generates the CC of $ad({\cal{D}}_r)$ while $ad(D_r)$ generates the CC of $ad(D_{r+1})$, {\it a second highly nontrivial result} (See examples in [27]). \\ 
 
$\bullet$ We now make a few comments on the relationship existing between these groups.  \\

As a Lie pseudogroup, the Poincar\'{e} group is defined by the system ${\cal{R}}_1\subset {\Pi}_1$ with the $n(n+1)/2$ equations ${\omega}_{kl}(y)y^k_iy^l_j={\omega}_{ij}(x)$. After linearization, $({\delta}^k_i,-{\gamma}^k_{ij})$ is the only existing symmetric $R_1$-connection for the Killing system $R_1\subset J_1(T)$ but $\gamma$ may also be considered as a geometric object of order $2$ with well known transition laws. As $g_2=0$, ${\pi}^2_1:R_2 \rightarrow R_1$ is an epimorphism but $R_1$ is {\it not involutive} and $R_2$ is {\it involutive} whenever the non-degenerate metric $\omega$ has constant riemannian curvature ([10], [20]). In actual practice, $n=4$ and $\omega$ is the Minkowski metric in the local coordinates $(x^1,x^2, x^3, x^4=ct)$. The fact that the Poincar\'{e} group could have something to do with the Galil\'{e}e group through a kind of limiting deformation procedure with $1/c\rightarrow 0$ is not correct because of a few general results on the {\it normalizer} $\tilde{\Gamma}=N(\Gamma)$ of $\Gamma$ in $aut(X)$ which are not so well known as their study involves a quite delicate use of the Spencer $\delta$-cohomology that we explain now (See [29] for more details).  \\

   In 1953 the physicists E. Inon\"{u} and E.P. Wigner (1963 Nobel prize) introduced the concept of {\it contraction of a Lie algebra} by considering the composition law $(u,v)\rightarrow (u+v)/(1+(uv/c^2))$ for speeds in special relativity (Poincar\' {e} group) when $c$ is the speed of light, claiming that the limit $c\rightarrow \infty$ or $1/c\rightarrow 0$ should produce the composition law $(u,v)\rightarrow u+v $ used in classical mechanics (Galil\' {e}e group) ([11]). However, this result is not correct indeed as $1/c\rightarrow 0$ has no meaning independently of the choice of length and time units. Hence, one has to consider the dimensionless numbers $\bar{u}=u/c,\bar{v}=v/c$ in order to get $(\bar{u},\bar{v})\rightarrow (\bar{u}+\bar{v})/(1+\bar{u}\bar{v})$ with no longer any perturbation parameter involved ([18]). Nevertheless, this idea brought the birth of the theory of {\it deformation of algebraic structures}, culminating in the use of the Chevalley-Eilenberg cohomology of Lie algebras ([6], [29]) and one of the first applications of computer algebra in the seventies because a few counterexamples can only be found for Lie algebras of dimension $\geq 11$ and have thus more than $500$ structure constants. Finally, it must also be noticed that the main idea of general relativity is to deform the Minkowski metric $dx^2+dy^2+dz^2-c^2dt^2$ of space-time by means of the small dimensionless parameter $\phi /c^2$ where $\phi=GM/r$ is the gravitational potential at a distance $r$ of a central attractive mass $M$ with gravitational constant $G$.\\

It has been the clever discovery of Ernest Vessiot (1865-1952) in 1903 ([35]), {\it still not known or even acknowledged today after more than a century} (Compare MR0720863 (85m:12004) to MR954613 (90e:58166)), to associate a {\it natural bundle} $\cal{F}$ over $X$ with any Lie pseudogroup $\Gamma \subset aut(X)$, both with a section $\omega$ of $\cal{F}$ called {\it geometric object} or {\it structure} on $X$ as we now explain by introducing a copy $Y$ of $X$ and considering the trivial fiber manifold $X\times Y \rightarrow X$. For this purpose, Vessiot noticed that any {\it horizontal vector field} $\xi={\xi}^i(x)\frac{\partial}{\partial x^i}$ commutes with any {\it vertical vector field} $\eta={\eta}^k(y)\frac{\partial}{\partial y^k}$ on $X\times X$. Using the chain rule for derivatives up to order $q$ with $\bar{x}=x+t{\xi}(x)+...$ or $\bar{y}=y+t{\eta}(y)+...$ where $t $ is a small parameter, we may work out the respective prolongations at order $q$ on jet coordinates, obtaining therefore the same commutation property on ${\Pi}_q$. As $[\Theta,\Theta]\subset \Theta$, we may use the Frobenius theorem on the target in order to find generating {\it differential invariants} $\{ {\Phi}^{\tau}(y_q)\}$ such that ${\Phi}^{\tau}({\bar{y}}_q)={\Phi}^{\tau}(y_q)$ whenever $\bar{y}=g(y)\in \Gamma$ acting now on the target copy $Y$ of $X$. Accordingly, {\it prolongations of source transformations exchange the differential invariants between themselves}, that is any (local) transformation $\bar{x}=\varphi(x)$ can be lifted to a (local) transformation of the differential invariants between themselves of the form $u\rightarrow \lambda(u,j_q(\varphi)(x))$ allowing to introduce a {\it natural bundle} $\cal{F}$ over $X$ by patching changes of coordinates $\bar{x}=\varphi(x), \bar{u}=\lambda(u,j_q(\varphi)(x))$. A section $\omega$ of $\cal{F}$ is called a {\it geometric object} or {\it structure} on $X$ and transforms like ${\bar{\omega}}(f(x))=\lambda(\omega(x),j_q(f)(x))$ or simply $\bar{\omega}=j_q(f)(\omega)$ whenever $y=f(x)$ is a reversible map. This is a way to generalize vectors and tensors ($q=1$), connections ($q=2$) or even higher order objects. As a byproduct we have $\Gamma=\{f\in aut(X){\mid} {\Phi}_{\omega}(j_q(f))\equiv j_q(f)^{-1}(\omega)=\omega\}$ as a new way to write out the finite Lie equations of $\Gamma$ and we may say that $\Gamma$ {\it preserves} $\omega$. Replacing $j_q(f)$ by $f_q$, we also obtain ${\cal{R}}_q=\{f_q\in {\Pi}_q{\mid} f_q^{-1}(\omega)=\omega\}$. Coming back to the infinitesimal point of view and setting $f_t=exp(t\xi)\in aut(X), \forall \xi\in T$, we may define the {\it ordinary Lie derivative} with value in the vector bundle $F_0={\omega}^{-1}({\cal{F}}_0)$ over $X$, {\it pull back} by $\omega$ of the vector bundle ${\cal{F}}_0=V({\cal{F}})$ over $\cal{F}$, by the formula :\\
\[   {\cal{D}}\xi={\cal{D}}_{\omega}\xi={\cal{L}}(\xi)\omega=\frac{d}{dt}j_q(f_t)^{-1}(\omega){\mid}_{t=0} \Rightarrow \Theta=\{\xi\in T{\mid}{\cal{L}}(\xi)\omega=0\}      \]
We have $x\rightarrow \bar{x}=x+t\xi(x)+...\Rightarrow u^{\tau}\rightarrow {\bar{u}}^{\tau}=u^{\tau}+t{\partial}_{\mu}{\xi}^kL^{\tau\mu}_k(u)+...$ where $\mu=({\mu}_1,...,{\mu}_n)$ is a multi-index and we may write down the system of infinitesimal Lie equations in the {\it Medolaghi form}:\\
\[     {\Omega}^{\tau}\equiv ({\cal{L}}(\xi)\omega)^{\tau}\equiv -L^{\tau\mu}_k(\omega(x)){\partial}_{\mu}{\xi}^k+{\xi}^r{\partial}_r{\omega}^{\tau}(x)=0    \]
as a way to state the invariance of the section $\omega$ of ${\cal{F}}$. Finally, replacing $j_q(\xi)$ by a section ${\xi}_q\in J_q(T)$ over $\xi\in T$, we may define $R_q\subset J_q(T)$ {\it on sections} by the purely linear equations:ÊÊ\\
\[  {\Omega}^{\tau}\equiv (L({\xi}_q)\omega)^{\tau}\equiv - L^{\tau\mu}_k(\omega (x)){\xi}^k_{\mu}+ {\xi}^r{\partial}_r{\omega}^{\tau}(x)=0  \]
By analogy with "special" and "general" relativity, we shall call the given section {\it special} and any other arbitrary section {\it general}. The problem is now to study the formal properties of the linear system just obtained with coefficients only depending on $j_1(\omega)$. In particular, if any expression involving $\omega$ and its derivatives is a scalar object, it must reduce to a constant whenever $\Gamma$ is assumed to be transitive and thus cannot be defined by any zero order equation. \\

\noindent
{\bf EXAMPLE 3.8}: Coming back to the affine and projective examples already presented, we show that the Vessiot structure equations may even exist when $n=1$. For this, we notice that the only generating differential invariant $\Phi\equiv y_{xx}/y_x$ of the affine case transforms like $u=\bar{u}{\partial}_xf+({\partial}_{xx}f/{\partial}_xf)$ while the only generating differential invariant $\Psi\equiv (y_{xxx}/y_x)-\frac{3}{2}(y_{xx}/y_x)^2$ of the projective case transforms like $v=\bar{v}({\partial}_xf)^2+({\partial}_{xxx}f/{\partial}_xf)-\frac{3}{2}({\partial}_{xx}f/{\partial}_xf)^2$ when $\bar{x}=f(x)$. If now $\gamma$ is the geometric object of the affine group $y=ax+b$ and $0\neq \alpha=\alpha (x)dx \in T^*$ is a $1$-form, we consider the object $\omega=(\alpha,\gamma)$ and get at once {\it one first order and one second order general Medolaghi equations}:\\  
\[  {\cal{L}}(\xi)\alpha\equiv \alpha {\partial}_x\xi + \xi {\partial}_x\alpha =0, \hspace{1cm} {\cal{L}}(\xi)\gamma\equiv {\partial}_{xx}\xi+\gamma {\partial}_x\xi+ \xi {\partial}_x\gamma =0  \]
Differentiating the first equation and substituting the second, we get the zero order equation:  \\
\[  \xi (\alpha {\partial}_{xx}\alpha-2({\partial}_x\alpha)^2+\alpha \gamma {\partial}_x\alpha-{\alpha}^2{\partial}_x\gamma)=0\hspace{5mm} \Leftrightarrow \hspace{5mm} \xi {\partial}_x(\frac{{\partial}_x\alpha}{{\alpha}^2} - \frac{\gamma}{\alpha} )=0  \]
and the {\it Vessiot structure equation} ${\partial}_x\alpha-\gamma \alpha=c{\alpha}^2$ where $c$ is an arbitrary constant. With $\alpha=1, \gamma=0 \Rightarrow c=0$ we get the translation subgroup $y=x+b$ while, with $\alpha=1/x, \gamma=0 \Rightarrow c=-1$ we get the dilatation subgroup $y=ax$. Similarly, if $\nu$ is the geometric object of the projective group and we consider the new geometric object $\omega=(\gamma,\nu)$, we get at once {\it one second order and one third order general Medolaghi equations}:  \\
\[  {\cal{L}}(\xi)\gamma\equiv {\partial}_{xx}\xi+\gamma {\partial}_x\xi+ \xi {\partial}_x\gamma =0 , \hspace{1cm}  {\cal{L}}(\xi)\nu\equiv {\partial}_{xxx}\xi+2\nu{\partial}_x\xi + \xi {\partial}_x\nu = 0  \]
and the only {\it Vessiot structure equation} is ${\partial}_x\gamma-\frac{1}{2}{\gamma}^2-\nu=0$, without any structure constant.\\

\noindent
{\bf EXAMPLE 3.9}: ({\it Riemann structure}) If $\omega=({\omega}_{ij}={\omega}_{ji})\in S_2T^*$ is a metric on a manifold $X$ with $dim(X)=n$ such that $det(\omega)\neq 0$, the Lie pseudogroup of transformations preserving $\omega$ is $\Gamma=\{f\in aut(X){\mid}j_1(f)^{-1}(\omega)=\omega \}$ and is a Lie group with a maximum number of $n(n+1)/2$ parameters. A special metric could be the Euclidean metric when $n=1,2,3$ as in elasticity theory or the Minkowski metric when $n=4$ as in special relativity [18]. The {\it  first order general Medolaghi equations}:\\
\[ {\Omega}_{ij}\equiv ({\cal{L}}(\xi)\omega)_{ij}\equiv {\omega}_{rj}(x){\partial}_i{\xi}^r+{\omega}_{ir}(x){\partial}_j{\xi}^r+{\xi}^r{\partial}_r{\omega}_{ij}(x)=0 \]
are also called {\it classical Killing equations} for historical reasons. The main problem is that {\it this system is not involutive} unless we prolong it to order two by differentiating once the equations. For such a purpose, introducing ${\omega}^{-1}=({\omega}^{ij})$ as usual, we may define the {\it Christoffel symbols}:\\
\[ {\gamma}^k_{ij}(x)=\frac{1}{2}{\omega}^{kr}(x)({\partial}_i{\omega}_{rj}(x) +{\partial}_j  {\omega}_{ri}(x) -{\partial}_r{\omega}_{ij}(x))=
{\gamma}^k_{ji}(x) \]
This is a new geometric object of order $2$ providing the Levi-Civita isomorphism $j_1(\omega)=(\omega,\partial \omega)\simeq (\omega,\gamma)$ of affine bundles and allowing to obtain the {\it second order general Medolaghi equations}:\\
\[ {\Gamma}^k_{ij}\equiv ({\cal{L}}(\xi)\gamma)^k_{ij}\equiv {\partial}_{ij}{\xi}^k+{\gamma}^k_{rj}(x){\partial}_i{\xi}^r+{\gamma}^k_{ir}(x){\partial}_j{\xi}^r-{\gamma}^r_{ij}(x){\partial}_r{\xi}^k+{\xi}^r{\partial}_r{\gamma}^k_{ij}(x)=0   \]
Surprisingly, the following expression, called {\it Riemann tensor}:\\
\[ {\rho}^k_{lij}(x)\equiv {\partial}_i{\gamma}^k_{lj}(x)-{\partial}_j{\gamma}^k_{li}(x)+{\gamma}^r_{lj}(x){\gamma}^k_{ri}(x)-{\gamma}^r_{li}(x){\gamma}^k_{rj}(x)  \]
is still a first order geometric object and even a $4$-tensor with $n^2(n^2-1)/12$ independent components satisfying the purely algebraic relations :\\
 \[   {\rho}^k_{lij}+{\rho}^k_{ijl}+{\rho}^k_{jli}=0, \hspace{5mm} {\omega}_{rl}{\rho}^r_{kij}+{\omega}_{kr}{\rho}^r_{lij}=0  \]
Accordingly, the IC must express that the new first order equations $R^k_{lij}\equiv ({\cal{L}}(\xi)\rho)^k_{lij}=0$ are only linear combinations of the previous ones and we get the {\it Vessiot structure equations}:\\
    \[  {\rho}^k_{lij}(x)=c({\delta}^k_i{\omega}_{lj}(x)-{\delta}^k_j{\omega}_{li}(x))   \]
with the only {\it structure constant} $c$ describing the constant Riemannian curvature condition of Eisenhart ([10], [20], [22], [23]). One can proceed similarly for the {\it conformal Killing system} ${\cal{L}}(\xi)\omega=A(x)\omega$ and obtain that the {\it Weyl tensor} must vanish, without any structure constant ([20], p 132). Though this result, first found by the author of this paper as early as in 1978 ([20]) is still not acknowledged, {\it there is no conceptual difference at all between the unique structure constant c appearing in this example and the previous one}. Moreover, the structure constants have in general nothing to do with the structure constants of any Lie algebra. \\ 

More generally, any generating set $\{{\Phi}^{\tau}\}$ of differential invariants must satisfy quasi-linear CC of the symbolic form $v\equiv I(u_1) \equiv A(u)u_x+B(u)=0$ where $u_1=(u,u_x)$, allowing to define an affine subfibered manifold ${\cal{B}}_1\subset J_1({\cal{F}})$ over ${\cal{F}}$ and a natural bundle ${\cal{F}}_1=J_1({\cal{F}})/{\cal{B}}_1$ over ${\cal{F}}$ with local coordinates $(x,u,v)$. The {\it Vessiot structure equations} $I(u_1)=c(u)$ are defined by an equivariant section $c:{\cal{F}} \rightarrow {\cal{F}}_1:(x,u) \rightarrow (x,u,v=c(u))$ depending, as we just saw, on a finite number of constants (See [20] and [23] for details and other examples). {\it The form of the Vessiot structure equations is invariant under any change of local coordinates}. The following result, already known to Vessiot in 1903 ([35], p 445), is still ignored today. For this, let us consider two sections $\omega$ and $\bar{\omega}$ of ${\cal{F}}$ giving rise, through the corresponding Medolaghi equations, to the systems $R_q$ and ${\bar{R}}_q$. We define the equivalence relation:  \\

\noindent
{\bf DEFINITION 3.10}: $\bar{\omega} \sim \omega  \Leftrightarrow {\bar{R}}_q=R_q$.  \\

The following result is not evident at all ([20], [29]):  \\

\noindent
{\bf PROPOSITION 3.11}: $\bar{\omega}$ is obtained from $\omega$ by a Lie group of transformations acting on the fibers of ${\cal{F}}$, namely the {\it reciprocal} of the Lie group of transformations describing the natural structure of ${\cal{F}}$. These finite transformations of the form $\bar{u}=g(u,a)$ will be called {\it label transformations} and the number of parameters $a$ is $\leq  dim (J_q(T)/R_q)=dim (F_0)$. \\

\noindent
{\bf COROLLARY 3.12}: Any finite label transformation $\bar{u}=g(u,a)$ induces a finite transformation $\bar{c}=h(c,a)$ and we say that $\bar{\omega} \sim \omega  \rightarrow \bar{c}\sim c$. \\

\noindent
{\bf DEFINITION 3.13}: The {\it normalizer} ${\tilde{\Gamma}}=N(\Gamma)$ of $\Gamma$ in $aut(X)$ is the biggest Lie pseudogroup in which $\Gamma$ is {\it normal}, that is (roughly) $N(\Gamma)={\tilde{\Gamma}}=\{ {\tilde{f}}\in aut(X){\mid} {\tilde{f}}\circ f\circ {\tilde{f}}^{-1}\in \Gamma, \forall f\in \Gamma \}$ and we write $\Gamma \lhd N(\Gamma)\subset aut(X)$.  \\ 

Of course, $N(\Theta)=\{\eta\in T {\mid} [\xi,\eta]\subset \Theta, \forall \xi\in \Theta\}$ will play the part of a Lie algebra for $N(\Gamma)$ exactly like $\Theta$ did for $\Gamma$. However, $N(\Gamma)$ may have many components different from the connected component of the identity. For example, when $n=2$ and $\Gamma$ is defined by the system $\{y^1_2=0, y^2_1=0\}$, then $N(\Gamma)$ is defined by the system $\{y^1_2=0, y^2_1=0\}\cup \{ y^1_1=0, y¬2_2=0\}$ as it contains the permutation $y^1=x^2, y^2=x^1$. In actual practice, ${\tilde{\Gamma}}=\{{\tilde{f}}\in aut(X){\mid} \bar{\omega}=j_q({\tilde{f}})^{-1}(\omega)=g(\omega,a),h(c,a)=c\}$ is defined by the system ${\tilde{{\cal{R}}}}_{q+1}=\{ {\tilde{f}}_{q+1}\in {\Pi}_{q+1}{\mid}{\tilde{f}}_{q+1}(R_q)=
R_q\}$ with linearization ${\tilde{R}}_{q+1}=\{ {\tilde{\xi}}_{q+1}{\mid} L({\tilde{\xi}}_{q+1}){\eta}_q\in R_q, \forall {\eta}_q\in R_q\}$, that is to say $\{{\tilde{\xi}}_{q+1},{\eta}_{q+1}\}+i({\tilde{\xi}})D{\eta}_{q+1}\in R_q \Leftrightarrow \{{\tilde{\xi}}_{q+1},{\eta}_{q+1}\}\in R_q$. Accordingly, the system ${\tilde{R}}_{q+1}$ defining ${\tilde{\Theta}}=N(\Theta)$ can be obtained by purely algebraic techniques from the system defining $\Theta$. We have ([20], p 390; [21], p 715; [22], p 548;[29]):\\

\noindent
{\bf PROPOSITION 3.13}: If $R_q$ is formally integrable and $g_q$ is $2$-acyclic, then ${\tilde{R}}_{q+1}$ is formally integrable with ${\tilde{g}}_{q+1}=g_{q+1}$.  \\

\noindent
{\bf EXAMPLE 3.15}: In the previous Example with $\omega = (\alpha, \gamma)$, we obtain by substraction $\bar{\omega}\sim \omega \Leftrightarrow (\bar{\alpha}=a\alpha, \bar{\gamma}=\gamma + b\alpha) \Rightarrow \bar{c}=\frac{1}{a}c-\frac{b}{a}$. The condition $\bar{c}=c$ provides $b=(1-a)c$, that is to say $b=0$ if $c=0$ and $b=a-1$ if $c=-1$. Hence, in both cases the corresponding Lie pseudogroup is of codimension $1$ in its normalizer. Indeed, the normalizer of $y=ax$ is $y=ax^b$ while the normalizer of $y=x+b$ is $y=ax+b$ with different meanings for the constants $a$ and $b$. Similarly, in the case of the Riemann structure, we let the reader prove as an exercise that $\bar{\omega}\sim \omega \Leftrightarrow \bar{\omega}=a\omega \rightarrow \bar{c}=\frac{1}{a}c$ because $\bar{\gamma}=\gamma$. Accordingly, the corresponding Lie pseudogroup is of codimension zero in its normalizer if $c\neq 0$ and of codimension $1$ if $c=0$, a result explaining why the normalizer of the Poincar\'{e} group is the Weyl group, obtained by adding a {\it unique dilatation for space and time}, contrary to the Galil\'{e}e group which is of codimension $2$ in its normalizer, obtained by adding {\it separate dilatations for space and time} ([22], [29]). We invite the reader to treat similarly the examples provided in the first section in order to understand how tricky are the computations involved or to look at the example fully treated in ([21], p726).  \\

$\bullet$ We now study each group separately, in relation with applications. \\

\noindent
{\bf EXAMPLE 3.16}: (Poincar\'{e} group) Changing slightly the notations while restricting for simplicity the formulas to the plane with $n=2$ and local coordinates $(x^1,x^2)$ instead of space with $n=3$ and local coordinates $(x^1,x^2, x^3)$ or space-time with $n=4$ and local coordinates $(x^1,x^2, x^3, x^4=ct)$, we may copy the equations $(12)$ of ([7], p 14) and $(12')$ of ([7], p 19) side by side in the following way:  \\
\[   \left\{ \begin{array}{lcl}
\frac{dF}{ds}   & = & L \\
& &   \\
\frac{dG}{ds}  & = & M  \\
& &  \\
\frac{dH}{ds} + \frac{dx^1}{ds}G-\frac{dx^2}{ds}F &=& N 
\end{array}  \right. 
\Longleftrightarrow
\left\{  \begin{array}{lcl}
 \frac{dF}{ds}  & = &L    \\
 &  &  \\
 \frac{dG}{ds}  & =& M  \\
&  &  \\
\frac{d}{ds}(H + x^1G - x^2F) & = & N + x^1 M - x^2L  
\end{array}
\right.
  \]
We notice that the left members of the equations on the right hand side are only made by the derivative of an expression with respect to the curvilinear abscissa $s$ along the curve considered in the plane with local coordinates $(x^1(s),x^2(s))$. Equivalently, we may use the linear transformations $(F,G,H) \rightarrow (F'=F, G'=G, H'=H+x^1G-x^2F)$ and $(L,M,N) \rightarrow (L'=L, M'=M, N'=N+x^1M-x^2L)$ {\it with the same underlying $3\times 3$ matrix of full rank $3$}, namely: \\
\[    (F',G',H')=(F,G,H) \left( \begin{array}{ccr}
1 & 0 & -x^2 \\
0 & 1 & x^1 \\
0 & 0 & 1\hspace{2mm} 
\end{array} 
\right), \hspace{2mm}
(L',M',N')=(L,M,N) \left( \begin{array}{ccr}
1 & 0 & -x^2 \\
0 & 1 & x^1 \\
0 & 0 & 1 \hspace{2mm}
\end{array}
\right)    \]
but {\it this result is not intrinsic at all and just looks like a pure coincidence}. It is important to notice that, while these formulas have been exhibited in the study of the (static) {\it deformation theory of a line} (Chapter II of [7], p 14 and 19), similar formulas also exist in the study of the (static) {\it deformation theory of a surface} (Chapter III of [7], p 76 and 91) and in the study of the (static) {\it deformation theory of a medium} (Chapter IV of [7], p 137 and 140). We shall not insist on these points which have already been treated elsewhere with full details ([22], [28]) and that we have recovered in this paper by means of other methods, but invite the reader to look at the amount of calculations provided by the brothers E. and F. Cosserat. However, in order to establish a link between this example and the use of the Spencer operator, we now consider the Killing system for $n=2$ and the euclidean metric. The dual of the Spencer operator is provided by the integration by parts of the contraction $2$-form while raising or lowering the indices by means of the metric, : \\
\[ {\sigma}^{1,1}({\partial}_1{\xi}_1- {\xi}_{1,1})+{\sigma}^{2,1}({\partial}_1{\xi}_2-{\xi}_{2,1})+{\sigma}^{1,2}({\partial}_2{\xi}_1-{\xi}_{1,2})+{\sigma}^{2,2}({\partial}_2{\xi}_2-{\xi}_{2,2})+{\mu}^{12,r}({\partial}_r{\xi}_{1,2}-0)     \]
provided that ${\xi}_{1,1}=0, {\xi}_{1,2}+{\xi}_{2,1}=0,{\xi}_{2,2}=0$. Integrating by parts, the factors of ${\xi}_1,{\xi}_2,{\xi}_{1,2}$ furnishes (up to sign) the Cosserat equations where, of course, ${\sigma}^{1,2}$ may be different from ${\sigma}^{2,1}$:  \\
\[  {\partial}_r{\sigma}^{i,r}=f^i, \hspace{3mm}  {\partial}_r{\mu}^{12,r}+{\sigma}^{1,2}-{\sigma}^{2,1}=m^{12}  \]  
In arbitrary dimension, one should get similarly ([7], p 137 for $n=3$, p 167 for $n=4$) (See [28], Remark 7.1, p 25 for more details):\\
\[   {\partial}_r{\sigma}^{i,r}=f^i, \hspace{3mm}  {\partial}_r{\mu}^{ij,r}+{\sigma}^{i,j}-{\sigma}^{j,i}=m^{ij}, \hspace{3mm}  \forall i<j    \]
As a byproduct, we obtain: \\
\[  {\partial}_r({\mu}^{ij,r}+x^j{\sigma}^{i,r}-x^i{\sigma}^{j,r})=m^{ij}+x^jf^i-x^if^j , \hspace{4mm}  \forall 1\leq i<j\leq n \]
that is {\it exactly} the equation used in continuum mechanics in order to study the torsor equilibrium bringing the symmetry of the stress tensor when $\mu=0$ and $m=0$, where the left member is the Stokes formula applied to the total surface density of momentum while the right member is the total volume density of momentum (See [22], [24], [27] for more details and compare to [34]).\\
With the infinitesimal generators ${\theta}_1={\partial}_1,{\theta}_2={\partial}_2, { \theta}_3=x^1{\partial}_2-x^2{\partial}_1$, setting ${\xi}^k_{\mu}={\lambda}^{\tau}{\partial}_{\mu}{\theta}^k_{\tau}$ with ${\lambda}^{\tau}(x)\in {\cal{G}}$, we have:  \\
\[  {\xi}^1={\lambda}^1-x^2{\lambda}^3, {\xi}^2={\lambda}^2+x^1{\lambda}^3, {\xi}^2_1=-{\xi}^1_2={\lambda}^3  \]
and we find back {\it exactly} the $3\times 3$ matrix with full rank already exhibited. \\

In fact, our purpose is quite different now though it is also based on the combined use of group theory and the Spencer operator. The idea is to notice that the brothers are {\it always} dealing with the same group of rigid motions because the lines, surfaces or media they consider are all supposed to be in the same $3$-dimensional background/surrounding space which is acted on by the group of rigid motions, namely a group with $6$ parameters ($3$ {\it translations} + $3$ {\it rotations}). In 1909 it should have been strictly impossible for the two brothers to extend their approach to bigger groups, in particular to include the only additional {\it dilatation} of the Weyl group that will provide the virial theorem and, {\it a fortiori}, the {\it elations} of the conformal group considered later on by H.Weyl. In order to emphasize the reason for using Lie equations, we now provide the explicit form of the $n$ finite elations and their infinitesimal counterpart, namely:\\
\[  y=\frac{x-x^2b}{1-2(bx)+b^2x^2}  \Rightarrow  {\theta}_s= - \frac{1}{2} x^2 {\delta}^r_s{\partial}_r+{\omega}_{st}x^tx^r{\partial}_r   \Rightarrow \  
{\partial}_r{\theta}^r_s=n{\omega}_{st}x^t, \hspace{2mm}   \forall 1\leq r,s,t \leq n\]
where the underlying metric is used for the scalar products $x^2,bx,b^2$ involved.  \\

\noindent
{\bf EXAMPLE 3.17}:({\it Weyl group}) We may rewrite the infinitesimal Lie equations in the local form:  \\
\[    {\xi}_{i,j}+{\xi}_{j,i} - \frac{2}{n}{\omega}_{ij}{\xi}^r_r=0, \hspace{4mm}  {\xi}^k_{ij}=0, \hspace{4mm}  0\leq i,j,k \leq n   \]
The contraction form should be complemented by the terms $ {\nu}^r({\partial}_r{\xi}^1_1-0)$ and the integration by parts provides the additional dual equation ${\partial}_r{\nu}^r + {\sigma}^r_r= v  $ because ${\xi}^1_1=...={\xi}^n_n=\frac{1}{n}{\xi}^r_r$ (Compare to (74) in [36], p 288). As a byproduct, we get:  \\
\[   {\partial}_r({\nu}^r + x^i{\sigma}^r_i) = v+ x^if_i   \]
that is {\it exactly} the virial equation already presented for the symmetric stress used in continuum mechanics and gas dynamics, where the left member 
is the Stokes formula applied to the total surface density of virial while the right member is the total volume density of virial. \\
Introducing the additional infinitesimal generator ${\theta}_4=x^i{\partial}_i$, we now get:\\
\[  {\xi}^1={\lambda}^1-x^2{\lambda}^3+x^1{\lambda}^4, {\xi}^2={\lambda}^2+x^1{\lambda}^3+x^2{\lambda}^4, {\xi}^2_1={\lambda}^3, {\xi}^1_1={\xi}^2_2={\lambda}^4\]
and obtain the $4\times 4$ matrix of rank $4$: \\
\[   \left( \begin{array}{ccrr}
1 & 0 & -x^2 & x^1 \\
0 & 1 & x^1& x^2 \\
0 & 0 & 1\hspace{2mm} & 0 \hspace{2mm} \\
0 & 0 & 0 \hspace{2mm}& 1 \hspace{2mm}
\end{array} 
\right)  \]
describing the linear map $X\times {\tilde{\cal{G}}} \rightarrow {\tilde{R}}_1$ with $dim({\tilde{R}}_1)=4$ when $X={\mathbb{R}}^2$. \\

\noindent
{\bf EXAMPLE 3.18}: ({\it Conformal group}) First of all, we explain the confusion done by Weyl in ([36]) between {\it natural bundles} and {\it jet bundles}, recalling that both bundles have only been introduced fifty years later but that the formula of Weyl that we shall consider has been one of the key ingredients of gauge theory, also fifty years later but for a quite different reason (See ([22]), Chapter 5, p 321-343 for historical comments). Indeed, considering $\omega$ and $\gamma$ as geometric objects, we obtain at once the formulas:\\
\[  {\bar{\omega}}_{ij}=a(x){\omega}_{ij} \hspace{5mm}  \Rightarrow  \hspace{5mm}  {\bar{\gamma}}^r_{ri}={\gamma}^r_{ri} + \frac{1}{2a}{\partial}_ia  \]
Though looking like the key formula ($69$)in ([36], p 286), this transformation is quite different because the sign is not coherent and the second object has nothing to do with a $1$-form. Moreover, if we use $n=2$ and set ${\cal{L}}(\xi)\omega=A\omega$ for the standard euclidean metric, we should have $({\partial}_{11}+{\partial}_{22})A=0$, contrary to the assumption that $A$ is arbitrary which is {\it only} agreeing with the jet formulas:\\
\[ L({\xi}_1)\omega=A\omega\hspace{3mm} \Rightarrow \hspace{3mm} 2({\xi}^r_r+{\gamma}^r_{ri}{\xi}^i)=nA, \hspace{3mm}(L({\xi}_2\gamma)^r_{ri}=nA_i , \hspace{3mm} \forall {\xi}_2\in {\hat{R}}_2\]
Now, if we make a change of coordinates $\bar{x}=\varphi (x)$ on a function $a\in {\wedge}^0T^*$, we get:\\
\[  \bar{a}(\varphi(x))=a(x) \hspace{4mm} \Rightarrow \hspace{4mm} \frac{\partial\bar{a}}{\partial {\bar{x}}^j}\frac{\partial {\varphi}^j}{\partial x^i}=\frac{\partial a}{\partial x^i}  \]
We obtain therefore an isomorphism $J_1({\wedge}^0T^*)\simeq {\wedge}^0T^*{\times}_XT^*$, a result leading to the following commutative diagram: \\
\[  \begin{array}{rcccccl}
0 \longrightarrow & R_2 & \longrightarrow & {\hat{R}}_2 & \longrightarrow  & J_1({\wedge}^0T^*) & \longrightarrow 0  \\
  & \hspace{3mm}\downarrow D  &  & \hspace{3mm}  \downarrow D  & &\hspace{3mm} \downarrow D  &   \\
0 \longrightarrow & T^* \otimes R_1 & \longrightarrow & T^* \otimes {\hat{R}}_1 & \longrightarrow  & T^* & \longrightarrow 0     
\end{array}    \]
where the rows are exact by counting the dimensions. The operator on the right is $D:(\frac{1}{2}A,A_i) \longrightarrow (\frac{1}{2}{\partial}_iA-A_i)$ and is induced by the central Spencer operator, a result that could not have been even imagined by Weyl and followers.\\
Though striking it may loo like, this result provides a good transition towards the conformal origin of electromagnetism. The nonlinear aspect has been already presented in ([22], [23], [31]) and we restrict our study to the linear framework. A first problem to solve is to construct vector bundles from the various components of the image of $D_1$. For this purpose, let us introduce $(B^k_{l,i}=X^k_{l,i}+{\gamma}^k_{ls}X^s_{,i}) \in T^*\otimes T^*\otimes T$ with $(B^r_{r,i}=B_i)\in T^*$ and $(B^k_{lj,i}=X^k_{lj,i}+{\gamma}^k_{sj}X^s_{l,i}+{\gamma}^k_{ls}X^s_{j,i}-{\gamma}^s_{lj}X^k_{s,i}+X^r_{,i}{\partial}_r{\gamma}^k_{lj}) \in T^*\otimes S_2T^*\otimes T$ with $(B^r_{ri,j}-B^r_{rj,i}=F_{ij})\in {\wedge}^2T^*$,
We obtain from the relations ${\partial}_i{\gamma}^r_{rj}={\partial}_j{\gamma}^r_{ri}$ and the previous identities:  \\
\[ \begin{array}{rcl}
F_{ij}=B^r_{ri,j}-B^r_{rj,i} & = & X^r_{ri,j}-X^r_{rj,i}+{\gamma}^r_{rs}X^s_{i,j}-{\gamma}^r_{rs}X^s_{j,i}+X^r_{,j}{\partial}_r{\gamma}^s_{si}-X^r_{,i}{\partial}_r{\gamma}^s_{sj}  \\
  &  =  & {\partial}_iX^r_{r,j}-{\partial}_jX^r_{r,i}+{\gamma}^r_{rs}(X^s_{i,j}-X^s_{j,i})+X^r_{,j}{\partial}_i{\gamma}^s_{sr}-X^r_{,i}{\partial}_j{\gamma}^s_{sr} \\
    &  =  & {\partial}_i(X^r_{r,j}+{\gamma}^r_{rs}X^s_{,j})-{\partial}_j(X^r_{r,i}+{\gamma}^r_{rs}X^s_{s,i})  \\
      &  =  &  {\partial}_iB_j-{\partial}_jB_i
      \end{array}   \]
Now, we have:\\
\[ \begin{array}{rcl}
 B_i & =  & ({\partial}_i{\xi}^r_r - {\xi}^r_{ri})+{\gamma}^r_{rs}({\partial}_i{\xi}^s - {\xi}^s_i)\\
    &  =  &{\partial}_i{\xi}^r_r + {\gamma}^r_{rs}{\partial}_i{\xi}^s+
{\xi}^s {\partial}_s{\gamma}^r_{ri} - nA_i \\
  &  =  &{\partial}_i({\xi}^r_r + {\gamma}^r_{rs}{\xi}^s) - nA_i \\
    &  =  &n (\frac{1}{2}{\partial}_iA - A_i) 
\end{array}   \]  
and we finally get $F_{ij}=n({\partial}_jA_i-{\partial}_iA_j)$, a result fully solving the dream of Weyl. Of course, when $n=4$ and $\omega$ is the Minkowski metric, then we have $\gamma=0$ in actual practice and the previous formulas become particularly simple. \\     
As ${\tilde{C}}_r={\wedge}^rT^*\otimes {\tilde{R}}_2\subset {\wedge}^rT^*\otimes {\hat{R}}_2={\hat{C}}_r$ and ${\hat{R}}_2/{\tilde{R}}_2\simeq T^*$, we get ${\hat{C}}_r/C_r\simeq {\wedge}^rT^*\otimes T^*$ and the conformal Spencer sequence projects onto the sequence $T^*\rightarrow T^*\otimes T^*\rightarrow {\wedge}^2T^*\otimes T^*\rightarrow ...$. Finally, {\it the Spencer sequence projects with a shift by one step onto the Poincar\'{e} sequence} $T^*\stackrel{d}{\rightarrow} {\wedge}^2T^* \stackrel{d}{\rightarrow} {\wedge}^3T^*\rightarrow ... $ obtained by applying the Spencer map $\delta$, because these two sequences are only made by first order involutive operators and the successive projections can therefore be constructed inductively. The short exact sequence $0\rightarrow S_2T^*\stackrel{\delta}{\rightarrow} T^*\otimes T^* \stackrel{\delta}{\rightarrow} {\wedge}^2T^*\rightarrow 0$ has already been used in ([22], [23], [28], [30]) for exhibiting the Ricci tensor and the above result brings for the first time a conformal link between electromagnetism and gravitation by using second order jets. \\
As for duality, using standard notations, we have the possible additional terms:\\
\[      ... + {\cal{J}}^i ({\partial}_i{\xi}^r_r - {\xi}^r_{ri}) + {\sum}_{i<j} {\cal{F}}^{ij}({\partial}_i{\xi}^r_{rj} - {\partial}_j{\xi}^r_{ri})  \]
Developping the sum, we get:  \\
\[... - {\cal{J}}^1{\xi}^r_{r1}- ... - {\cal{J}}^n{\xi}^r_{rn} +{\cal{F}}^{12}({\partial}_1{\xi}^r_{r2}-{\partial}_2{\xi}^r_{r1})+ ... + {\cal{F}}^{1n}({\partial}_1{\xi}^r_{rn}-{\partial}_n{\xi}^r_{rn}) + ... \] 
and the integration by part provides therefore the equations: \\
\[            {\partial}_j{\cal{F}}^{ij}-{\cal{J}}^i= w^i                                                \]

The link with the virial theorem is provided by the formulas ${\xi}^r_{ri}= n {\xi}^1_{1i}= ...=n {\xi}^n_{ni}$ and ${\nu}^i=n{\cal{J}}^i$. Accordingly, when the second members of the inductions equations vanish, we have ${\partial}_i{\cal{J}}^i={\partial}_{ij}{\cal{F}}^{ij}=0 \Rightarrow {\sigma}^r_r=0$ in a coherent way with a well known property of the so-called impulsion-energy tensor in electromagnetism.     \\

We have therefore obtained the following crucial theorem and striking corollary:   \\

\noindent
{\bf THEOREM 3.19}: In the Spencer sequence for the conformal Killing system with $n=4$, the {\it field equations} $C_1 \stackrel{D_2}{\rightarrow} C_2$ projects onto the first set of Maxwell equations ${\wedge}^2T^* \stackrel{d}{\longrightarrow} {\wedge}^3T^*$ and the correponding {\it potential parametrization} $C_0 \stackrel{D_1}{\longrightarrow} C_1$ projects onto the usual parametrization ${\wedge}^1T^* \stackrel{d}{\longrightarrow} {\wedge}^2T^*$ by the electromagnetic $4$-potential, according to the following commutative diagram:\\
\[  \begin{array}{ccccc}
  {\wedge}^0T^*\otimes {\hat{R}}_2 & \stackrel{D_1}{\longrightarrow} & {\wedge}^1T^*\otimes {\hat{R}}_2 & \stackrel{D_2}{\longrightarrow} & {\wedge}^2T^*\otimes {\hat{R}}_2 \\
 \downarrow & & \downarrow  & & \downarrow  \\
 T^* & \stackrel{d}{\longrightarrow} & {\wedge}^2T^* & \stackrel{d}{\longrightarrow} & {\wedge}^3T^* \\ 
 \downarrow &  & \downarrow  &  & \downarrow  \\
   0  &  &  0  &  &  0  
 \end{array}  \]
 By duality, the second set of Maxwell equations ${\wedge}^{n-2}T^* \stackrel{ad(d)}{\longrightarrow} {\wedge}^{n-1}T^*$ is induced by the {\it induction equations} ${\wedge}^nT^*\otimes C^*_1 \stackrel{ad(D_1)}{\longrightarrow} {\wedge}^nT^*\otimes C^*_0$ and a similar property holds for the corresponding {\it pseudo-potential parametrizations}, according to the following commutative diagram: \\
\[  \begin{array}{ccccc}
 0 & & 0 & &0 \\
  \downarrow & & \downarrow & & \downarrow \\
{\wedge}^{n-1}T^* & \stackrel{ad(d)}{\longleftarrow} & {\wedge}^{n-2}T^* & \stackrel{ad(d)}{\longleftarrow} & {\wedge}^{n-3}T^*  \\
 \downarrow & & \downarrow & & \downarrow  \\
{\wedge}^nT^*\otimes {\hat{R}}^*_2 & \stackrel{ad(D_1)}{\longleftarrow} &{\wedge}^{n-1}T^*\otimes {\hat{R}}^*_2 & \stackrel{ad(D_2)}{\longleftarrow} & {\wedge}^{n-2}T^*\otimes {\hat{R}}^*_2  
\end{array}  \]
Similar comments can be done for the Clausius and Cosserat equations because $R_2 \subset {\tilde{R}}_2 \subset {\hat{R}}_2$.  \\

\noindent
{\bf REMARK 3.20}: The key formulas $(76)$ in ([36], p 289) are based on a confusion between the Janet and Spencer sequences. Indeed, using only components of the Spencer operator when ${\xi}_2\in {\hat{R}}_2$, we have on one side in $S_2T^*$:  \\
\[   \begin{array}{rcl}
{\omega}_{rj}X^r_{,i}+{\omega}_{ir}X^r_{,j} & = & {\omega}_{rj}({\partial}_i{\xi}^r - {\xi}^r_i) + {\omega}_{ir}({\partial}_j{\xi}^r - {\xi}^r_j)   \\
& =& ({\omega}_{rj}{\partial}_i{\xi}^r+{\omega}_{ir}{\partial}_j{\xi}^r+{\xi}^r{\partial}_r{\omega}_{ij})-\frac{2}{n}({\xi}^r_r + {\gamma}^r_{rs}{\xi}^s){\omega}_{ij} \end{array}  \]
Similarly, we have on the other side in $T^*$:  \\
\[  \begin{array}{rcl}
X^r_{r,i}+{\gamma}^r_{rs}X^s_i & = & ({\partial}_i{\xi}^r_r - {\xi}^r_{ri})+{\gamma}^r_{rs}({\partial}_i{\xi}^s-{\xi}^s_i)  \\
                                                 & = &  {\partial}_i({\xi}^r_r+{\gamma}^r_{rs}{\xi}^s) - ({\xi}^r_{ri}+{\gamma}^r_{rs}{\xi}^s_i+{\xi}^s{\partial}_s{\gamma}^r_{rri}) \\
                                                &  =  &  n(\frac{1}{2}{\partial}_iA - A_i)
\end{array}   \]
In the nonlinear framework, using the variational formulas already establihed at the beginning of this section, it follows that:  \\
\[  ({\alpha}_i={\chi}^r_{r,i} + {\gamma}^r_{rs}{\chi}^s_{,i}) \in T^*\hspace{5mm} \Rightarrow \hspace{5mm}  \delta {\alpha}_i=n(\frac{1}{2}{\partial}_iA - A_i) + ({\xi}^r{\partial}_r{\alpha}_i + {\alpha}_s{\partial}_i{\xi}^s-n{\chi}^s_{,i}A_s)  \]
a formula that {\it cannot} be fully identified with $(76)$ because $A^k_i\neq 0$. \\

\noindent
{\bf REMARK 3.21}: As another confusion, we revisit a basic result of classical gauge theory. First of all, we recall that the classical lagrangian of a free particle of mass $m$ and charge $e$ in an EM field $F=dA$ is $L(t,x,\dot{x})=\frac{1}{2}m{\omega}_{ij}{\dot{x}}^i{\dot{x}}^j + e{\dot{x}}^iA_i$ and let the reader check that the corresponding Euler-Lagrange equations are $\frac{d\vec{v}}{dt}=\vec{v}\wedge \vec{B}$ where $\vec{v}=({\dot{x}}^i)$ and the right member is the Lorentz force. Introducing the {\it momentum} $p_i=m{\omega}_{ir}{\dot{x}}^r+eA_i$ and substituting, we obtain easily the hamiltonian $H(t,x,p)=\frac{1}{2m}{\omega}^{ij}(p_i-eA_i)(p_j-eA_j)$ which is obtained from the hamiltonian $H=\frac{1}{2m}{\omega}^{ij}p_ip_j$ of the uncharged particle by the transformation $p_i \longrightarrow p_i-eA_i$. A main idea of gauge theory has been to transform $p_i$ to $- i\hbar {\partial}_i$ according to the correspondence principle of quantum mechanics and $A_i$ to a connection. Quite contrary to this point of view, we have:\\

\noindent
{\bf COROLLARY 3.22}: The transformation $p_i \longrightarrow p_i - A_i$ is provided by the dual action of the second order jets of the conformal system and has a purely group-theoretical origin.  \\

\noindent
{\bf Proof}: The nonlinear system ${\hat{\cal{R}}}_2\subset {\Pi}_2$ of conformal finite Lie equations is:  \\
\[  {\omega}_{kl}f^k_if^l_j=a(x){\omega}_{ij} \Rightarrow g^k_l(f^l_{ij} + {\gamma}^l_{rs}f^r_if^s_j)= {\gamma}^k_{ij} + {\delta}^k_ia_j(x)+{\delta}^k_ja_i(x) - {\omega}_{ij}{\omega}^{kr}a_r(x) \]
Using the formulas of the last proposition describing the transformation ${\eta}_1=f_2 ({\xi}_1)$ with $f_2 \in {\hat{\cal{R}}}_2$ and ${\xi}_1, {\eta}_1 \in{\hat{R}}_1$, we obtain by contraction ${\eta}^k=f^k_i{\xi}^i, {\eta}^k_k={\xi}^r_r + na_i{\xi}^i$. The inverse transformation allowing to describe 
${\hat{R}}^*_1$ is $ {\xi}^i=g^i_k{\eta}^k, {\xi}^r_r={\eta}^k_k - nb_k{\eta}^k$ if we set $a_i=f^k_ib_k$ as in ([23], p 448). Hence, according to the variational formula $\delta {\chi}_1=f^{-1}_2 \circ D{\eta}_2 \circ j_1(f_1)$ of the first lemma of this section, if the dual finite field of $({\chi}^k_{,i},{\chi}^r_{r,i})$ {\it over the source} is $({\cal{X}}^{,i}_k, {\cal{X}}^i)$, we obtain for the dual field $({\cal{Y}}^{,k}_l, {\cal{Y}}^k)$ {\it over the target} the dual formulas replacing the ones used for the adjoint operator:  \\
\[  \Delta {\cal{Y}}^{,k}_l= g^s_l {\cal{X}}^{,i}_s{\partial}_if^k - nb_l{\cal{X}}^i{\partial}_if^k, \hspace{3mm} \Delta {\cal{Y}}^k={\cal{X}}^i{\partial}_if^k   \]
This result can also be found by a direct computation showing that:\\
\[ \delta {\chi}^r_{r,i} = [ (\frac{\partial {\eta}^r_r}{\partial y^k} - {\eta}^r_{rk}) - n b_l (\frac{\partial{\eta}^l}{\partial y^k} - {\eta}^l_k)]\frac{\partial f^k}{\partial x^i} \]
{\it The action of the second order jets only is} $ {\cal{Y}}^{,k}_l \longrightarrow {\cal{Y}}^{,k}_l - na_i {\cal{Y}}^k$ if we set $f(x)=x,f^k_i={\delta}^k_i$. The corollary follows when $k=4$ because ${\cal{Y}}^4={\cal{J}}^4$ is the charge density.  \\
\hspace*{12cm}  Q.E.D.   \\

\noindent
{\bf REMARK 3.23}: As we have already seen, the Helmholtz analogy establishes a parallel between analytical mechanics and thermodynamics by setting $T=\dot{q}$ and $L= - F$, that is by trying to describe the absolute temperature as a classical " {\it field} " that should be obtained from the derivative of a "{\it potential} ", in order to bring the possibility of an integration by part. Accordingly, it should not be just a scalar quantity but should transform like a time derivative and we recognize the definition of a jet. On the other side, it follows from the preceding results that the Poincar\'{e} group is of codimension $1$ in its normalizer which is the Weyl group, the "{\it difference} " being the unique space-time dilatation that must be added and which is not accessible to intuition as it must be {\it at the same time} a space dilatation {\it and} a time dilatation, that is a concept where space cannot be separated from time, contrary to the standard examples of thermostatics that we have provided. Also, we must not forget that the measure of $V$ has to do with the translation of a piston ring while the measure of $T$ has similarly to do with the translation of a liquid or a solid in a thermometer. The link with thermodynamics will be obtained by defining as usual the so-called " {\it dimensionless speed} " $v^k/c={\partial}_4f^k$ and the " {\it normalized speed} " $u^k=(v^k/c)/\sqrt{1-\frac{v^2}{c^2}}$ over the target with ${\omega}_{kl}(y)u^ku^l=-1 \Rightarrow u^k\frac{\partial}{\partial y^l}({\omega}_{rk}(y)u^r)=0$. Also, introducing $A_{ij}={\omega}_{rs}(x)A^r_iA^s_j$ while setting ${\Theta}^2(f(x))=1/a(x)$ where $\Theta$ is a gauging of the dilatation subgroup over the target as in section $2$, we may introduce:  \\
\[  \left\{  \begin{array}{l}
\tilde{\theta} = \sqrt{ -A_{44}}=\Theta \sqrt{1-\frac{v^2}{c^2}}{\partial}_4f^4 = \Theta \theta  \\
\tilde{\rho} = \sqrt{-A_{44}}/det(A)=(1/{\Theta}^3)\sqrt{1-\frac{v^2}{c^2}}\partial(x^1,x^2,x^3)/\partial(y^1,y^2,y^3)=\rho/{\Theta}^3
\end{array}
\right.  \]
It follows that $\rho u^k=\frac{1}{\Delta}{\partial}_4f^k$ and we obtain therefore the {\it identity} $\partial (\rho u^k)/\partial y^k \equiv 0$ over the target from the more general {\it identity} $\partial(\frac{1}{\Delta} {\partial}_if^k)/\partial y^k \equiv 0$ over the target ([25], Lemma 3.94, p 490). We have thus a reason for introducing the {\it constraints} $A^1_4=A^2_4=a^3_4=0$ on the moving frame like in ([7], \S 28, \S 45), obtaining the so-called {\it temperature vector} $f^k_4(x)=\frac{u^k(f(x))}{\Theta(f(x))}$ by gauging the Lorentz subgroup. For a fluid at rest, we have $\tilde{\theta}=\Theta{\partial}_4f^4$ and {\it we cannot separate} $\Theta$ {\it from} ${\partial}_4f^4$ in the nonlinear description of the Spencer operator because ${\chi}_0=A-id$, explaining therefore the Helmholtz analogy because the {\it evaluation} $y^4=x^4$ on time does not commute with the integration by part on time. Finally, introducing a density of action $w(A)$ over the source, we get for the variation:  \\
\[ \delta W=\delta \int w(A)dx=\int \frac{\partial w}{\partial A^r_i}\delta A^r_i dx=\int {\cal{X}}^i_rg^r_k(\frac{\partial {\eta}^k}{\partial y^l}-{\eta}^k_l){\partial}_if^ldx=\int {\cal{Y}}^l_k(\frac{\partial {\eta}^k}{\partial y^l}-{\eta}^k_l)dy  \]
where ${\cal{Y}}^l_k=\frac{1}{\delta}g^r_k\frac{\partial w}{\partial A^r_i}{\partial}_if^l$ over the target. With a density $w(\tilde{\rho},\tilde{\theta})=\bar{w}(\rho,\theta,\Theta)$, we get:  \\
\[ \pi=\frac{1}{\Delta}\tilde{\rho}\frac{\partial w}{\partial \tilde{\rho}}, \hspace{2mm} \sigma=\frac{1}{\Delta}\tilde{\theta}\frac{\partial w}{\partial \tilde{\theta}} \hspace{2mm} \Rightarrow \hspace{2mm}  {\cal{Y}}^l_k=-\pi({\delta}^l_k+{\omega}_{rk}u^ru^l)-\sigma{\omega}_{rk}u^ru^l  \]
For example, we may set $- w=\frac{1}{3}\alpha \hspace{1mm}det(A)=\frac{1}{3}\alpha (\tilde{\theta}/\tilde{\rho})=\frac{1}{3}\alpha {\Theta}^4\Delta$ for the black body and recover the pressure $\pi= \frac{1}{3}\alpha {\Theta}^4$.  \\
Finally, using a $R_1$-connection $(\delta, -\gamma)$ with $\gamma=0$ when $n=4$, we may introduce ${\xi}_1=({\xi}=(0,0,0,1), {\xi}^k_i=0)\in R_1\subset {\hat{R}}_1\subset J_1(T)$ over the source and get ${\eta}^k(f(x))={\partial}_4f^k(x), {\eta}^k_u(f(x))f^u_i(x)={\partial}_4f^k_i(x)$ over the target. The following tricky proposition revisits the first principle of relativistic thermodynamics in this new framework while taking into account the Helmholtz analogy:  \\

\noindent
{\bf PROPOSITION 3.24}: We have the formula:  \\
\[   \frac{\partial {\cal{Y}}^l_k}{\partial y^l}{\eta}^k + {\cal{Y}}^l_k{\eta}^k_l=\rho u^l\frac{\partial}{\partial y^l}(\tilde{\theta}\frac{\partial w}{\partial \tilde{\theta}} - w)   \]

\noindent
{\bf Proof}: First of all, as ${\eta}^k=\theta u^k$, $\theta =\rho \Delta$ and $u^k({\delta}^l_k + {\omega}_{rk}u^ru^l)=0$, we have:  \\
\[ \begin{array}{rcl}
\theta u^k \frac{\partial}{\partial y^l}(\pi ({\delta}^l_k+{\omega}_{rk}u^ru^l)) & = & -\pi({\delta}^l_k+{\omega}_{rk}u^ku^l)(u^k\frac{\partial \theta}{\partial y^l} + \theta \frac{\partial u^k}{\partial y^l})\\
   &  =  &  - \pi({\delta}^l_k+{\omega}_{rk}u^ru^l)\theta \frac{\partial u^k}{\partial y^l}  \\
      &   =   & - \pi \theta \frac{\partial u^l}{\partial y^l}  \\
        &  =  & - \tilde{\rho}\frac{\partial w}{\partial \tilde{\rho}}\rho \frac{\partial u^l}{\partial y^l} \\
           &  = &   \frac{\partial \bar{w}}{\partial \rho}\rho u^l\frac{\partial \rho}{\partial y^l}
           \end{array}   \]
Similarly, we have:   \\
\[  \begin{array}{rcl}
\theta u^k \frac{\partial}{\partial y^l} (\sigma {\omega}_{rk}u^ru^l)  & = &  - \theta \rho u^l\frac{\partial}{\partial y^l}(\frac{1}{\theta}\tilde{\theta}\frac{\partial w}{\partial  \tilde{\theta}} )  \\
  &  =  &  -\rho u^l\frac{\partial}{\partial y^l}(\tilde{\theta}\frac{\partial w}{\partial \tilde{\theta}})+\frac{\partial \bar{w}}{\partial \theta}\rho u^l\frac{\partial \theta}{\partial y^l}
\end{array}   \]
Finally, multiplying the conformal Killing equations for ${\eta}_1$ over the source by $u^ru^l$, we obtain ${\omega}_{rk}u^ru^l{\eta}^k_l= - \frac{1}{n}{\eta}^k_k=(\theta/\Theta) u^l \frac{\partial \Theta}{\partial y^l}$ and we have thus for $n=4$:  \\
\[  \Theta \frac{\partial \bar{w}}{\partial \Theta}= - 3\tilde{\rho}\frac{\partial w}{\partial \tilde{\rho}} + \tilde{\theta}\frac{\partial w}{\partial \tilde{\theta}} = \Delta (-3\pi +\sigma)\hspace{2mm} \Rightarrow \hspace{2mm}  {\cal{Y}}^l_k{\eta}^k_l=\frac{1}{n}( - 3\pi + \sigma){\eta}^k_k= -  \frac{\partial \bar{w}}{\partial \Theta} \rho u^l\frac{\partial \Theta}{\partial y^l}               \]
The proposition follows from the fact that $d\bar{w}=\frac{\partial\bar{w}}{\partial \rho} d\rho + \frac{\partial \bar{w}}{\partial \theta} d\theta +\frac{\partial\bar{w}}{\partial \Theta}d\Theta$ along the trajectory. The extension to continuum mechanics could be done by using the $6$ quantities $A^*_{ij}=A_{ij} - \frac{A_{i4}A_{j4}}{A_{44}}$ for $i,j=1,2,3$ which do not contain $A^4_i$ for $i=1,2,3,4$, in place of $\tilde{\rho}$ and we have $det(A^*_{ij})=1/{\tilde{\rho}}^2$. (Compare to [17]). \\

\hspace*{12cm}   Q.E.D.    \\

From the infinitesimal point of view, we may use (up to sign) $({\partial}_1{\xi}^1+{\partial}_2{\xi}^2+{\partial}_3{\xi}^3) - ({\xi}^1_1+{\xi}^2_2+{\xi}^3_3)$ in place of $\tilde{\rho}$, ${\partial}_4{\xi}^4 - {\xi}^4_4$ in place of $\tilde{\theta}$ and ${\partial}_i{\xi}^r_r -{\xi}^r_{ri}$ in place of ${\chi}^r_{r,i}$. Identifying the speed with a Lorentz gauging, we may introduce the constraint ${\partial}_4{\xi}^i-{\xi}^i_4=0$ in order to have ${\partial}_{44}{\xi}^i - {\xi}^i_{44}={\partial}_4{\xi}^i_4-{\xi}^i_{44}$ and we may consider by substraction the components ${\partial}_i{\xi}^4_4 - {\partial}_{44}{\xi}^i$ because ${\xi}^4_{4i}={\xi}^i_{44}=A_i$ for $i=1,2,3$, a result explaining why it is possible to combine the gradient of temperatue with the (small) dimensionless object obtained when dividing the acceleration by the square $c^2$ of the speed of light as in ([9], (37), p 922). Also, considering the coupling bilinear term $({\partial}_4{\xi}^i-{\xi}^i_4)({\partial}_i{\xi}^4_4 - {\xi}^4_{4i})$ as in ([24]), we understand how it is possible to introduce a spatial heat flow in the components ${\sigma}^i_4$ of the stress $\sigma$ for $i=1,2,3$, contrary to the phenomenological approach of ([9], (24)+(37)).  \\

\noindent
{\bf EXAMPLE 3.25}: ({\it Projective group}) With $n=1$, let us consider the projective transformations of the real line. For this, we may introduce the third order system $R_3\subset J_3(T)$ of infinitesimal Lie equations defined by ${\xi}_{xxx}=0$ and we let the reader check easily that $[R_3,R_3]\subset R_3$. We may then exhibit a basis made by the three infinitesimal generators $\{ {\theta}_1={\partial}_x, {\theta}_2=x{\partial}_x, {\theta}_3=\frac{1}{2}x^2{\partial}_x \}$ while introducing ${\xi}_{\mu}={\lambda}^{\tau}{\partial}_{\mu}{\theta}_{\tau}$ leading to $\xi={\lambda}^{\tau}{\theta}_{\tau}, {\xi}_x={\lambda}^{\tau}{\partial}_x{\theta}_{\tau}, {\xi}_{xx}={\lambda}^{\tau}{\partial}_{xx}{\theta}_{\tau}$ along with the following corresponding linear transformation:      \\
\[  (L,M,N)\left( \begin{array}{l}
\xi  \\
{\xi}_x  \\
{\xi}_{xx}
\end{array}  \right)=(L,M,N)
\left(  \begin{array}{lcc}
1 & x & \frac{1}{2}x^2 \\
0 & 1 & x \\
0 & 0 & 1
\end{array} \right)
\left( \begin{array}{lll}
{\lambda}^1  \\
{\lambda}^2  \\
{\lambda}^3
\end{array}  \right)        \]
bringing the relations:  \\
\[ X\equiv {\partial}_x\xi-{\xi}_x= ({\partial}_x{\lambda}^{\tau}){\theta}_{\tau}, Y\equiv {\partial}_x{\xi}_x-{\xi}_{xx}=({\partial}_x{\lambda}^{\tau}){\partial}_x{\theta}_{\tau},Z\equiv {\partial}_x{\xi}_{xx}-{\xi}_{xxx}=({\partial}_x{\lambda}^{\tau}){\partial}_{xx}{\theta}_{\tau}   \]
Finally, it just remains to integrate by parts the expression/contraction:  \\
\[  FX+GY+HZ \equiv F({\partial}_x\xi-{\xi}_x)+G({\partial}_x{\xi}_x-{\xi}_{xx})+H({\partial}_x{\xi}_{xx}-{\xi}_{xxx})  \]
while taking into account the fact that ${\xi}_{xxx}=0$ in order to find the dual Cosserat equations:  \\
\[   \left\{ \begin{array}{lclcl}
\xi & \rightarrow &{\partial}_xF   & = & L \\
& & & &   \\
{\xi}_x & \rightarrow & {\partial}_xG + F & = & M  \\
& & & &  \\
{\xi}_{xx} & \rightarrow &{\partial}_xH + G &=& N 
\end{array}  \right. 
\Longleftrightarrow
\left\{  \begin{array}{lcl}
 {\partial}_x(F) & = &L    \\
  &  &  \\
  {\partial}_x(G+xF)  & =& M + xL \\
&  &  \\
{\partial}_x(H + xG+\frac{1}{2}x^2F) & = & N + xM + \frac{1}{2}x^2L 
\end{array}
\right.     \]
involving the formal adjoint of the first Spencer operator and the above linear transformation acting on both sides of the equations. The study of the conformal group is quite similar to that of the projective group because the symbol at order $3$ is equal to zero in both cases. We may therefore just replace ${\partial}_r{\xi}^r-{\xi}^r_r, {\partial}_i{\xi}^r_r-{\xi}^r_{ri},{\partial}_i{\xi}^r_{rj}-0$ by ${\partial}_x\xi - {\xi}_x, {\partial}_x{\xi}_x - {\xi}_{xx}, {\partial}_x{\xi}_{xx}-0$ respectively.\\

\noindent
{\bf 4)  CONCLUSION}  \\

Considering a Lie group of transformations as a Lie pseudogroup of transformations, we have revisited in this new framework the mathematical foundations of both thermodynamics and gauge theory. As a byproduct, we have proved that the methods known for Lie groups cannot be adapted to Lie pseudogroups and that the two approaches are thus not compatible on the purely mathematical level. In particular, the electromagnetic field, which is a $2$-form with value in the Lie algebra of the unitary group $U(1)$ according to classical gauge theory, becomes part of a $1$-form with value in a Lie algebroid in the new conformal approach. More generally, {\it shifting by one step the interpretation of the differential sequences involved},  the "field" is no longer a $2$-form with value in a Lie algebra but must be a $1$-form with value in a Lie algebroid. Meanwhile, we have proved that the use of Lie equations allows to avoid any explicit description of the action of the underlying group, a fact particularly useful for the elations of the conformal group. However, a main problem is that the formal methods developped by Spencer and coworkers around 1970 are still not acknowledged by physicists and we don't even speak about the Vessiot structure equations for pseudogroups, not even acknowledged by mathematicians after more than a century. Finally, as a very striking fact with deep roots in homological algebra, the Clausius/Cosserat /Maxwell/Weyl equations can be parametrized, contrary to Einstein equations. We hope this paper will open new trends for future theoretical physics, based on the use of new differential geometric methods.  \\

\noindent
{\bf REFERENCES}\\

\noindent
[1] V. ARNOLD: M\'{e}thodes Math\'{e}matiques de la M\'{e}canique Classique, Appendice 2 (G\'{e}od\'{e}siques des m\'{e}triques invariantes \`{a} gauche sur des groupes de Lie et hydrodynamique des fluides parfaits), MIR, moscow, 1974,1976. \\
\noindent
[2] G.  BIRKHOFF: Hydrodynamics, Princeton University Press, 1954.  \\
\noindent
[3] S. BORDONI: Routes Towards an Abstract Thermodynamics in the Late Nineteeth Century, Eur. Phys. J. H, 38 (2013) 617-660.  \\
http://dx.doi.org/10.1140/epjh/e2013-40028-7   \\
\noindent
[4] M. BORN: Physik Z., 22 (1921) 218, 249, 282.  \\
\noindent
[5] C. CARATHEODORY: Math. Annalen, 67 (1909) 355. \\
\noindent
[6] C. CHEVALLEY, S. EILENBERG: Cohomology Theory of Lie Groups and Lie Algebras, Trans. American Math. Society, 63, 1 (1948) 85-124. \\
\noindent
[7] E. COSSERAT, F. COSSERAT: Th\'{e}orie des Corps D\'{e}formables, Hermann, Paris, 1909.\\
\noindent
[8] P. DUHEM: Commentaires aux Principes de la Thermodynamique, Journal de Math\'{e}matiques Pures et AppliquŽes, Part I, 8 (1892) 269-330; Part II, 9(1893) 293-359; Part III, 10 (1894) 207-285.  \\
\noindent
[9] C. ECKART: The Thermodynamics of Irreversible Processes III, Physical review, 58, (1940) 919-924.\\
\noindent
[10] L. P. EISENHART: Riemannian Geometry, Princeton University Press, Princeton, 1926.  \\
\noindent
[11] E. INONU, E.P. WIGNER: On the Contraction of Lie Groups and Lie Algebras, Proc. Nat. Acad. Sci. USA, 39 (1953) 510.  \\
\noindent
[12] M. JANET: Sur les Syst\`{e}mes aux D\'{e}riv\'{e}es Partielles, Journal de Math., 8 (1920) 65-151. \\
\noindent 
[13] R. KUBO: Thermodynamics, North Holland, Amsterdam, 1968.  \\
\noindent
[14] A. KUMPERA, D.C. SPENCER: Lie Equations, Ann. Math. Studies 73, Princeton University Press, Princeton, 1972.\\
\noindent
[15] E. KUNZ: Introduction to Commutative Algebra and Algebraic Geometry, BirkhaŸser, 1985.  \\
\noindent
[16]  F.S. MACAULAY: The Algebraic Theory of Modular Systems, Cambridge, 1916.  \\
\noindent
[17] G. NORDSTR\"{O}M: Einstein's Theory of Gravitation and Herglotz's Mechanics of Continua, Proc. Kon. Ned. Akad. Wet., 19 (1917) 884-891. \\
\noindent
[18] V. OUGAROV: Th\'{e}orie de la Relativit\'{e} Restreinte, MIR, Moscow, 1969 ( french, 1979).\\
\noindent
[19] H. POINCARE: Sur une Forme Nouvelle des Equations de la M\'{e}canique, C. R. Acad\'{e}mie des Sciences Paris, 132 (7) (1901) 369-371.  \\
\noindent
[20] J.-F. POMMARET: Systems of Partial Differential Equations and Lie Pseudogroups, Gordon and Breach, New York, 1978; Russian translation: MIR, Moscow, 1983.\\
\noindent
[21] J.-F. POMMARET: Differential Galois Theory, Gordon and Breach, New York, 1983.\\
\noindent
[22] J.-F. POMMARET: Lie Pseudogroups and Mechanics, Gordon and Breach, New York, 1988.\\
\noindent
[23] J.-F. POMMARET: Partial Differential Equations and Group Theory, Kluwer, 1994.\\
http://dx.doi.org/10.1007/978-94-017-2539-2    \\
\noindent
[24] J.-F. POMMARET: Group Interpretation of Coupling Phenomena, Acta Mechanica, 149 (2001) 23-39.\\
http://dx.doi.org/10.1007/BF01261661  \\
\noindent
[25] J.-F. POMMARET: Partial Differential Control Theory, Kluwer, Dordrecht, 2001.\\
\noindent
[26] J.-F. POMMARET: {\it Algebraic Analysis of Control Systems Defined by Partial Differential Equations}, Advanced Topics in Control Systems Theory, Springer, Lecture Notes in Control and Information Sciences 311 (2005) Chapter 5, pp. 155-223.\\
\noindent
[27] J.-F. POMMARET: Parametrization of Cosserat Equations, Acta Mechanica, 215 (2010) 43-55.\\
http://dx.doi.org/10.1007/s00707-010-0292-y  \\
\noindent
[28] J.-F. POMMARET: Spencer Operator and Applications: From Continuum Mechanics to Mathematical Physics, in "Continuum Mechanics-Progress in Fundamentals and Engineering Applications", Dr. Yong Gan (Ed.), ISBN: 978-953-51-0447--6, InTech, 2012, Available from: \\
http://www.intechopen.com/books/continuum-mechanics-progress-in-fundamentals-and-engineering-applications/spencer-operator-and-applications-from-continuum-mechanics-to-mathematical-physics  \\
\noindent
[29] J.-F. POMMARET: Deformation Cohomology of Algebraic and Geometric Structures, 2012, Preprint.  \\
http://arxiv.org/abs/1207.1964  \\
\noindent
[30] J.-F. POMMARET: The Mathematical Foundations of General Relativity Revisited, Journal of Modern Physics, 4 (2013) 223-239. \\
 http://dx.doi.org/10.4236/jmp.2013.48A022   \\
 \noindent
[31] J.-F. POMMARET: The Mathematical Foundations of Gauge Theory Revisited, Journal of Modern Physics, 5 (2014) 157-170.  \\
http://dx.doi.org/10.4236/jmp.2014.55026  \\
\noindent
[32] J. J. ROTMAN: An Introduction to Homological Algebra, Academic Press, 1979.\\
\noindent
[33] D. C. SPENCER: Overdetermined Systems of Partial Differential Equations, Bull. Am. Math. Soc., 75 (1965) 1-114.\\
\noindent
[34] P.P. TEODORESCU: Dynamics of Linear Elastic Bodies, Editura Academiei, Bucuresti, Romania; Abacus Press, Tunbridge, Wells, 1975.\\
\noindent
[35] E. VESSIOT: Sur la Th\'{e}orie des Groupes Infinis, Ann. Ec. Norm. Sup., 20 (1903) 411-451.\\
\noindent
[36] H. WEYL: Space, Time, Matter, Springer, 1918, 1958; Dover, 1952. \\
\noindent

\end{document}